\documentclass[12pt,english]{smfart}
\usepackage{fullpage, amsmath, amsthm,amsfonts,amssymb,stmaryrd, mathrsfs, mathtools}
\usepackage{hyperref}
\usepackage{mathdots}
\usepackage[pdftex]{color}
\usepackage[all]{xy}
\usepackage{xcolor}
\usepackage{hyperref}
\usepackage{bbm}
\usepackage{xspace}
\usepackage[T1]{fontenc}

\usepackage[autolanguage]{numprint}
\usepackage{relsize}
\usepackage[bbgreekl]{mathbbol}
\usepackage{tikz-cd}
\usepackage[lite,abbrev,msc-links,backrefs]{amsrefs}

\numberwithin{equation}{subsection}

\newtheorem{theorem}[subsection]{Theorem}

\newtheorem{corollary}[subsection]{Corollary}
\newtheorem{lemma}[subsection]{Lemma}
\newtheorem{proposition}[subsection]{Proposition}

\theoremstyle{definition}

\newtheorem{definition}[subsection]{Definition}

\newtheorem{remark}[subsection]{Remark}

\renewcommand{\proofname}{Proof}

\DeclareSymbolFontAlphabet{\mathbbl}{bbold}
\newcommand{\prism}{{\mathlarger{\mathbbl{\Delta}}}}

\def\AAA{\mathbb{A}}

\def\NN{\mathbb{N}}

\def\QQ{\mathbb{Q}}
\def\RR{\mathbb{R}}

\def\ZZ{\mathbb{Z}}

\def\scrC{\mathscr{C}}

\def\gothS{\mathfrak{S}}
\def\gothm{\mathfrak{m}}

\newcommand{\rd}{\mathrm{d}}

\newcommand{\rA}{\mathrm{A}}

\newcommand{\rR}{\mathrm{R}}

\newcommand{\qproet}{\mathrm{qproet}}

\def\calC{\mathcal{C}}

\def\calE{\mathcal{E}}
\def\calF{\mathcal{F}}

\def\calI{\mathcal{I}}
\def\calJ{\mathcal{J}}
\def\calK{\mathcal{K}}
\def\calL{\mathcal{L}}
\def\calM{\mathcal{M}}

\def\calO{\mathcal{O}}

\DeclareMathOperator{\End}{End}
\DeclareMathOperator{\Gal}{Gal}

\DeclareMathOperator{\coker}{coker}

\DeclareMathOperator{\Spa}{Spa}
\DeclareMathOperator{\Ainf}{\mathrm{A}_{\mathrm{inf}}}

\DeclareMathOperator{\Spec}{Spec}
\DeclareMathOperator{\Spf}{Spf}
\DeclareMathOperator{\Mod}{\mathbf{Mod}}

\newcommand{\cO}{\mathcal{O}}

\newcommand{\colim}{\mathrm{colim}}
\newcommand{\CR}{\mathbf{CR}}
\newcommand{\CA}{\mathrm{\check{C}A}}
\newcommand{\ad}{\mathrm{ad}}
\newcommand{\qOmega}{\mathrm{q}\Omega}
\newcommand{\rel}{\mathrm{rel}}

\newcommand{\hra}{\hookrightarrow}
\newcommand{\ra}{\rightarrow}

\newcommand{\xra}{\xrightarrow}

\newcommand{\et}{\mathrm{\acute{e}t}}

\newcommand{\cris}{\mathrm{cris}}

\newcommand{\id}{\mathrm{id}}

\newcommand{\Kos}{\mathrm{Kos}}

\newcommand{\Strat}{\mathbf{Strat}}
\newcommand{\ev}{\mathrm{ev}}

\newcommand{\Vect}{\mathrm{Vect}}
\newcommand{\perf}{\mathrm{perf}}
\newcommand{\Shv}{\mathrm{Shv}}

\newcommand{\opp}{\mathrm{opp}}

\newcommand{\Loc}{\mathrm{Loc}}

\newcommand{\proet}{\mathrm{proet}}
\newcommand{\anfr}{\mathrm{an-fr}}
\newcommand{\an}{\mathrm{an}}

\newcommand{\laurent}{\cO_{\prism}[1/\calI_{\prism}]^{\wedge}_p}

\newcommand{\MIC}{\mathbf{MIC}}

\newcommand{\tf}{\mathrm{tf}}
\newcommand{\gp}{\mathrm{gp}}
\newcommand{\tAinf}{\tilde{\mathrm{A}}_{\inf}}

\newcommand{\CRhat}{\widehat{\CR}}
\newcommand{\qHiggs}{\mathbf{qHiggs}}
\newcommand{\Higgs}{\mathbf{Higgs}}
\newcommand{\uA}{\underline{A}}

\newcommand{\uR}{\underline {R}}
\newcommand{\uX}{\underline {X}}
\newcommand{\ufS}{\underline{\gothS}}
\newcommand{\urA}{\underline{\rA}}
\newcommand{\coh}{\mathrm{coh}}
\newcommand{\uO}{\underline{\cO}}
\newcommand{\topnil}{\mathrm{nil}}
\newcommand{\st}{\mathrm{st}}
\newcommand{\fX}{\mathfrak{X}}

\newcommand{\ufX}{\underline{\fX}}
\begin{document}

	\title{A prismatic-etale comparison theorem in the semi-stable case}
	
	\author{Yichao Tian}
	\address{Morningside Center of Mathematics,  State Key Laboratory of Mathematical Sciences\\ AMSS, Chinese Academy of Sciences \\
			55 Zhong Guan Cun East Road\\
			100190, Beijing, China}
	\date{}

	\begin{abstract}
	Let $K|\QQ_p$ be a complete discrete valuation field with perfect residue field, $\cO_K$ be its ring of integers. Consider a semistable $p$-adic formal scheme $X$ over $\Spf(\cO_K)$ with smooth generic fiber $X_{\eta}$. Du--Liu--Moon--Shimizu showed recently that the category of analytic  prismatic $F$-crystals on the absolute log prismatic site  $\uX_{\prism}$ is equivalent to the category of semistable \'etale $\ZZ_p$-local systems on the adic generic fiber $X_{\eta}$. In this article, we prove a comparison between the Breuil--Kisin cohomology of an analytic log prismatic $F$-crystal on $\uX_{\prism}$ and the \'etale cohomology of its corresponding \'etale $\ZZ_p$-local system. This generalizes Guo--Reneicke's prismatic--\'etale comparison for crystalline $\ZZ_p$-local systems to the semi-stable case. 
		
	\end{abstract}
	\maketitle	
	\tableofcontents
	\setcounter{section}{-1}
	
	\section{Introduction}
	
\subsection{Crystalline comparison theorems via prismatic cohomology}	In their ground-breaking work \cite{BS},   Bhatt and Scholze introduced the prismatic cohomology theory for $p$-adic formal schemes, based on their joint work with Morrows \cite{BMS1} and \cite{BMS2}. The prismatic theory provides a uniform framework to study other $p$-adic cohomology theories such as de Rham cohomology, crystalline cohomology and $p$-adic \'etale cohomology, and  it suggests a new approach to $p$-adic comparison theorems, which is one of  central topics in $p$-adic Hodge theory. 

	 Recently, using the prismatic cohomology,  Guo and Reneicke \cite{GR} established a crystalline comparison theorem for cyrstalline  $\ZZ_p$-local systems.
	 More precisely,  let $K$ be a  complete discrete valuation field over $\QQ_p$ with a perfect residue field $k$,  and let $\cO_K$ be its ring of integers. Consider   a  proper smooth $p$-adic formal scheme $X$ over $\Spf(\cO_K)$ with generic  adic fiber $X_{\eta}$.
	  Guo--Reneicke \cite{GR} and  Du--Liu--Moon--Shimizu \cite{DLMS1} establish independently an equivalence between the category of crystalline \'etale $\ZZ_p$-local systems on $X_{\eta}$ and that of analytic prismatic $F$-crytsals on the absolute prismatic site $X_{\prism}$. When  $X=\Spf(\cO_K)$, crystalline \'etale $\ZZ_p$-local systems on $X_{\eta}$  are nothing but Galois stable $\ZZ_p$-lattices in classical crystalline representations of $\Gal(\overline K/K)$, and this result was a theorem of Bhatt--Scholze \cite{BS2}. 
	 Once such an equivalence is established,  Guo--Reneicke's approach to a  comparison theorem for crystalline $\ZZ_p$-local systems  can be divided into two parts: the first part is  a prismatic-\'etale comparison theorem which relates the cohomology of an analytic prismatic $F$-crystal on $X_{\prism}$ to the  geometric \'etale cohomology of the corresponding crystalline local system on $X_{\eta}$; and the second part is a prismatic-crystalline comparison that compares the cohomology of an analytic prismatic $F$-crystal with the crystalline cohomology of some crystalline prismatic crystals. 
	    In this article, we try to generalize Guo--Reneicke's prismatic-\'etale comparison theorem for semi-stable local systems on the adic generic fiber of  $p$-adic formal schemes over $\Spf(\cO_K)$. 
	 
	\subsection{Semistable local systems and analytic prismatic $F$-crystals}
	 Let $X$ be a separated semi-stable $p$-adic formal scheme over $\Spf(\cO_K)$  (Def.~\ref{D:semi-stable}), i.e. $X$ is \'etale locally isomorphic to $$\Spf(\cO_K\langle T_0, \cdots, T_r, T_{r+1}^{\pm1}, \cdots, T_d^{\pm1}\rangle/(T_0\cdots T_r-\pi),$$ where $\pi$ is a fixed uniformizer of $\cO_K$. We equip $X$ with the canonical log structure $\calM_X=\cO_{X_{\et}}\cap \cO_{X_{\et}}[1/p]^{\times}$, and let $\uX:=(X, \calM_X)$ denote the resulting log formal scheme. 
	 
	 On the adic generic fiber $X_{\eta}$, there is a natural notion of semistable \'etale $\ZZ_p$-local systems that goes back to the work of Faltings \cite{Fal} and Andreatta--Iovita \cite{AI}:  they are  those $\ZZ_p$-local systems on $X_{\eta}$ associated to some   $F$-crystals on the absolute log crystalline site of $\uX_1:=\uX\otimes_{\cO_K}\cO_K/p$ in the sense of \cite[Def. 3.40]{DLMS}. 
	     When $X=\Spf(\cO_K)$,  this definition gives  exactly Galois stable $\ZZ_p$-lattices in Fontaine's semi-stable representations of $\Gal(\overline K/K)$. We denote by  $\Loc^{\st}_{\ZZ_p}(X_{\eta, \et})$  the category of semistable \'etale $\ZZ_p$-local systems on $X_{\eta}$. 
	 	     
	  On the other hand, Koshikawa \cite{Kos} and Koshikawa--Yao \cite{KY}  generalized the prismatic theory to the logarithmic setting, and defined the  absolute log prismatic site $\uX_{\prism}$. The objects of $\uX_{\prism}$ consist of log prisms $(A,I, M_{\Spf(A)})$  together with a strict morphism of log formal schemes $(\Spf(A/I), M_{\Spf(A/I)})\to \uX$, and we equip $\uX_{\prism}$ with the strict flat topology (See Def.~\ref{D:abs-prismatic-site}). There is a natural  structure ring sheaf $\cO_{\prism}$
on $\uX_{\prism}$ equipped with an Frobenius endomorphism $\phi_{\cO_{\prism}}$. 	 

	Simiarly to the smooth case, Du--Liu--Moon--Shimizu \cite[Def. 3.3]{DLMS} defined the category of analytic  prismatic  $F$-crystals on  $\uX_{\prism}$ as 
	   \[\Vect^{\an}(\uX_{\prism}, \cO_{\prism})^{\phi=1}:=\varprojlim_{(A,I, M_{\Spf(A)})\in \uX_{\prism}} \Vect(\Spec(A)\backslash V(p, I))^{\phi=1},
	   \] 
	   where $\Vect(\Spec(A)\backslash V(p, I))^{\phi=1}$ means the category of vector bundles $\calE_A$ on $\Spec(A)\backslash V(p, I)$ together with an isomorphism $\phi_A^*\calE_A[\frac{1}{I}]\simeq \calE_A[\frac{1}{I}]$. 
	   They prove \cite[Cor. 5.2]{DLMS} that there exists an equivalence of categories  
	 \[
	 T\colon \Vect^{\an}(\uX_{\prism}, \cO_{\prism})^{\phi=1}\xra{\sim} \Loc^{\st}_{\ZZ_p}(X_{\eta, \et}).
	 \]

 Our main result in this article is then a comparison theorem between the cohomology of an analytic prismatic $F$-crystal on $\uX_{\prism}$ and  the cohomology of its attached semistable $\ZZ_p$-local system.

\subsection{Breuil--Kisin cohomology of analytic prismatic $F$-crystals} Let $\gothS=W(k)[[u]]$ and $E(u)\in \gothS$ be the Eisenstein polynomial of $\pi$ over $W(k)[1/p]$. Equip $\gothS$ with the $\delta$-structure with $\delta(u)=0$  and the $\delta_{\log}$-structure $\alpha, \delta_{\log}\colon \NN\to \gothS$ with $\alpha(1)=u$ and $\delta_{\log}=0$. Then we have  the absolute Breuil--Kisin log prism $\ufS=(\gothS,  (E(u)), \NN)^a$ (cf. \S~\ref{S:BK-log-prism}), and the relative log prismatic site $(\uX/\ufS)_{\prism}$ (cf. \S~\ref{S:relative-prismatic-site}), which is the subcategory of $\uX_{\prism}$ consisting of log prisms over $\ufS$.
 
First, we remark that every analytic prismatic $F$-crystal $(\calE, \phi_{\calE})$ on $\uX_{\prism}$ has a canonical extension $(j_*\calE, \phi_{\calE})$ to a complete prismatic crystal over $\uX_{\prism}$ (Prop.~\ref{P:can-extension-analytic}). We define the Breuil--Kisin cohomology of $(\calE, \phi_{\calE})$ as  
$$R\Gamma_{\gothS}(X, \calE):=R\Gamma((\uX/\ufS)_{\prism}, j_*\calE).$$
This is a derived $(p, E(u))$-complete object in the derived category  of $\gothS$-modules $D(\gothS)$, and it has with a natural Frobenius action induced from $\phi_{\calE}$. If $X$ is proper over $\Spf(\cO_K)$, then each cohomology group $H^i_{\gothS}(X, \calE)$ is a finitely generated $\gothS$-module and vanishes if $i\notin [0, 2d]$ with $d$ being the relative dimension of $X$ over $\cO_K$ (Thm. \ref{T:BK-analytic-cohomology}).

Let $C=\widehat{\overline K}$ be the $p$-adic completion of an algebraic closure of $K$, $\cO_C$ be its ring of integers with tilt $\cO_C^{\flat}$, and  $A_{\inf}=W(\cO_C^{\flat})$. Fix a compatible system of primitive $p^n$-th root of unity $(\zeta_{p^n})_{n\geqslant 1}$ in $C$, and a compatible system $(\pi^{1/p^n})_{n\geqslant 1}$ of $p^n$-th root of $\pi$. Put $\epsilon=(1,\zeta_p, \zeta_{p^2}, \cdots)\in \cO_C^{\flat}$,  $\mu:=[\epsilon]-1$, and $\pi^{\flat}=(\pi, \pi^{1/p}, \cdots)\in \cO_C^{\flat}$. We have then a $\phi$-equivariant embedding $\iota\colon \gothS\to A_{\inf}$ given by $u\mapsto [\pi^{\flat}]^p$.

\begin{theorem}\label{T:comp-prismatic-etale-analytic}
	Let $X$ be a separated semistable $p$-adic formal scheme over $\Spf(\cO_K)$, and let $X_C:=X_{\eta}\times_{\Spa(K, \cO_K)}\Spa(C,\cO_C)$ denote the geometric generic fiber. 
		Let $(\calE, \phi_{\calE})$ be an object of  $\Vect^{\an}(\uX_{\prism},\cO_{\prism})^{\phi=1}$, and $T(\calE)=T(\calE, \phi_{\calE})$ be the corresponding semistable local system on $X_{\eta}$.
	\begin{enumerate}
		\item[(1)]There exists a canonical isomorphism 
		\[
		(R\Gamma_{\gothS}(X,  \calE)\otimes^L_{\gothS}A_{\inf})^{\wedge}[\frac{1}{\mu}]\simeq R\Gamma(X_{C, \proet}, T(\calE)\otimes_{\ZZ_p} \AAA_{\inf,X_C})[\frac{1}{\mu}]
		\]
		equivariant under the natural Frobenius and $\Gal(\overline K/K)$-actions, 
		where  $X_{C, \proet}$ is the pro-\'etale site of $X_C$,  $\AAA_{\inf,X_C}:=W(\widehat{\cO}_{X}^{+, \flat})$ is the  sheaf version of $A_{\inf}$ on $X_{C,\proet}$, and $(-)^{\wedge}$ means the derived $(p, \mu)$-completion. 
		\item[(2)] Suppose that $X$ is proper over $\Spf(\cO_K)$.  Then there exists a canonical isomorphism 
		\[
		R\Gamma_{\gothS}(X, \calE)\otimes_{\gothS}^LA_{\inf}[\frac{1}{\mu}]\simeq R\Gamma(X_{C,\et}, T(\calE))\otimes_{\ZZ_p}^LA_{\inf}[\frac{1}{\mu}] 
		\]
		equivariant under the natural Frobenius and $\Gal(\overline K/K)$-actions on both sides.
		In particular,    one has a canonical isomorphism 
		\[
		H^i_{\gothS}(X, \calE)\otimes_{\gothS}A_{\inf}[\frac{1}{\mu}]\simeq H^i(X_{C,\et}, T(\calE))\otimes_{\ZZ_p}A_{\inf}[\frac{1}{\mu}]
		\]
		for all integers $i\in \ZZ$.
		
	\end{enumerate}

\end{theorem}
If $X$ is smooth and proper over $\Spf(\cO_K)$, then statement (2) of this theorem is a special case of Guo--Reneicke's strong \'etale comparison theorem \cite[Thm 9.1]{GR}. For the constant coefficient $\calE=\cO_{\prism}$, statement (2) was  proved in \cite[Prop. 8.3]{KY} for more general log formal schemes. Statement (1) seems  new even when $X$ is smooth (but non-proper) over $\Spf(\cO_K)$.
Note that  our approach to Theorem~\ref{T:comp-prismatic-etale-analytic} is different from that of  Guo--Reneicke \cite{GR}, even though we still need their \cite[Lemma 3.12]{GR} in our argument.  In fact,  they prove first a weaker version of prismatic--\'etale comparison and deduce the strong version from the Poincar\'e duality for prismatic cohomology, while our approach leads to a more direct proof of the strong \'etale comparison without using the Poincar\'e duality.

 \subsection{Comparison theorems for formal schemes over $\cO_C$} In order to prove Theorem~\ref{T:comp-prismatic-etale-analytic}, we will actually establish a prismatic--\'etale  comparison theorem for some complete prismatic $F$-crystals  on the absolute prismatic site of a semi-stable formal scheme over $\Spf(\cO_C)$. 
 
 Let $\fX$ be a semistable $p$-adic formal scheme over $\Spf(\cO_C)$ (cf. Def. \ref{D:semi-stable}), and $\ufX=(\fX, \calM_{\fX})$ be the associated log formal scheme with the canonical log structure $\calM_{\fX}:=\cO_{\fX_{\et}}\cap \cO_{\fX_{\et}}[1/p]^{\times}$. Even though the log formal scheme $\ufX$ is no longer fs, we still have the absolute log prismatic site $\ufX_{\prism}$ (Def. \ref{D:abs-prismatic-site}). Let  $\tAinf(\uO_C)$ denote the log prism $(A_{\inf}, \ker(\theta_{\cO_C}\circ\phi^{-1}), M_{\cO^{\flat}_C})^a$ defined in \S~\ref{S:setup of notation}, where $\theta_{\cO_C}\colon A_{\inf}\to \cO_C$ is Fontaine's surjection $[x]\mapsto x^{\sharp}$. Then $\tAinf(\uO_C)$  is a final object in the log prismatic site attached to $\Spf(\cO_C)$ equipped with the canonical log structure, and we have a natural equivalence of sites $\ufX_{\prism}\simeq (\ufX_{\prism}/\tAinf(\uO_C))_{\prism}$.
  
 
 We will consider the category $\CRhat^{\anfr}(\ufX_{\prism}, \cO_{\prism})^{\phi=1}$ of complete prismatic $F$-crystals on $\ufX_{\prism}$ whose restriction to the ``analytic locus''  are vector bundles (see Def. \ref{D:complete-prismatic-F-crystal}). This can be viewed as a more practical analogue of the category of  analytic prismatic $F$-crystals in this non-Noetherian  case. 
 If $\fX=X\times_{\Spf(\cO_K)}\Spf(\cO_C)$ for some semi-stable formal scheme $X$ over $\Spf(\cO_K)$, then there exists a natural functor 
 \[
 \Vect^{\an}(\uX_{\prism}, \cO_{\prism})^{\phi=1}\to \CRhat^{\anfr}(\ufX_{\prism}, \cO_{\prism})^{\phi=1}
 \]given by the pullback to $\ufX_{\prism}$ of the canonical extension of analytic prismatic $F$-crystals. Via the equivalence $\ufX_{\prism}\simeq (\ufX_{\prism}/\tAinf(\uO_C))_{\prism}$, we see that, for an object $(\calE, \phi_{\calE})$ of $\CRhat^{\anfr}(\ufX_{\prism}, \cO_{\prism})^{\phi=1}$, $R\Gamma(\ufX_{\prism}, \calE)$ is an object in $D(A_{\inf})$,  the derived category  of $A_{\inf}$-modules.
 
  We have also an \'etale realization functor (see \S~\ref{S:etale-realization-complete})
 \[
 T\colon  \CRhat^{\anfr}(\ufX_{\prism}, \cO_{\prism})^{\phi=1}\to \Loc_{\ZZ_p}(\fX_{\eta, \et}),
 \]
 where the right hand side denotes the category of \'etale $\ZZ_p$-local systems on the adic generic fiber $\fX_{\eta}$. This functor is compatible with the equivalence between analytic prismatic $F$-crystals and semistable local systems in an evident sense.   Theorem~\ref{T:comp-prismatic-etale-analytic} is a consequence of the following theorem together with a result on the base change of prismatic cohomology from $\cO_K$ to $\cO_C$ (Prop. \ref{P:base-change-coh-O_C}).
 
 \begin{theorem}[Theorem~\ref{T:etale-comparison-O_C}]\label{T:etale-O_C-intro}
 Under the above notation, let $(\calE, \phi_{\calE})$ be an object of $\CRhat(\ufX_{\prism}, \cO_{\prism})$, and $T(\calE)$ be the associated local system on $\fX_{\eta, \et}$. 
 \begin{enumerate}
 \item[(1)] There exists a canonical isomorphism in  $D(A_{\inf})$
 \[
 R\Gamma(\ufX_{\prism}, \calE)[\frac{1}{\mu}]\simeq R\Gamma(\fX_{\eta, \proet}, T(\calE)\otimes_{\ZZ_p}\AAA_{\inf, \fX_{\eta}})[\frac{1}{\mu}],
 \]
 compatible with the natural Frobenius action on both sides.
 
 \item[(2)] If $\fX$ is proper over $\Spf(\cO_C)$, then there exists a canonical isomorphism in $D(A_{\inf})$
 \[
 R\Gamma(\ufX_{\prism}, \calE)[\frac{1}{\mu}]\simeq R\Gamma(\fX_{\eta, \et}, T(\calE))\otimes^L_{\ZZ_p}A_{\inf}[\frac{1}{\mu}]
 \]
 equivariant under  the natural Frobenius action on both sides.
 \end{enumerate}
 \end{theorem}

 When $\calE=\cO_{\prism}$,  statement (2) recovers  \cite[Theorem~2.3]{CK}  after identifying    the $A_{\inf}$-cohomology $R\Gamma_{A_{\inf}}(\fX)$ in \emph{loc. cit.}   with our $R\Gamma(\ufX_{\prism},\cO_{\prism})$ (cf. Remark~\ref{R:Ainf-coh}(2)). 

We now give some ideas on the proof of Theorem~\ref{T:etale-O_C-intro}. Let $\ufX^{\perf}_{\prism}$ be the site of the full subcategory of $\ufX_{\prism}$ consisting of log prisms $(A, I, M_{\Spf(A)})$ such that the underlying prism $(A,I)$ is perfect. Surprisingly, there exists an equivalence of sites between $\ufX^{\perf}_{\prism}$ and the non-logarithmic perfect prismatic site $\fX^{\perf}_{\prism}$ (Prop.~\ref{P:equiv-perfect-log-site}, see also \cite[Prop. 2.18]{MW}). 
Let $(\calE^{\perf}, \phi_{\calE^{\perf}})$ be the restriction of  the complete prismatic $F$-crystal $(\calE, \phi_{\calE})$ to $\fX^{\perf}_{\prism}\simeq \ufX^{\perf}_{\prism}$. 
As perfect prisms are equivalent to integral perfectoid rings, it is much easier to relate the perfect prismatic site $\fX^{\perf}_{\prism}$ to Scholze's $v$-site of the diamond attached to $\fX_{\eta}$ and hence to the pro-\'etale site $\fX_{\eta, \proet}$ (cf. ~\eqref{E:commutatiive-diag-topoi}). 
Then a crucial step towards Theorem~\ref{T:etale-O_C-intro} is to show that the canonical map 
\begin{equation}\label{E:isom-perfect-non-perfect}
R\Gamma(\ufX_{\prism}, \calE[1/\mu])\xra{\sim} R\Gamma(\fX^{\perf}_{\prism}, \calE^{\perf}[1/\mu])
\end{equation} 
 is an isomorphism.  

To establish this isomorphism, one needs a detailed local study on the prismatic cohomology $R\Gamma(\ufX_{\prism}, \calE)$. When $\fX$ is smooth over $\Spf(\cO_C)$,  Tsuji \cite{Tsuji} gives a local description of $R\Gamma(\fX_{\prism}, \calE)$ in terms of the de Rham complex of a certain $q$-Higgs module attached to the complete prismatic crystal $\calE$. In this article, we will generalize Tsuji's results to the small semistable case. More precisely, assume that $\fX=\Spf(R)$ is affine such that $R$ admits a framing $\square \colon R^{\square}\to R$ as in \eqref{E:framing-over-C}. Such a framing allows us to construct an explicit object $\urA(R)$ in $\ufX_{\prism}$, which is a  cover of the final object of the topos $\Shv(\ufX)$ (cf. \S~\ref{S:frames}).
 Then to each $\mu$-torsion free complete prismatic crystal $\calE$ on $\fX_{\prism}$, one can attach a topologically quasi-nilpotent $q$-Higgs module $(E, \Theta)$ over  $\rA(R)$ (cf. Prop. ~\ref{P:top-nilpotence}). Moreover, we prove that the prismatic cohomology   $R\Gamma(\ufX_{\prism}, \calE)$ is computed by the de Rham complex of the $q$-Higgs module $(E, \Theta)$ (Thm.~\ref{T:crystal-qHiggs}). 
 Using this result, we can also express $R\Gamma(\ufX_{\prism}, \calE)$ in terms of some Galois cohomology, which allows us finally to prove \eqref{E:isom-perfect-non-perfect} (see Thm.~\ref{T:coh-crystal-Galois-coh}).

\subsection{} The organization of this article is as follows. In Section~\ref{Section1}, we give some preliminaries on the log prismatic  site of a semistable $p$-adic formal scheme over $\Spf(\cO_L)$, where $L$ is either $K$ or $C$. Section~\ref{Section2} is devoted to recalling some basic facts on the $\alpha$-derivation and modules with an integrable connection with respect to such an $\alpha$-derivation. The results of this section will  be only needed in the next two sections.
In Section~\ref{S:local computation}, for an affine small semistable formal scheme $\fX$ over $\Spf(\cO_C)$, we give a local description of the cohomology  a  complete prismatic crystal $\calE$ on $\ufX_{\prism}$ in terms of a certain $q$-Higgs module $E$. 
In Section~\ref{Section4}, we will prove a relation between the cohomology of a complete prismatic crystal on $\ufX_{\prism}$ and some Galois cohomology, and deduce from it the isomorphism \eqref{E:isom-perfect-non-perfect}. In Section~\ref{Section5}, we will  use \eqref{E:isom-perfect-non-perfect} to prove our comparison Theorem \ref{T:etale-O_C-intro} (Theorem~\ref{T:etale-comparison-O_C}). In Section~\ref{S:analytic-cyrstals}, we will come back to the arithmetic situation, i.e. the case of a semistable $p$-adic formal scheme $X$ over $\Spf(\cO_K)$. We will define the Breuil--Kisin cohomology of an analytic prismatic ($F$-)crystal, and prove some finiteness results for it when $X$ is proper. In the last Section~\ref{Section:comparison}, we put all the previous results together to prove Theorem~\ref{T:comp-prismatic-etale-analytic}.

\subsection{Acknowledgement}
I would like to express my gratitude to Ahmed Abbes, Haoyang Guo,  Teruhisa Koshikawa, Tong Liu for  helpful discussions. This work is supported by National Key R\&D Program of China No. 2023YFA1009701, Chinese Natural Science Foundation (No. 12225112, 12231001, 12288201) and CAS Project for Young Scientists in Basic Research (No. YSBR-033).

	\subsection{Notations and Conventions}
	

	\begin{itemize}
		
		\item  All monoids in this article are commutative and unitary. 
		
		\item For a morphism of monoids $u\colon M\to N$, we denote by $u^{\gp}\colon M^{\gp}\to N^{\gp}$ the induced map on the associated abelian groups. 
		\item 	A surjective morphism of monoids  $u:M\to N$ is called \emph{exact} if it induces an isomorphism $M/M^{\times}\xra{\sim}N/N^{\times}$. A surjection of integral monoids $u\colon M\to N$ admits a unique factorization $M\hra \widetilde M:=u^{\gp, -1}(N)\xra{\tilde u} N$, where $\tilde u$ is exact. We call $\tilde u:\widetilde M\to N$ the \emph{exactification} of $u$.
		
		\item 	 \emph{A chart} on a log (formal) scheme $(X,\calM_X)$ is a morphism of monoids $P\to \Gamma(X,\calM_X)$ such that $\underline P^a\xra{\sim} \calM_X$, where $\underline P$ denotes the constant \'etale sheaf and $\underline P^a\to \cO_{X_{\et}}$ is the induced log structure.

			\item 	A log (formal) scheme $(X,\calM_X)$ is called 
		\begin{itemize}
			\item \emph{quasi-coherent}  if it admits a chart \'etale locally, 
			\item \emph{integral} if $\calM_X$ is integral, 
			\item \emph{fs} if $\calM_X$ admits a chart with monoids which are finitely generated, integral and  saturated. 
		\end{itemize}  
		All log (formal) schemes considered in this article will be quasi-coherent and integral.

		\item For a prelog ring $(A,M)$, we write $(A,M)^a=(A,M^a)$ for the associated log ring, and $(\Spec(A), M)^a=(\Spec(A), \underline{M}^a)$ for the assocaited log scheme. If $A$ is $I$-adically complete for some finitely generated ideal $I$, we also write $(\Spf(A), M)^a$ for the associated log formal scheme.

		\item A morphism of log (formal) schemes $f\colon (X,\calM_X)\to (Y, \calM_Y)$ is called \emph{strict}, if the natural map of log structures $f^*(\calM_Y)\xra{\sim} \calM_X$ is an isomorphism. Obviously, a strict morphism of quasi-coherent log (formal) schemes admits a chart \'etale locally.


		\item	For   a category $\scrC$, we write $\scrC^{\opp}$ for the opposite category.
		
		\item For a site $\calC$, we denote by $\Shv(\calC)$, the associated  topos, i.e. the category of sheaves of sets on $\calC$. 
		
		\item By a perfectoid ring, we mean an integral perfectoid ring in the sense of \cite[Def. 3.5]{BMS1}.  For a perfectoid ring $R$, we denote by 
		\[
		R^{\flat}:=\varprojlim_{x\mapsto x^p} R/p
		\]
		its tilt, which is a perfect ring in characteristic $p$. For $x=(x_0, x_1, \cdots)\in R^{\flat}$, we put 
		\[
		x^{\sharp}:=\lim_{n\to+\infty}\tilde x_n^{p^n},
		\]
		where $\tilde x_n\in R$ is an arbitrary lift of $x_n\in R/p$. 	Put $\Ainf(R):=W(R^{\flat})$, and  let  
		$
		\theta_{R}:\Ainf(R)\to R
		$
		denote Fontaine's surjection defined by $\theta_R([x])=x^{\sharp}$ for any $x\in R^{\flat}$. 
		
		
	\end{itemize}

	\section{Preliminaries on log prismatic site}\label{Section1}
	
	\subsection{Log prisms}\label{S:log prisms}  According to \cite[Def. 2.2]{Kos},  a $\delta_{\log}$-ring is a tuple $(A, \delta, \alpha: M\to A, \delta_{\log})$, where $(A,\delta)$ is a $\delta$-ring, $\alpha: M\to A$ is a prelog structure, and $\delta_{\log}:M\to A$ is a map such that 
	\begin{itemize}
		\item $\delta_{\log}(e)=0$, 
		\item $\delta(\alpha(m))=\alpha(m)^p\delta_{\log}(m)$,
		\item $\delta_{\log}(mm')=\delta_{\log}(m)+\delta_{\log}(m')+p\delta_{\log}(m)\delta_{\log}(m')$
	\end{itemize}
	hold for all $m,m'\in M$. 
	Morphisms of $\delta_{\log}$-rings are defined in an evident way. If there is no confusion, we write a $\delta_{\log}$-ring simply by $(A,M)$ or by $(A,\alpha: M\to A)$.

	Recall also that \cite[Def. 3.3]{Kos}, a \emph{(bounded) prelog prism} is a triple $(A,I, M)$ where
	\begin{itemize}
		\item $(A,I)$ is a (bounded) prism in the sense of  \cite[Def. 3.2]{BS}, 
		\item $(A,M)$ is a $\delta_{\log}$-ring.
	\end{itemize}

	A (bounded) \emph{log prism} is a bounded prelog prism up to taking the associated log structure, i.e.  a (bounded) log prism consists of a (bounded) prism $(A,I)$, together with a  log structure $\alpha \colon M_{\Spf(A)}\to \cO_{\Spf(A)_{\et}}$  and a (sheaf version of) $\delta_{\log}$-structure $\delta_{\log}\colon M_{\Spf(A)}\to \cO_{\Spf(A)_{\et}}$ such that \'etale locally  $\alpha$ and $\delta_{\log}$  come from some prelog prism $(A,I,M)$. We write 
	\[
	(A,I, M_{\Spf(A)})=(A,I,M)^a
	\]
	for  the log prism attached to a prelog prism $(A,I,M)$. 
	For a log prism $(A,I, M_{\Spf(A)})$, we write $ M_{\Spf(A/I)}$ for the pullback to $\Spf(A/I)$ of the log structure $M_{\Spf(A)}$.

	A  morphism of bounded log prisms $f\colon (A,I, M_{\Spf(A)})\to (B, J, M_{\Spf(B)})$ is  a map of formal schemes $(\Spf(B), M_{\Spf(B)})\to (\Spf(A),M_{\Spf(A)})$ that is compatible with the $\delta_{\log}$-structures and induces a  moprhism of prisms $(A,I)\to (B,J)$.
By the rigidity of morphisms of prisms, we have $J=IB$. 

	A prelog prism $(A,I,M)$ (resp. a log prism $(A,I,M_{\Spf(A)})$) is called \emph{integral} if  $M$ (resp.  $M_{\Spf(A)}$) is integral.

	\begin{remark}\label{R:morphism log prisms}
	(1)	Let $(A,I,M_{\Spf(A)})$ be an integral log prism. Assume that $(\Spf(A/I), M_{\Spf(A/I)})$ admits a chart $M_{\overline A}\to \Gamma(\Spf(A/I), M_{\Spf(A/I)})$. Then, by \cite[Lemma~2.1]{DLMS}, 
		\[
		M_{A}:=M_{\overline A}\times_{\Gamma(\Spf(A/I), M_{\Spf(A/I)})}\Gamma(\Spf(A), M_{\Spf(A)})\to \Gamma(\Spf(A), M_{\Spf(A)})
		\]
		is a chart of  $(\Spf(A), M_{\Spf(A)})$ and the natural $\delta_{\log}$-structure
		\[
		M_A\to \Gamma(\Spf(A), M_{\Spf(A)})\xra{\delta_{\log}} A
		\]
		makes $(A,I,M_A)$ a prelog prism such that $(A,I,M_{\Spf(A)})=(A,I,M_A)^a$. Note that $M_A\to M_{\overline A}$ is a torsor under the group $1+I\subset A^{\times}$. If $(m_i)_{i\in I}$  are elements of $M_A$ such that their images in $M_{\overline A}$ generate the monoid $M_{\overline A}$, then the submonoid $P\subset M_A$ generated by  $(m_i)_{i\in I}$ is also a chart for $(\Spf(A), M_{\Spf(A)})$. 
		
	(2)	If $f\colon (A, I, M_{\Spf(A)})\to (B, J, M_{\Spf(B)})$ is a morphism of integral log prisms such that $(\Spf(B/J), M_{\Spf(B/J)})\to (\Spf(A/I), M_{\Spf(A/I)})$ admits a chart:  
		\[
		\xymatrix{M_{\overline A}\ar[r]\ar[d] &M_{\overline B}\ar[d]\\
			\Gamma(\Spf(A/I), M_{\Spf(A/I)})\ar[r] &	\Gamma(\Spf(B/J),M_{\Spf(B/J)}).}
		\]
	By the discussion in (1),   $f$ is induced by a morphism of prelog prisms $(A, I, M_A)\to (B, J, M_B)$. 
	\end{remark}

	Let  $(A,M)$ be a $\delta_{\log}$-ring such that $p$ lies in the Jacobson radical of $A$ and  $\alpha: M\to A$ is a log-structure, i.e. $\alpha$ induces an isomorphism  $\alpha^{-1}(A^{\times})\xra{\sim}A^{\times}$. Then the map  \[m\mapsto m^p\alpha^{-1}(1+p\delta_{\log}(m))\] defines  an endomorphism 
	$
	\phi_M:M\to M.	$ Together with the lift of Frobenius $\phi_A(x):=x^p+p\delta(x)$ on $A$, one gets an endomorphism $\phi:=(\phi_A, \phi_M)$ of the log ring $(A,M)$.
	If $(A, I, M_{\Spf(A)})$ is a log prism, then the sheafification of this construction gives an Frobenius endomorphism on the affine log formal scheme $(\Spf(A), M_{\Spf(A)})$.
	

	\subsection{Semistable formal schemes over completion valuation rings}\label{S:semistable-setup}
	Let $L|\QQ_p$ be either a perfectoid field  or a complete discrete valuation field with perfect residue field.  Let $\cO_L$ be the ring of integers of $L$, and  $v_p\colon L\to \RR\cup \{+\infty\}$ denote the $p$-adic valuation. 
	
	Consider the prelog ring  $\uO_L:=(\cO_L, M_{\cO_L})$ with  $M_{\cO_L}:=\cO_L\backslash\{0\}$ naturally embeded in $\cO_L$.  Let $\Spf(\uO_L):=(\Spf(\cO_L), M_{\cO_L})^a$ be the associated  log formal scheme.  If $L$ is a perfectoid field, we also put $M_{\cO^{\flat}_L}:=\cO_{L}^{\flat}\backslash\{0\}$ and consider the prelog ring $M_{\cO_L^{\flat}}\to \cO_L$ given by $x\mapsto x^{\sharp}$. As $x\mapsto x^{\sharp}$ induces an isomorphism of monoids  $M_{\cO^{\flat}_L}/M_{\cO^{\flat}_L}^{\times}\xra{\sim} M_{\cO_L}/M_{\cO_L}^{\times}$, we see that $\Spf(\uO_L)=(\Spf(\cO_L), M_{\cO_{L}^{\flat}})^a$.

\begin{definition}\label{D:semi-stable}
	(1)	An  affine $p$-adic formal scheme $\Spf(R)$ over $\Spf(\cO_L)$ is called \emph{small semi-stable} if there exists a 
	$p$-completely \'etale map (called a framing)
	\begin{equation}\label{E:small-framing}
		\square\colon R^{\square}:=\cO_L\langle T_0,\cdots, T_r, T_{r+1}^{\pm 1}, \cdots T_{d}^{\pm 1}\rangle/(T_0\cdots T_r-\pi)\to R
	\end{equation}
	for some integers $0\leqslant r\leqslant d$, where  $\pi\in L$ is an element with $0<v_p(\pi)\leqslant 1$ if $L$ is perfectoid, and it is  a uniformizer if $L$ has discrete valuation.
	
	(2) A $p$-adic formal scheme $X$ over $\Spf(\cO_L)$ is called semi-stable if it has an \'etale cover by affine small semistable formal schemes over $\Spf(\cO_L)$. 
	
	\end{definition}
	\begin{remark}\label{R:pseudo-uniformizer}
		If $L$ is a perfectoid field and $\pi\in L$ with $0<v_p(\pi)\leqslant 1$, then  \cite[Lemma 3.9]{BMS1} implies that there exists a unit $u\in \cO_L^{\times}$ and $\varpi^{\flat}\in \cO_L^{\flat}$ such that $\pi=\varpi^{\flat, \sharp}u$. Note that  
				\[
\cO_L\langle T_0,\cdots, T_r, T_{r+1}^{\pm 1}, \cdots T_{d}^{\pm 1}\rangle/(T_0\cdots T_r-\varpi^{\flat , \sharp})\simeq \cO_L\langle T_0,\cdots, T_r, T_{r+1}^{\pm 1}, \cdots T_{d}^{\pm 1}\rangle/(T_0\cdots T_r-\pi).
		\]
		Up to replacing $\pi$ by $\varpi^{\flat, \sharp}$, we may assume that $\pi$ in \eqref{E:small-framing} lies in the image of $M_{\cO_L^{\flat}}\xra{x\mapsto x^{\sharp}} \cO_L$. 

	\end{remark}

	Let $X$ be a  semistable $p$-adic formal scheme over $\Spf(\cO_L)$.
We equip $X$ with the canonical log structure $\calM_X:=\cO_{X_\et}\cap \cO_{X_{\et}}[1/p]^{\times}$, and write $\uX:=(X,\calM_X)$. We call $\uX$ a semistable $p$-adic log formal scheme over $\Spf(\cO_L)$.

		\subsection{Local chart}\label{S:local-chart} Assume that $X=\Spf(R)$  is affine small semistable with a framing \eqref{E:small-framing}. We will give an explicit chart for $\uX$. We distinguish the two cases of $L$:
		
		\begin{enumerate}
			\item[(1)] Assume first that $L$ has discrete valuation. Then $\uX$ admits a chart 
			\[
			\alpha_R\colon M_{R}:=\NN^{r+1}\to \Gamma(X, \calM_X)=R\cap R[1/p]^{\times}
			\]
			given by $e_i\mapsto T_i$ for $0\leqslant i\leqslant r$, where $\{e_i:0\leqslant i\leqslant r\}$ denotes the canonical basis of $\NN^{r+1}$.
			\item[(2)] Assume  now that $L$ is a perfectoid field. According to Remark~\ref{R:pseudo-uniformizer}, we may assume that $\pi$ in \eqref{E:small-framing} has the form $\pi=\pi^{\flat, \sharp}$ for some $\pi^{\flat}\in M_{\cO_L^{\flat}}$.  Consider the monoid $M_{R}: =M_{\cO_L^{\flat}}\sqcup_{\NN}\NN^{r+1}$, where $\NN\to \NN^{r+1}$ is the diagonal embedding and $\NN\ra M_{\cO_L^{\flat}}$ is given by $1\mapsto \pi^{\flat}$.  Then $\uX$ admits a chart 
			\[
			\alpha_R\colon 	M_{R} =M_{\cO^{\flat}_L}\sqcup_{\NN}\NN^{r+1}\to \Gamma(X, \calM_X)=R\cap R[1/p]^{\times},
			\]
			such that $\alpha_R(x)=x^{\sharp}$ for $x\in M_{\cO_L^{\flat}}$ and $\alpha_R(e_i)=T_i$ for $0\leqslant i\leqslant r$. 
			In particular, one has a cartesian diagram of integral quasi-coherent log formal schemes 
			\[
			\xymatrix{\uX\ar[r]\ar[d] &\uX_0 \ar[d]\\
				\Spf(\uO_L)\ar[r] & \Spf(\uO_L)_0,
			}
			\]
			where  $\Spf(\uO_L)_0$ is the formal scheme attached to the prelog ring  $\NN\xra{1\mapsto \pi} \cO_L$, and $\uX_0:=(\Spf(R), \NN^{r+1})^a$ is attached to the prelog ring $\NN^{r+1}\to R$ defined by $e_i\mapsto T_i$. 
			
		\end{enumerate}

We recall the absolute prismatic site of $\uX$ introduced in \cite[Def. 7.34]{KY} and \cite[Def. 2.3]{DLMS}.
\begin{definition}\label{D:abs-prismatic-site}
Let $\uX$ be a semistable   $p$-adic log formal scheme over $\Spf(\cO_L)$ as above.
 An object of  the \emph{absolute log prismatic site} $\underline X_{\prism}$ is a diagram  of  morphisms of log $p$-adic formal schemes 
		\[
		(\Spf(A), M_{\Spf(A)})\hookleftarrow  (\Spf(A/I), M_{\Spf(A/I)})\xra{f} \uX,
		\]
		where $(A,I, M_{\Spf(A)})$ is an integral bounded log prism, $M_{\Spf(A/I)}$ is the pullback of $M_{\Spf(A)}$ to $\Spf(A/I)$, and  $f$ is a strict morphism of log formal schemes, i.e. $M_{\Spf(A/I)}\cong f^*\calM_X$. 
		We write such an object as $(\Spf(A), I, M_{\Spf(A)}; f)$ or simply $(A, I, M_{\Spf(A)})$ if there is no confusion. We also say  that   $(A, I, M_{\Spf(A)})$ is an object of $\uX^{\opp}$ if we want to invert morphisms in $\uX_{\prism}$. 
		
		A	morphism $(\Spf(B), J, M_{\Spf(B)};g)\to (\Spf(A), I, M_{\Spf(A)};f)$ in  $\uX_{\prism}$ is a  commuative diagram 
		\[
		\xymatrix{(\Spf(B), M_{\Spf(B)})\ar[d] & (\Spf(B/J), M_{\Spf(B/J)})\ar[l]\ar[r]^-{g}\ar[d] & \uX\ar@{=}[d] \\
			(\Spf(A), M_{\Spf(A)}) & (\Spf(A/I), M_{\Spf(A/I)})\ar[l]\ar[r]^-{f} & \uX}
		\]	
		compatible with all structures. In this case, we will also say that  $(A, I,M_{\Spf(A)})\to (B, J, M_{\Spf(B)})$ is a morphism in $\uX^{\opp}$ for simplicity. Note that by rigidity of prisms, we have $J=IA$.

	We equip $\uX_{\prism}$ with the  flat topology, i.e.	a morphism $(\Spf(B), J, M_{\Spf(B)};g)\to (\Spf(A), I, M_{\Spf(A)};f)$ in $\underline X_{\prism}$ (equivalently a map $(A,I,M_{\Spf(A)})\to (B, J, M_{\Spf(B)})$ in $\uX^{\opp}$)  is a \emph{flat cover}  if the underlying map of $\delta$-rings $A\to B$ is  $(p,I)$-completely   faithfully flat.
		
			If $\underline R:=(R,M)$ is a $p$-adically complete  prelog ring, we write usually  $\underline R_{\prism}$  for the absolute prismatic site of its associated log formal scheme $(\Spf(R), M)^a$. 
		
		\end{definition}
		
			\begin{remark}\label{R:chart for log prism}
			(1) Assume that  $\uX$ admits a global chart $P\to \Gamma(X,\calM_X)$. Then for any object $(\Spf(A),I, M_{\Spf(A)};f)$ in $\uX_{\prism}$, $P\to \Gamma(X, \calM_X)\xra{f^*}\Gamma(\Spf(A/I), M_{\Spf(A/I)})$ defines a chart for $(\Spf(A/I), M_{\Spf(A/I)})$. Applying the construction in Remark~\ref{R:morphism log prisms} with $M_{\overline A}=P$, we get a  prelog prism $(A,I,M_A)$ such that $(A,I,M_A)^a=(A,I, M_{\Spf(A)})$. 
			This construction is  funtorial in $(A, I, M_{\Spf(A)})$ so that any morphism  $(\Spf(B), IB, M_{\Spf(B)};g)\to (\Spf(A), I, M_{\Spf(A)};f)$ in $\underline X_{\prism}$ is induced by a map of integral prelog prisms $(A,I,M_A)\to (B,IB, M_B)$. 
			
			(2) Consider a diagram of morphisms in $\uX^{\opp}$:
			\[
			(C, IC, M_{\Spf(C)})	\xleftarrow{g}(A, I,M_{\Spf(A)})\xra{f}(B, IB, M_{\Spf(B)})
			\]
			where  one of $f$ and $g$ is a flat cover. Then the coproduct of $f$ and $g$ in $\uX^{\opp}_{\prism} $ exists\footnote{Note that this is part of the requirements for $\uX_{\prism}$ to be a site.} and can be described explicitly as follows. We may assume that $f$ is $(p,I)$-completely faithfully flat.  Let $D=(B\otimes^L_{A}C)^{\wedge}$ be the derived $(p,I)$-completion of $B\otimes_A^LC$. By \cite[Lemma 3.7]{BS},  $D$ coincides with the classical $(p,I)$-adic completion of $B\otimes_AC$ and is equipped with the canonical tensor $\delta$-strucutre so that $(B,IB)\to(D,ID)$ is a $(p,I)$-completely faithfully flat map of prisms. 
			
			Let $M_{\Spf(D)}$ be the log structure on $\Spf(D)$ attached to the prelog structure  $q_B^{-1}M_{\Spf(B)}\sqcup_{q_A^{-1}M_{\Spf(A)}}q_C^{-1}M_{\Spf(C)}\to \cO_{\Spf(D)_{\et}}$ where $q_{?}: \Spf(D)\to \Spf(?)$ with $?\in \{A,B, C\}$ is the canonical projection. Then there is a canonical (sheaf version of) $\delta_{\log}$-structure $\delta_{\log}\colon M_{\Spf(D)}\to \cO_{\Spf(D)_{\et}}$.
			We need to show that, \'etale locally, there exists a prelog prism $(D, ID, M_D)$ such that $(D, ID, M_{\Spf(D)})=(D, ID, M_D)^a$. The problem being local for the \'etale topology on $X$, we may assume that $\uX$ admits a global chart $P\to \Gamma(X, \calM_X)$. By part (1) of this remark, $f$ and $g$ come from morphisms of prelog prisms $(A, I, M_A)\to (B, IB, M_B)$ and $(A,I,M_A)\to (C, IC, M_C)$. Put $M_{D}=M_B\sqcup_{M_A}M_C$. Then there exist a natural $\delta_{\log}$-structure on $(D,M_D)$ so that $(D,ID, M_D)^a=(D, ID, M_{\Spf(D)})$. 
		\end{remark}
		
	When $X=\Spf(R)$ is affine and small semistable, 	we  have the following description for objects of $\uX_{\prism}$.

	 \begin{lemma}\label{L:local objects}
		
		Assume that $X=\Spf(R)$ is affine small semistable with framing \eqref{E:small-framing}, and $\alpha_R\colon M_R\to \Gamma(X, \calM_X)$ be the  chart constructed in \S~\ref{S:local-chart}. Let $(\Spf(A), I, M_{\Spf(A)};f)$ be an object of $\uX_{\prism}$. Then the chart  $M_R\xra{\alpha_R} \Gamma(X, \calM_X)\xra{f^*} \Gamma(\Spf(A/I), M_{\Spf(A/I)})$ lifts to  a chart  $\alpha_A\colon M_R\to \Gamma(\Spf(A), M_{\Spf(A)})
		$ for $M_{\Spf(A)}$ so that  the diagram of monoids 
		\[
		\xymatrix{M_R\ar[rr]^-{\alpha_R}\ar@{=}[d] && R\ar[d]^{f^*}\\
		M_R\ar[r]^{\alpha_A} & A\ar[r] & A/I}
		\]	
		is commutative, and  $(A,M_R)$ has a natural induced $\delta_{\log}$-structrure   such that  $(A,I, M_R)^a=(A,I, M_{\Spf(A)})$.  
	\end{lemma}
	\begin{proof}
		By  Remark~\ref{R:chart for log prism}(1),  the log prism $(A, I, M_{\Spf(A)})$ comes from a prelog prism $(A,I,M_A)$ with $M_A=M_R\times_{\Gamma(\Spf(A/I), M_{\Spf(A/I)})}\Gamma(\Spf(A), M_{\Spf(A)})$.
		Note that $M_A\to M_R$ is  a torsor under the group $1+I\subset A^{\times}$. We now distinguish two cases:
		
		(1) When $L$ has discrete valuation, then   $M_R=\NN^{r+1}$ is a free monoid  and we  can easily choose a section $s\colon M_R\to M_A$ to $M_A\to M_R$. We define $\alpha_A$ as the composite map 
		\[\alpha_A\colon M_R\xra{s}M_A\to \Gamma(\Spf(A), M_{\Spf(A)}),
		\]
		and equip $(A, M_R)$ with the $\delta_{\log}$-structure  induced from  $(A, \Gamma(\Spf(A), M_{\Spf(A)}))$ by pre-composing with $\alpha_A$. Then it is easy to see that $\alpha_A$ satisfies the requirement of the Lemma. 
		
		(2) Consider now the case when $L$ is a perfectoid field. Similar to the previous case, it suffices to construct a section $s:M_R=M_{\cO_L^{\flat}}\sqcup_{\NN}\NN^{r+1}\to M_A$ to the canonical surjection $M_A\to M_R$. 
		
		Since the prism $(\Ainf(\cO_L), \ker(\theta_{\cO_L}))$ is an initial object of the non-logarithmic prismatic site  $(\cO_L)_{\prism}^{\opp}$, there exists a unique map of prisms $\iota\colon (\Ainf(\cO_L), \ker(\theta_{\cO_L}))\to (A,I)$. We claim that $\iota$ extends to a map of prelog prisms 
		$$
		\iota\colon (\Ainf(\cO_L), \ker(\theta_{\cO_L}), M_{\cO_{L}^{\flat}})\to (A,I, M_A).$$ 
		Indeed,   for $x\in M_{\cO^{\flat}_L}=\cO_L^{\flat}\backslash \{0\}$, we put 
		\[\iota(x)=\lim_{n}\tilde x_n^{p^n}\in M_A,\]
		where $\tilde x_n\in M_A$ is an arbitrary lift of $x^{1/p^n}\in M_{R}$. As $(1+I)=\varprojlim_{n\to +\infty}(1+I)/(1+I)^{p^n}$,  it is easy to see that such a definition makes sense and $\iota$ is compatible with the $\delta_{\log}$-structures on both sides. 
		
		Now, let $m_i\in M_A$ be a lift of $e_i\in \NN^{r+1}\subset M_R$ with $0\leqslant i\leqslant r$. 
		Since both  $\prod_{i=0}^{r}m_i$ and $\iota(\pi^{\flat})$ have  the same image in $M_{R}$, there exists $u\in 1+I$ such that $\iota(\pi^{\flat})=u\prod_{i=0}^{r} m_i$. Up to replacing $m_0$ by $m_0u$, we may assume that $\iota(\pi^{\flat})=\prod_{i=0}^rm_i$. 
		Then the map  $s\colon M_R\to M_A$  given by 
		$$s(x\sqcup \sum_{i=0}^ra_ie_i)= \iota(x)\prod_{i=0}^{r}m_i^{a_i}, \quad \text{for all } x\in M_{\cO_L^{\flat}}, a_i\in \NN.$$
		defines a section to  $M_A\to M_R$. This finishes the proof of the Lemma.

	\end{proof}

  				\subsection{Relative Variant}\label{S:relative-prismatic-site}	The absolute log prismatic site $\uX_{\prism}$ has an obvious relative variant: Let  $\underline A:=(A,I,M_{\Spf(A)})$ be a fixed integral bounded log prism in $\Spf(\uO_L)_{\prism}$.
  			The \emph{relative log prismatic site} $(\underline X/\underline A)_{\prism}$ is  the opposite of the category of bounded  log prisms $(B,IB,\Spf(M_B))$ over $(A,I,M_{\Spf(A)})$ together with a commuative  diagram of log formal schemes 
  			\[\xymatrix{(\Spf(B/IB), M_{\Spf(B/IB)})\ar@{^(->}[r]\ar[d]^{f} & (\Spf(B), M_{\Spf(B)})\ar[dd]\\
  				\uX\ar[d] \\
  				(\Spf(A/I), M_{\Spf(A/I)})^a\ar@{^(->}[r] & (\Spf(A), M_{\Spf(A)}),
  			}
  			\] 
  			where $f$ is strict. 
  			We equip $(\underline X/\underline A)_{\prism}$ with the   flat topology as in the absolute case. We denote also by $\cO_{\prism}$ the structural sheaf on $(\underline X/\underline A)_{\prism}$.

  			\begin{remark}\label{R:abs-relative-prismatic-site}
  			Assume $L$ is perfectoid. Consider the prelog ring $M_{\cO_L^{\flat}}\to \Ainf(\cO_L)$ defined by $x\mapsto [x]$. Then the associated log prism $\Ainf(\uO_L)=(\Ainf(\cO_L), \ker(\theta_{\cO_L}), M_{\cO_L^{\flat}}\xra{x\mapsto [x] }\Ainf(\cO_L))^a$ is naturally an object of $\Spf(\uO_L)_{\prism}$. 
  	By the proof of Lemma~\ref{L:local objects},  
  				 for every object $(B,J,M_{\Spf(B)})$ in $\underline X_{\prism}$, there exists a unique  morphism of log prisms $\Ainf(\uO_L)\to (B,J,M_{\Spf(B)})$ such that $(B, J,M_{\Spf(B)})$ becomes an object of $(\underline X/ \Ainf(\uO_L))_{\prism}$.  Therefore,   the site $\underline X_{\prism}$ is equivalent to the relative log prismatic site $(\underline X/ \Ainf(\uO_L))_{\prism}$ in this case. More generally, for any integer $n\in\ZZ$, the log prism 
				 \[
				 \mathrm{A}_{\inf}^{(n)}(\uO_L):=(\Ainf(\cO_L), \ker(\theta_{\cO_L}\circ\phi^{-n}), M_{\cO_L^{\flat}}\xra{x\mapsto [x^{p^n}] }\Ainf(\cO_L))^a
				 \]
				 is isomorphic to $\Ainf(\uO_L)$. Hence, we have also an equivalence of sites $\uX_{\prism}\simeq (\uX/\mathrm{A}_{\inf}^{(n)}(\uO_L))_{\prism}$.

  			\end{remark}

	\subsection{Structural sheaves and crystals} By $(p,I)$-complete faithfully flat descent, the presheaf on  $
	\underline X_{\prism}$
	\[\cO_{\prism}=\cO_{\underline X_{\prism}}: (\Spf(B),J,M_{\Spf(B)};f)\mapsto  B\]
	is a sheaf of rings.  
	There exists a canonical lift of Frobenius  endomorphism $\phi_{\cO_{\prism}}$ on  $ \cO_{\prism}$.  As in the non-log case \cite[Remark 2.4]{BS2}, the topos $\Shv(\underline X_{\prism})$ is replete, because an inductive limit of strict faithfully flat maps of bounded log prisms is again strict faithfully flat. 
	
	Let 
	$\calI_{\prism}\subset \cO_{\prism}$ the ideal sheaf sending $(B,J, M_{\Spf(B)})$ to $J$.  
	For an integer $n\geqslant1$, put 
	\[
	\cO_{\prism,n}:=\cO_{\prism}/(p,\calI_{\prism})^n.
	\]
	Then by the fpqc-descent, we have, for any object $(B,J,M_{\Spf(B)})$ of $\underline X_{\prism}$,
	\[
	\cO_{\prism, n}(B,J,M_{\Spf(B)})=B/(p,J)^n,
	\]
	and $H^i((B,J,M_{\Spf(B)}), \cO_{\prism,n})=0$ for all $i>0$. It follows from \cite[Lemma 3.18]{Sch1} that 
	\[\cO_{\prism}\simeq \varprojlim_n \cO_{\prism,n},\quad R^i\varprojlim_n \cO_{\prism,n}=0,\; \text{ for all $i>0$.}\]
	Similarly, let $\cO_{\prism}[{1}/{\calI_{\prism}}]^{\wedge}_p$ be the $p$-adic completion of $\cO_{\prism}[{1}/{\calI_{\prism}}]$.
	Then the Frobenius endomorphism extends naturally to $\laurent$, and we have 
	\[\laurent\simeq  \varprojlim_n \cO_{\prism}[1/\calI_{\prism}]/(p^n), \quad 
	R^i\varprojlim_n \cO_{\prism}[1/\calI_{\prism}]/(p^n)=0,\; \text{for $i>0$.}\]

	As in the non-log case (cf. \cite[Def. 11.1]{Tsuji}), we have the following notion of  complete crystals on $\underline X_{\prism}$.
	\begin{definition}\label{D:crystals}
		(1) For an integer $n\geqslant1$, a \emph{crystal of $\cO_{\prism, n}$-modules} on $\underline X_{\prism}$ is a presheaf $\calF$ of $\cO_{\prism,n}$-modules such that for every morphism $u:(B,J_B,M_{\Spf(B)})\to (C,J_C,M_{\Spf(C)})$ in $\underline X_{\prism}^{\opp}$,  $\calF(u)$ induces  an isomorphism 
		\[
		\calF(B,J_B,M_{\Spf(B)})\otimes_{B,u}C\xra{\sim} \calF(C,J_C,M_{\Spf(C)}).
		\]
		Let $\CR(\underline X_{\prism}, \cO_{\prism,n})$ be the category of crystals of $\cO_{\prism,n}$-modules. Note that any crystal of  $\cO_{\prism,n}$-modules is automatically a sheaf by fpqc-descent of quasi-coherent sheaves. 
		
		(2) \emph{A complete crystal of $\cO_{\prism}$-modules} on $\underline X_{\prism}$ is a presheaf of $\cO_{\prism}$-modules such that $\calF/(p,\calI_{\prism})^n$ is a crystal of $\cO_{\prism,n}$-modules for each $n\geqslant1$ and the canonical map  $$\calF\simeq \varprojlim_n\calF/(p,\calI_{\prism})^n$$
		is an isomorphism of presheaves. 
		A complete crystal $\calF$ of $\cO_{\prism}$-modules on $\underline X_{\prism}$ is already a sheaf, and for every morphism  $u:(B,J_B,M_{\Spf(B)})\to (C,J_C,M_{\Spf(C)})$ in $\underline X_{\prism}$, the canonical map 
		\begin{equation}\label{E:Crystal-isom}
			c_{\calF}(u)\colon \calF(B,J_B,M_{\Spf(B)})\widehat\otimes_{B,u}C\xra{\sim} \calF(C,J_C,M_{\Spf(C)})
		\end{equation}
		is an isomorphism, where the left hand side is the classical $(p,J_B)$-adic completion of $\calF(B,J_B,M_{\Spf(B)})\otimes_{B,u}C$.
		We write  $\widehat{\CR}(\underline X_{\prism},\cO_{\prism})$ for the category of complete $\cO_{\prism}$-crystals on $\underline X_{\prism}$.
		
	\end{definition}
	
	\begin{remark}\label{R:adic-crystal}

		Following  \cite[Remark 11.2]{Tsuji},  we  define  an adic $\cO_{\prism}$-crystal on $\underline X_{\prism}$ to be    an  inverse system of $\cO_{\prism,n}$-crystals $(\calF_n)_{n\in \NN}$ such that the transition map $\calF_{n+1}\to \calF_n$ induces an isomorphism $$\calF_{n+1}\otimes_{\cO_{\prism,n+1}}\cO_{\prism,n}\xra{\sim}\calF_n.$$ Let $\CR^{\ad}(\underline X_{\prism}, \cO_{\prism})$ be the category of adic $\cO_{\prism}$-crystals on $\underline X_{\prism}$.
		Then  the natural functor 
		$$\CRhat(\underline X_{\prism}, \cO_{\prism})\xra{\sim} \CR^{\ad}(\underline X_{\prism}, \cO_{\prism}), \quad \calF\mapsto (\calF\otimes_{\cO_{\prism}}\cO_{\prism,n})_{n\in \NN}$$
		is  an equivalence of categories (cf. \cite[Lemma~4.9, Remark~11.2]{Tsuji}).
		
	\end{remark}
	
	The following notion is a  more practical version of analytic prismatic crystals considered in \cite{GR} and \cite{DLMS}.

	\begin{definition}\label{D:complete-prismatic-F-crystal}

		(1) Let $\CRhat^{\anfr}(\uX_{\prism}, \cO_{\prism})$ be the full subcategory of $\CRhat(\uX_{\prism}, \cO_{\prism})$ consisting of objects $\calF$ such that for every object $(B, J_B,M_{\Spf(B)})$ of $\uX_{\prism}$, the restriction to $\Spec(B)\backslash V(p, J_B)$ of the quasi-coherent sheaf on $\Spec(B)$ attached to $\calF(B, J_B,M_{\Spf(B)})$  is a vector bundle. 
		
		(2) A \emph{Frobenius structure} on an object $\calF$ of $\CRhat^{\anfr}(\uX_{\prism}, \cO_{\prism})$ is  an isomorphism  
		\[
		\phi_{\calF}:(\phi^*_{\cO_{\prism}}\calF)[\frac{1}{\calI_{\prism}}]\xra{\sim}\calF[\frac{1}{\calI_{\prism}}], 
		\]
		where $\phi^*_{\cO_{\prism}}\calF=\calF\widehat{\otimes}_{\cO_{\prism}, \phi_{\cO_{\prism}}}\cO_{\prism}$. We denote by $\CRhat^{\anfr}(\uX_{\prism}, \cO_{\prism})^{\phi=1}$ the category of objects of $\CRhat^{\anfr}(\uX_{\prism}, \cO_{\prism})$ equipped with a Frobenius structure. 
	\end{definition}

	We will give now some general discussion on the local description of objects of $\CRhat(\uX_{\prism}, \cO_{\prism})$. 
	
	\begin{definition}\label{D:versal object}
		An object $\uA:=(A, I, M_{\Spf(A)})$  of $\uX_{\prism}$
		is a \emph{strong covering object} if it satisfies the following condition: For every object $(B, J, M_{\Spf(B)})$ of $\uX_{\prism}$, the coproduct 
		\[
		(D_B, ID_B, M_{\Spf(D_B)}):=(A, I,M_{\Spf(A)})\coprod (B,J,M_{\Spf(B)})
		\]
		in the category  $\uX_{\prism}^{\opp}$	exists, and the canonical map  of log prisms $(B,J,M_{\Spf(B)})\to (D_B, ID_B, M_{\Spf(D_B)})$ is a  flat cover. In particular, $\uA$ is a cover of the final object of the topos $\Shv(\uX_{\prism})$. 
	\end{definition}

	Let $\uA$ be a strong covering object of $\uX_{\prism}$, and 
	let $\uA^{\bullet}: =(A^\bullet, IA^{
		\bullet}, M_{\Spf(A^{\bullet})})$ be the \v{C}ech nerve of $\uA$ over the final object of $\Shv(\uX_{\prism})$. This is a simplicial object of $\uX_{\prism}$ whose $n$-th term is  the $(n+1)$-fold self-product of $\uA$ over  the final object of $\Shv(\uX_{\prism})$. 
	For an abelian sheaf $\calF$ on $\underline{X}_{\prism}$, let $\calF(\uA^{\bullet})$ denote the evaluation of $\calF$ on the simplicial object $\uA^{\bullet}$ of $\underline{X}_{\prism}$.  We denote by $\CA(\uA^{\bullet}, \calF)$ the simple  complex associated to the cosimplicial abelian group $\calF(\uA^{\bullet})$. 
	
	\begin{lemma}\label{L:CA-complex-prismatic}
		Under the above notation, and let  $\calE$ be an object of $\CRhat(\uX_{\prism}, \cO_{\prism})$. Then  $R\Gamma(\underline {X}_{\prism}, \calE)$ is computed by the \v{C}ech--Alexander complex $\CA(\uA^{\bullet}, \calE)$. 
	\end{lemma}
	
	\begin{proof}
		For each $m\geqslant 1$, let $\calE_{m}:=\calE\otimes_{\cO_{\prism}}\cO_{\prism,m}$. Then each  $\calE_m$ is an object of $\CR(\uX_{\prism}, \cO_{\prism,m})$, and $\calE\simeq \varprojlim_m \calE_m$ (Remark \ref{R:adic-crystal}). It follows from the crystal property of $\calE$ and   the fpqc-descent for quasi-coherent sheaves that  $H^i(\uA^n, \calE_{m})=0$ for $i>0$ and all $n,m$. As $\uA$ is a strong covering object,  $R\Gamma(\underline X_{\prism}, \calE_m)$ is computed by $\CA(\uA^{\bullet}, \calE_m)$. 
		As the topos $\Shv(\uX_{\prism})$ is replete, we have $\calE=R\varprojlim_m \calE_m$ and 
		\begin{align*}
			R\Gamma(\uX_{\prism}, \calE)&\simeq R\varprojlim_m R\Gamma(\uX_{\prism}, \calE_m)\\
			&\simeq R\varprojlim_m \CA(\uA^{\bullet}, \calE_m)\\
			&\simeq \CA(\uA^{\bullet}, \calE),
		\end{align*}
		where the last isomorphism follows from  $R\varprojlim_m \calE_m(\uA^n)=\calE(\uA^n)$ for each $n\geqslant0$.

	\end{proof}
	\subsection{Simplicial notation}\label{S:simplicial notation}
	For each integer $n\geqslant 0$, write $[n]:=\{0,1, \cdots, n\}$. Let $\Delta\colon [n]\to [0]$ denote the unique morphism. Let $p_i$ with $i=0,1$ denote the map $[0]\to [1]$ sending $0$ to $i$, and  for $(i,j)\in \{(0,1), (0,2), (1,2)\}$, let $p_{i,j}: [1]\to [2]$ be the increasing map sending $[1]$ to $\{i,j\}\subset [2]$. 
	By abuse of notation, we still use the same notation  to denote  the corresponding morphsims $\Delta: \uA^n\to \uA$, $p_i:\uA\to \uA^1$  and $p_{i,j}\colon \uA^1\to \uA^2$ in the category $\uX_{\prism}^{\opp}$ as well as  the undelying maps of $\delta$-$A$-algebras.	
	We have in particular a diagram of $\delta$-$A$-algebras:
	\[\begin{tikzcd}
		A \arrow[r, shift left, "p_0"]
		\arrow[r, shift right, "p_1"'] & A^1
		\arrow[r, shift left=2ex, "p_{0,1}"]
		\arrow[r,"p_{1,2}" description]
		\arrow[r, shift right=2ex, "p_{0,2}"']
		& A^2.
	\end{tikzcd}
	\]

	For a continuous homomorphism  of  $(p, I)$-adically complete $A$-algebras $f:S\to S'$ and a $(p, I)$-adically complete $S$-module $N$ (resp. an homomorphsim of   $(p, I )$-adically complete $S$-modules  $\psi: N\to N'$), we denote by $f^*(N)$ (resp. $f^*(\psi)$) the $(p,I)$-adically complete scalar extension $N\widehat\otimes_{S,f}S'$ (resp. $\psi\widehat\otimes_{S,f}S'$) under $f$.
	
	\begin{definition}\label{D:stratified-modules}
		Let $A^{\bullet}$ be the undelying cosimplicial ring of $\uA^{\bullet}$ as above. 
		A \emph{(complete) stratified module over $A^{\bullet}$} is a pair $(E, \epsilon)$, where 
		\begin{itemize}
			\item $E$ is  a  $(p, I)$-adically complete $A$-module,
			\item $\epsilon$ is an isomorphism of $A^1$-modules 
			\[\epsilon\colon  p_1^*(E)\xra{\sim} p_0^*(E)\]
			satisfying the following cocycle conditions
			\begin{enumerate}
				\item $\Delta^*(\epsilon)=\id_E$,
				\item $p_{0,1}^*(\epsilon)\circ p_{1,2}^*(\epsilon)=p_{0,2}^*(\epsilon)$. 
			\end{enumerate}
		\end{itemize}
		A morphism $f:(E,\epsilon)\to (E',\epsilon')$ of  stratified $A$-modules is a continuous morphism $f:E\to E'$ of $A$-modules compatible with stratifications: $p_0^*(f)\circ \epsilon=\epsilon'\circ p_1^*(f)$. Let $\Strat(A^{\bullet})$ denote the category of stratified  modules over $A^{\bullet}$.
	\end{definition}
	
	Let $\calE$ be an object of $\CRhat(\uX_{\prism},\cO_{\prism})$.  Then the crystal property of $\calE$ gives rise to a composite isomorphism:
	\[
	\epsilon_{\calE}=c_{\calE}(p_0)^{-1}\circ c_{\calE}(p_1)\colon \quad p_1^*\calE(\uA)\xra{\sim} \calE(\uA^1)\xra{\sim} p_0^*\calE(\uA)
	\]
	where  the isomorphisms $c_{\calE}(p_i)$ are defined in \eqref{E:Crystal-isom}.
	It is easy to verify that $(\calE(\uA), \varepsilon_{\calE})$ is a stratified module over $A^{\bullet}$.  
	\begin{proposition}\label{P:crystals-strat}
		Let $\uA$ be a strong covering  object of $\uX_{\prism}$ (Definition~\ref{D:versal object}).
		The construction $\calE\mapsto (\calE(\uA), \epsilon_{\calE})$ induces an equivalence of categories 
		\[
		\ev_{\uA^{\bullet}}\colon\quad  \CRhat(\uX_{\prism},\cO_{\prism})\xra{\sim} \Strat(A^{\bullet})
		\]
	\end{proposition}
	\begin{proof}
		In view of Remark~\ref{R:adic-crystal}, it suffices to show that for every integer $n\geqslant 1$, the construction $\calE\mapsto (\calE(\uA), \epsilon_{\calE})$ establishes an equivalence of categories between $\CR(\uX_{\prism}, \cO_{\prism, n})$ and $\Strat(A^{\bullet}_n)$, which is the full subcategory of $\Strat(A^{\bullet})$ consisting of $(E,\epsilon)$ with $E$ annihilated by $(p,I)^n$. 
		Then the arguments  are exactly the same as  \cite[11.10]{Tsuji}. 
	\end{proof}
	
	\subsection{Perfect log prismatic site}
Following \cite{MW}, we say that an object $(A, I, M_{\Spf(A)})$ of $\uX_{\prism}$ is perfect if $(A,I)$ is a perfect prism \footnote{Note that this is different from the convention  in \cite[Def. 2.35]{KY}, where  a ``log prism'' $(A, I, M_A)$ is perfect if the Frobenius endomorphism  of $(A,I, M_{A})$ is surjective in both rings and in monoids.}. Let $\uX^{\perf}_{\prism}$ denote the full subcategory of perfect objects $(A, I, M_{\Spf(A)})$ in $\uX_{\prism}$ equippd still with the  flat topology, and $X^{\perf}_{\prism}$ be the non-log version of perfect prismatic site. 

We recall the following useful Lemma.

\begin{lemma}[Min--Wang]\label{L:Min-Wang}
	Let $R$ be a perfectoid ring over $\cO_L$, and  $\varpi\in \cO_L\backslash\{0\}$. Let $a_1,a_2\in R^{\flat}$ and $x\in R$ such that $a_1^{\sharp}=a_2^{\sharp}(1+\varpi x)$. Then there exists $u\in R^{\flat, \times}$ such that $u^{\sharp}=1+\varpi x$.
	
	\end{lemma}
	\begin{proof}
When $\varpi=p$, this is exactly \cite[Lemma 2.15]{MW}. But as $\sqrt{pR}=\sqrt{\varpi R}$, the same argument as \emph{loc. cit.} applies to general $\varpi$.
	\end{proof}

	An important observation is the following Proposition due to Yu Min and Yupeng Wang. 
	\begin{proposition}[\cite{MW}, Prop. 2.18]\label{P:equiv-perfect-log-site}
		The forgetful functor $(A,I, M_{\Spf(A)})\mapsto (A,I)$ induces an equivalence of sites 
		\[
		\uX_{\prism}^{\perf}\simeq X^{\perf}_{\prism}.
		\]
	\end{proposition}

	\begin{proof}
		The case when $L$ has discrete valuation is exactly \cite[Prop. 2.18]{MW}. Now we adapt their arguments to the case when $L$ is perfectoid. 
		Let $(A,I)$ be a perfect prism in $X^{\perf}_{\prism}$. We need to prove that there exists a unique $\delta_{\log}$-structure $M_{\Spf(A)}$ on $\Spf(A)$ so that $(A,I, M_{\Spf(A)})$ becomes an object  of  $\uX_{\prism}^{\perf}$. 
		We may assume that $X=\Spf(R)$ is small semistable equipped with a framing as in \eqref{E:small-framing}.  Recall the chart $\alpha_R\colon M_R=M_{\cO_{L}^{\flat}}\sqcup_{\NN}\NN^{r+1}\to R$ for $\uX$ constructed  in \S~\ref{S:local-chart}.  
		First of all,  there exists a unique  map of perfect prisms $(\Ainf(\cO_L), \ker(\theta_{L}))\to (A,I)$ corresponding to the structural map $\cO_L\to A\to A/I$.

		We show first the uniqueness of $M_{\Spf(A)}$. The proof is similar to \cite[Lemma~2.17]{MW}. Let $M_{\Spf(A),1}$ and $M_{\Spf(A), 2}$ be two $\delta_{\log}$-structures on $(A,I)$ such that $(A, I, M_{\Spf(A), 1})$ and $(A, I, M_{\Spf(A), 2})$  both become objects of $\uX_{\prism}$.
		 By Lemma~\ref{L:local objects}, the chart $M_{R}=M_{\cO_L^{\flat}}\sqcup_{\NN}\NN^{r+1}\xra{\alpha_R} R\to A/I$ lifts to  charts $\alpha_{A,i}\colon M_R\to A$  for $i=1,2$ and  there exist $\delta_{\log}$-structures $ \delta_{\log, i}:M_R\to  A$ such that $(A, I, M_{\Spf(A), i})=(A,I, M_{R}, {\alpha_{A,i}, \delta_{\log, i}})^a$. Since $A$ is an integral domain,  $\delta_{\log, i}$ is uniquely determined by $\alpha_{A, i}$ because of the  requirement that 
		$$\delta(\alpha_{A,i}(m))=\alpha_{A,i}(m)^p\delta_{\log,i}(m)$$ for all $m\in M_R$. 
		Therefore, to finish the proof of the uniqueness, it suffices to show that $\alpha_{A,1}$ and $\alpha_{A,2}$ induce the same log structure on $\Spf(A)$. 
		 By the proof of Lemma~\ref{L:local objects}, the restrictions of $\alpha_{A,1}$ and $\alpha_{A,2}$ to $M_{\cO^{\flat}_L}\subset M_R$ coincide with $M_{\cO^{\flat}_L}\xra{x\mapsto [x]}\Ainf(\cO_L)\to A$.  
		 Put $t_{j,i}=\alpha_{A,i}(e_j)\in A$ for $0\leqslant j\leqslant r$. By \cite[Lemma 2.12]{MW}, we have $t_{j,i}=[a_{j,i}](1+px_{j,i})$ for some $a_{j,i}\in (A/I)^{\flat}$ and $x_{j,i}\in A$. As  $t_{j,1},t_{j,2}\in A$ are  both lifts of the image of $T_i$ in $A/I$ via the structural map $R\to A/I$, we have 
		\[
		a_{j,1}^{\sharp}(1+p\theta_{A/I}(x_{j,1}))=a_{j,2}^{\sharp}(1+p\theta_{A/I}(x_{j,2})),
		\]
 		where $\theta_{A/I}:A\cong \Ainf(A/I)\to A/I$ is Fontaine's surjection. By Lemma~\ref{L:Min-Wang}, there exists $u_{j}\in (A/I)^{\flat}$ such that $u_j=a_{j,1}a_{j,2}^{-1}$. It follows that $t_{j,1}$ and $t_{j,2}$ differ  by unit in $A$.  Hence, the log structures on $\Spf(A)$ induced by the charts $\alpha_{A,1}$ and $\alpha_{A,2}$ coincide. This finishes the proof of the uniqueness of $M_{\Spf(A)}$. 
 		
 		We prove now the existence of $M_{\Spf(A)}$. By the local uniqueness of $M_{\Spf(A)}$, we may assume that $X=\Spf(R)$ is small semistable with framing \eqref{E:small-framing}. It suffices to prove that there exists a prelog structure  $\alpha_A\colon M_R\to A$ such that 
 		\begin{itemize}
 			\item  $M_R\xra{\alpha_A}A\to A/I$ induces the same log structure on $\Spf(A/I)$ as $M_{R}\xra{\alpha_R}R\to A/I$,
 			\item  and   we have $\delta(\alpha_A(m))\in \alpha_A(m)^p A$ for all $m\in M_R$.
 		\end{itemize} We define $\alpha_{A}|_{M_{\cO_L^{\flat}}}$ as the composite map $M_{\cO_L^{\flat}}\xra{x\mapsto [x]}\Ainf(\cO_L)\to A$. 
 		 Let $\varpi=(\pi^{\flat})^{1/p, \sharp}\in \cO_L$, and $t_i\in A/I$ for $0\leqslant i\leqslant r$ be the image of $T_i\in R$ under the structural map $R\to A/I$. Then by \cite[Lemma~3.9]{BMS1}, there exists $t_i^{\flat}\in (A/I)^{\flat}$ such that $t_i^{\flat, \sharp}\equiv t_i\mod p\varpi (A/I)$. Write  $t_i^{\flat, \sharp}=t_i+p\varpi x_i$ for some $x_i\in A/I$.  As $\prod_{j=0}^r t_j=\pi=\pi^{\flat,\sharp}$ which divides $p$ by  assumption, there exists $y_i\in A/I$ for each $0\leqslant i\leqslant r$ such that $t_i y_i=p$. It follows that $t_i^{\flat, \sharp}=t_i+p\varpi x_i=t_i(1+\varpi x_iy_i)$, and 
 		\[
 		\prod_{i=0}^rt_i^{\flat, \sharp}=\pi^{\flat, \sharp} \prod_{i=0}^r(1+\varpi x_iy_i). 
 		\]
 		By Lemma~\ref{L:Min-Wang}, there exists $u\in (A/I)^{\flat, \times }$ such that $\prod_{i=0}^r t_i^{\flat, \sharp}=\pi^{\flat, \sharp} u^{\sharp}$. Up to replacing $t_0^{\flat}$ by $t_0^{\flat} u^{-1}$, we  get elements $t_i^{\flat}\in (A/I)^{\flat}$  for $0
 		\leqslant i\leqslant r$ such that each $t_i^{\flat,\sharp}$ differs from $t_i$ by a unit in $(A/I)$ and $\prod_{i=0}^r t_i^{\flat}=\pi^{\flat}$. We can   extend $\alpha_{A}|_{M_{\cO_L^{\flat}}}$ to a map $\alpha_{A}\colon M_R\to A$ by setting $\alpha_{A}(e_i)=[t_i^{\flat}]$ for $0\leqslant i\leqslant r$. 
 		Then it is clear that $\delta(\alpha_{A}(m))=0$ for all $m\in M_R$ and the prelog structure  $M_R\xra{\alpha_{A}} A\to A/I$ induces the same log structure on $\Spf(A/I)$ as $M_R\xra{\alpha_R}R\to A/I$. This finishes the proof. 
		\end{proof}
		
\subsection{Perfect versus non-perfect prismatic site}
 The natural inclusion functor 
\[
  \uX^{\perf}_{\prism}\to \uX_{\prism}
\]
 is cocontinuous. Indeed, if $(A,I, M_{\Spf(A)})$ is a perfect log prism in $\uX_{\prism}^{\perf, \opp}$ and $(A,I, M_{\Spf(A)})\to (B, IB, M_{\Spf(B)})$ is a flat cover in $\uX_{\prism}^{\opp}$. Let $(B^{\perf}, IB^{\perf})$ denote the perfection of the prism $(B,IB)$, and  $M_{\Spf(B^{\perf})}$ be the natural $\delta_{\log}$-structure induced from $M_{\Spf(B)}$. Then the composite $(A,I, M_{\Spf(A)})\to (B, IB, M_{\Spf(B)})\to (B^{\perf}, IB^{\perf}, M_{\Spf(B^{\perf})})$
 is a flat cover. 
 The inclusion functor  induces thus a morphism of topoi 
 \begin{equation}\label{E:perfect-to-nonperfect}
 \varepsilon\colon \quad  \Shv(X_{\prism}^{\perf})\simeq  \Shv(\uX^{\perf}_{\prism})\to \Shv(\uX_{\prism}).
 \end{equation}
 where the first equivalence is given by Proposition~\ref{P:equiv-perfect-log-site}. 
 The pullback functor $\varepsilon^*$ is the restriction, and for  
 $\calF\in \Shv(X^{\perf}_{\prism})$ and an object $(A,I, M_{\Spf(A)})$ of $\uX_{\prism}$, we have 
 	$\varepsilon_*(\calF)((A,I, M_{\Spf(A)})=\calF(A^{\perf}, IA^{\perf})$.

	\begin{definition}\label{D:prismatic-Laurent-crystal}
		(1) \emph{A  Laurent prismatic cyrstal on $\underline X_{\prism}$} is a presheaf $\calF$ of $\laurent$-modules  such that
		\begin{itemize}
			\item[(a)] for each object $(B,J_B,M_{\Spf(B)})$ of $\underline X_{\prism}$, $\calF(B,J_B,M_{\Spf(B)})$ is a finite projective $B[\frac{1}{J_B}]^{\wedge}$-module,
			\item[(b)] and for each morphism $u\colon (B,J_B,M_{\Spf(B)})\to (C,J_C,M_{\Spf(C)})$ in $\uX^{\opp}$, the canonical map 
			\[
			\calF(B,J_B,M_{\Spf(B)})\otimes_{B[1/J_B]_p^{\wedge},u}C[1/J_C]_p^{\wedge}\xra{\sim} \calF(C,J_C,M_{\Spf(C)})
			\] 
			is an isomorphism.
		\end{itemize}  
		A  Laurent prismatic crystal is automatically a sheaf. We write $\Vect(\underline X_{
			\prism}, \laurent)$ for the category  of  Laurent prismatic crystals   on $\underline X_{\prism}$. 
		
		(2)
		A  \emph{  Laurent prismatic $F$-crystal} on $\underline{X}_{\prism}$ is an object $\calE\in \Vect(\underline X, \laurent)$ together an isomorphism 
		\[
		\phi_{\calE}: \phi^*_{\cO_{\prism}[1/\calI_{\prism}]^{\wedge}_p}\calE\xra{\sim}\calE.
		\] 
		We write   $\Vect(\underline X_{\prism}, \laurent)^{\phi=1}$ for  the category of    Laurent prismatic $F$-crystals on $\underline X_{\prism}$. 
		

	\end{definition}
	
	Similarly, one have evident similar categories $\Vect(X_{\prism}, \cO_{\prism}[1/\calI_{\prism}]^{\wedge}_p)$,  $\Vect(X^{\perf}_{\prism}, \cO_{\prism}[1/\calI_{\prism}]^{\wedge}_p)$ and $\Vect(\uX^{\perf}_{\prism}, \cO_{\prism}[1/\calI_{\prism}]^{\wedge}_p)$   for $X_{\prism}$,  and   the perfect prismatic sites $X^{\perf}_{\prism}$ and $\uX^{\perf}_{\prism}$. 
	Let $X_{\eta}$ denote the adic generic fiber of the formal scheme  $X$, and $\Loc(X_{\eta, \et}, \ZZ_p)$ be the category of \'etale $\ZZ_p$-local systems on $X_{\eta}$. 
	
	\begin{theorem}[Bhatt--Scholze, Min--Wang]\label{T:etale-realization}
	There exist natural equivalences of categories 
	\[
\Vect(\uX_{\prism}, \cO_{\prism}[1/\calI_{\prism}]^{\wedge}_p)\simeq \Vect(\uX^{\perf}_{\prism}, \cO_{\prism}[1/\calI_{\prism}]^{\wedge}_p)\simeq \Vect(X^{\perf}_{\prism}, \cO_{\prism}[1/\calI_{\prism}]^{\wedge}_p)\simeq \Loc(X_{\eta, \et}, \ZZ_p).
	\]
	\end{theorem}
	
	\begin{proof}
	As mentioned in \cite[Remark 4.9]{MW}, the first equivalence follows from the same arguments as \cite[Thm. 3.2]{MW2} (see also \cite[Cor. 3.7]{BS2}). The second equivalence is a direct consequence of Proposition~\ref{P:equiv-perfect-log-site}, and the third is \cite[Cor. 3.8]{BS2}.
	
	\end{proof}
	
	\begin{remark}
	Koshikawa and Yao proved in \cite[Thm. 7.36]{KY} a similar equivalence  of categories for more general log formal schemes using  quasi-pro-Kummer-\'etale site of a log diamond. In our  case, their theorem should be equivalent to this theorem, since the local structure on $\uX$ is trivial on the generic fiber. 
		\end{remark}

	\if false

	\subsection{Functoriality} Let $\underline Y=(Y, M_Y)\to \underline X=(X,M_X)$ be a morphism of $p$-adic log formal schemes. Then  $\underline Y_{\prism}$ is naturally identified with the sliced category $ \underline{X}_{\prism}/_{\underline Y}$. Therefore, there exists a canonical morphism of  ringed topoi 
	\[
	f_{\prism}\colon 	\Shv(\underline Y, \cO_{\underline Y_{\prism}}) \to \Shv(\underline X, \cO_{\underline X_{\prism}})
	\]
	such that the inverse image functor  $f^{*}_{\prism}$ is the natural restriction, and the direct image of  a sheaf $\calF$ on $\underline Y_{\prism}$ is given as follows: For an object $(B, J, M_B)$  of $\underline X_{\prism}$, if   $\underline Y_{B/J}:=\underline Y\times_{\underline X}(\Spf(B/J), M_B)^a$ denotes the fiber product  in the category of log formal schemes, then one has 
	\[
	f_{\prism,*}\calF(B,J,M_B)=\Gamma((\underline Y_{B/J}/B)_{\prism}, \calF|_{(\underline Y_{B/J}/B)_{\prism}}),
	\]
	where $\calF|_{(\underline Y_{B/J}/(B,M_B))_{\prism}}$ is the natural restriction of $\calF$ to the relative log  prismatic site $(\underline Y_{B/J}/(B,M_B))_{\prism}$.  
	Note that $f_{\prism}^*$ induces the natural restriction functors:
	\[
	f^*_{\prism}\colon \CRhat^{\anfr}(\underline X_{\prism}, \cO_{\underline X_{\prism}} )\to \CRhat^{\anfr}(\underline Y_{\prism}, \cO_{\underline Y_{\prism}})
	\]
	and similarly the pull-back functors  for complete prismatic crystals or prismatic Laurent crystals.
	\fi

\section{Preliminaries on derivations and connections}\label{Section2}

In this section,  we recall some  preliminaries  in \cite{Tsuji} on the $\alpha$-derivation and modules with integrable connection. These results will  be needed in the next two sections.

\subsection{Derivations and $\ZZ_p$-actions}\label{S:derivation}
Let $A$ be a topological $\ZZ_p$-algebra, and $\alpha\in A$ be an element such that $A$ is $(p,\alpha)$-adically complete.
For an $A$-algebra $B$,  an $A$-linear \emph{$\alpha$-derivation} on $B$ is an  endomorphism $\partial\in \End_A(B)$ such that 
\[
\partial(fg)=\partial(f)g+f\partial(g)+\alpha \partial(f)\partial(g), \quad (f,g\in B).
\]

\begin{lemma}\label{L:derivation-gamma-action}
	Let $B$ be a $(p,\alpha)$-adically complete $A$-algebra.
	\begin{enumerate}
		\item[(1)]	If $\partial$ is an $A$-linear $\alpha$-derivation on $B$, then $\gamma:=\id_B+\alpha\partial$
		is an automoprhism of the $A$-algebra $B$ such that for every $c\in \ZZ_p$, the formula 
		\[
		\gamma^c(f):=\sum_{n=0}^{+\infty} \binom{c}{n} \alpha^n\partial^n(f)
		\]
		defines a continuous action of the group ${\ZZ_p}$ on the $A$-algebra $B$.	
		\item[(2)] Conversely, if $\alpha$ is a nonzero divisor in $B$ and  $\gamma$ is an automoprhism of the $A$-algebra $B$ such that $\gamma(f)\equiv f\mod \alpha B$, then 
		\[
		\partial(f):=\frac{\gamma(f)-f}{\alpha}
		\]
		defines an $A$-linear $\alpha$-derivation on $B$. 	
	\end{enumerate} 
\end{lemma}
\begin{proof}
	This is straightforward from the definition of $A$-linear $\alpha$-derivation. 
\end{proof}

Let  $B$  be a $(p,\alpha)$-adically complete $A$-algebra equipped with an $A$-linear derivation $\partial$, and  write  $\gamma=\id_B+\alpha\partial$ for the corresponding automorphism.
Let $M$ be a $B$-module.
An \emph{$\alpha$-derivation on $M$ with respect to $\partial$} is an endomorphism $\theta_M\in \End_A(M)$ such that 
\begin{equation}\label{E:alpha-derivation}
	\theta_{M}(bx)=\partial(b)x+b\theta_M(x)+\alpha\partial(b) \theta_M(x)=\partial(b)x+\gamma(b)\theta_M(x)
\end{equation}
for all $b\in B$ and $x\in M$.
Clearly, Lemma~\ref{L:derivation-gamma-action} can be generalized  to the $B$-module $M$: 
\begin{itemize}
	\item if $\theta_M$ is an $\alpha$-derivation on $M$ over $\partial$, then 
	\[
	\gamma_M\colon x\mapsto x+\alpha\theta_M(x)
	\]
	is an automorphism of $M$ such that $\gamma_M(b x)=\gamma(b)\gamma_M(x)$ for $b\in B$ and $x\in M$. Moreover, if $M$ is $(p,\alpha)$-adically complete, then  $c\mapsto \gamma_M^c$ defines a semi-linear action of $\ZZ_p$ on the $B$-module $M$.

	\item if $\alpha$ is a nonzero divisor $M$ and $\gamma_M$ is an automorphism on $M$ with $\gamma_M(x)\equiv x\mod \alpha M $ and $\gamma_{M}(bx)=\gamma(b)\gamma_M(x)$ for all $b\in B$ and $x\in M$, then 
	\[
	\theta_M(x):=\frac{\gamma_M(x)-x}{\alpha}
	\]
	defines an $\alpha$-derivation on $M$ with respect to  $\partial$.
\end{itemize}

\subsection{A noncommutative differential graded algebra}\label{S:derivation-dga}
Let $I$ be a finite set, and  $B$ be an $(p,\alpha)$-adically complete $A$-algebra equipped with $A$-linear $\alpha$-derivations $\partial_{i}$ for every $i\in I$ such that $\partial_i\partial_j=\partial_j\partial_i$ for all $ i,j\in I$.  We write $\gamma_i=\id_B+\alpha\partial_i$ with $i\in I$. It is easy to see that $\gamma_i\gamma_j=\gamma_j\gamma_i$ for all $i,j\in I$, and   Lemma~\ref{L:derivation-gamma-action} implies that there exists a natural continuous action of  the group $\ZZ_p^I=\prod_{i\in I}\gamma_i^{\ZZ_p}$ on $B$. 

We recall a noncommutative  differential graded algebra introduced  in \cite[\S 8]{Tsuji}. Let $\Omega_{B,\underline \gamma}=\bigoplus_{i\in I}A\omega_{i}$ with $B$-bimodule structure given by 
\[
b(a\omega_i)c=ba\gamma_i(c)\omega_i,
\] 
for $a,b,c\in B$ and $i\in I$. Define an $A$-linear map $d\colon B\to \Omega_{B, \underline \gamma}$ by 
\[d(a)=\sum_{i\in I}\partial_{i}(a)\omega_i.
\] Then $\rd$ satisfies the Lebnitz rule $d(ab)=ad(b)+d(a)b$. We define $\Omega^{\bullet}_{B,\underline \gamma}$ to be the graded left $B$-module with $\Omega^{q}_{B,\underline \gamma}=\bigwedge^q\Omega_{B,\underline \gamma}$ as left $B$-module. For $q\in \NN$ and $J=(i_1,\cdots, i_q)\in I^{q}$, we put 
\[
\omega_{J}:=\omega_{i_1}\wedge \cdots \wedge\omega_{i_q}\in \Omega^q_{B,\underline \gamma}, \quad \gamma_J:=\prod_{k=1}^q\gamma_{i_k}\in \ZZ_p^I.
\]
Note that $\omega_J\neq0$ if and only if the $i_k$'s are distinct. 
We equip $\Omega^{\bullet}_{B, \underline \gamma}$ with a structure of associative graded $A$-algebra by defining 
\[
a\omega_{J}\wedge a'\omega_{J'}:=a\gamma_{J}(a')\omega_{J}\wedge \omega_{J'}	
\]
for $a,a'\in I$, $J\in I^q$ and $J'\in I^{q'}$ for some $q,q'\in \NN$. For $q\in \NN$, we also define a differential map  $d^q\colon \Omega^q_{B,\underline \gamma}\to \Omega^{q+1}_{B,\underline \gamma}$ to be 
\[
d^q(a\omega_J)=\sum_{i\in I} \partial_i(a)\omega_i\wedge \omega_{J}, \quad (a\in B, J\in I^q).
\]
Then we have $d^0=d$, and $(\Omega^{\bullet}_{B,\underline \gamma}, d^{\bullet})$ becomes a differential graded $A$-algebra, i.e. every $d^q$ is $A$-linear, and we have $d^{q+1}\circ d^q=0$ and \[d^{q+q'}(\omega\wedge \eta)=d^{q}\omega \cdot \eta+(-1)^q\omega\wedge d^{q'}\eta\]
for all $\omega\in \Omega^q_{B,\underline \gamma}$ and $\eta\in \Omega^{q'}_{B,\underline \gamma}$.

\begin{definition}[\cite{Tsuji}, Def. 8.8, 8.9]\label{D:connection}
	(1)	A \emph{module with connection   over $(B, d:B\to \Omega_{B, \underline \gamma})$ (or  over $(B, \Omega^{\bullet}_{B, \underline \gamma})$)} is a $B$-module $M$ equipped with an $A$-linear  map 
	\[
	\nabla_M\colon M\to M\otimes_{B}\Omega_{B,\underline\gamma}
	\]
	such that $\nabla_{M}(am)=\nabla_M(m)a+m\otimes d(a)$ for $a\in B$ and $m\in M$. Let   $\theta_{M,i}\in \End_A(M)$ for $i\in I$ be endomorphisms of $M$ such that 
	$\nabla_{M}(m)=\sum_{i\in I}\theta_{M,i}(m)\omega_i$. 
	Then the Lebnitz rule for $\nabla_M$ is equivalent to saying that each $\theta_{M,i}$ is an $\alpha$-derivation on $M$ with respect to $\partial_{i}$. Write $\gamma_{M,i}:=\id_M+\alpha\theta_{M,i}$ be the autormphism of $M$ corresponding to $\theta_{M,i}$. Then $\gamma_{M,i}$ is $\gamma_i$-semilinear. If $M$ is complete for the $(p,\alpha)$-adic topology, we also require that $\nabla_M$ is continuous for the $(p,\alpha)$-adic topology.
	Morphisms between modules with (continuous) connection over $(B,d:B\to \Omega_{B, \underline \gamma})$ are defined in an evident way. 
	
	(2) Let $(M,\nabla_M)$ be a module with connection over $(B,\Omega_{B, \underline \gamma})$. For $q\in \NN$, let \[\nabla_{M}^{q}\colon M\otimes_{B}\otimes\Omega^{q}_{B,\underline \gamma}\to M\otimes_{B}\otimes\Omega^{q+1}_{B,\underline \gamma}\]
	be the additive map  
	\[\nabla_{M}^q(m\otimes \omega)=\nabla_M(m)\wedge \omega+ m\otimes d_q(\omega), \quad (m\in M,\omega\in \Omega^q_{B,\underline \gamma}).\]
	We say that $(M,\nabla_M)$ is \emph{integrable} if $\nabla_M^1\circ \nabla_M=0$. Note that this is equivalent to  $\theta_{M,i}\circ\theta_{M,j}=\theta_{M,j}\circ\theta_{M,i}$ for all $i,j\in I$. Hence, we have  $\nabla^{q+1}_M\circ\nabla_{M}^q=0$ for all $q\in \NN$. We call the resulting  complex $(M\otimes_{B}\Omega^{\bullet}_{B, \underline \gamma}, \nabla_{M}^{\bullet})$, the \emph{de Rham complex} of $(M,\nabla_{M})$.

	We denote by $\MIC^{\wedge}(B,d)$ the category of $(p,\alpha)$-adically complete  $B$-modules with continuous integrable connection over $(B, d:B\to \Omega_{B, \underline \gamma})$. 
\end{definition}

\begin{proposition}\label{P:module with connection}
	Let $M$ be a  $(p,\alpha)$-adically complete $B$-module. Then:
	\begin{enumerate}
		\item 
		If  $\nabla_M\colon M\to M\otimes_B\Omega_{B, \underline \gamma}$ is a continuous integrable connection over $(B,\Omega^{\bullet}_{B, \underline{\gamma}})$,  then the automorphisms  $(\gamma_{M,i})_{i\in I}$  extend naturally to a semi-linear continuous action by the group $\ZZ_p^I$ on $M$ that is trivial modulo $\alpha$. 
		
		\item Conversely, if $M$ is  $\alpha$-torsion free and  equipped with a semi-linear continuous action by $\ZZ_p^I=\prod_{i\in I}\gamma_i^{\ZZ_p}$ such that $(\gamma-1)(M)\subset \alpha M$ for all $\gamma\in \ZZ_p^I$, then $\nabla_M(m):=\sum_{i\in I} \theta_{M,i}(m)\omega_i$ with 
		\[
		\theta_{M,i}(m):=\frac{\gamma_{i}(m)-m}{\alpha}
		\]
		for $m\in M$
		well defines an integrable connection on $M$ over $(B,d)$. Moreover, one has a canonical isomorphism in the derived category of $A$-modules:
		\[(M\otimes_{B}\Omega^{\bullet}_{B, \underline \gamma}, \nabla_{M}^{\bullet})\simeq L\eta_{\alpha}R\Gamma(\ZZ_p^{I}, M)\]
		where $L\eta_{\alpha}$ is the d\'ecalage operator studied in \cite[\S 6]{BMS1}.
	\end{enumerate}	
	
\end{proposition}
\begin{proof}(1) The integrability of $\nabla_{M}$ implies the commutativity of the $\gamma_{M,i}$'s. The  rest of (1) follows easily from the discussion below Lemma~\ref{L:derivation-gamma-action}.
	
	(2) The first part of the statement follows from the inversement of the arguments in (1). Since $R\Gamma(\ZZ_p^{I}, M)$  is computed by the Koszul complex (cf. \cite[Lemma 7.3]{BMS1})
	\[
	\Kos(M, (\gamma_i-1)_{i\in I}):=\big( M\xra{(\gamma_i-1)_{i\in I}}\bigoplus_{i\in I}M\omega_i\to \cdots \xra{} M\omega_I\big)
	\]
	and one has  an isomoprhism of complexes 
	\[
	(M\otimes_B\Omega_{B, \underline \gamma}, \nabla^{\bullet}_{M})\simeq \eta_{\alpha} \Kos(M, (\gamma_i-1)_{i\in I}),
	\]
	the second part follows immediately.
	
\end{proof}

\subsection{Scalar extensions of modules with connection }\label{S:scalar-extension-MIC}
Let $C$ be another $(p,\alpha)$-adically complete $A$-algebra equipped with $\alpha$-derivations $(\partial_{C,i})_{i\in I}$ such that $\partial_{C,i}\circ \partial_{C,j}=\partial_{C,j}\circ \partial_{C,i}$. We write  $(\Omega^{\bullet}_{C, \underline \gamma}, d_C^{\bullet})$ for  the associated differential graded algebra. 

Suppose that we have a continuous map of topological $A$-algebras $f: B\to C$ compatible with the $\alpha$-derivations, i.e. $f\circ \partial_i=\partial_{C,i}\circ f$ for all $i\in I$. Clearly, $f$ commutes with the $\ZZ_p^I$-actions on $B$ and $C$ determined by their $\alpha$-derivatives via Proposition~\ref{P:module with connection}, and it induces  a morphism of differential graded $A$-algebras $df\colon \Omega^\bullet_{B,\underline \gamma}\to \Omega_{C,\underline \gamma}^{
	\bullet} $.

Let $(M, \nabla_{M})$ be an object of $\MIC^{\wedge}(B,d)$. We put $M_C:=M\widehat{\otimes}_{B}C$, and define $\nabla_{M_C}\colon M_C\to M_C\otimes_{C}\Omega_{C,\underline \gamma}$ by 
$$
\nabla_{M_C}(m\otimes c)=\nabla_M(m)\otimes c+m\otimes d_C(c), \quad \forall m\in M, c\in C.$$ 
Then $\nabla_{M_C}$ is a continuous  connection on $M_C$. In terms of $\alpha$-derivations, one has 
\[
\theta_{M_C, i}(m\otimes c)=\theta_{M,i}(m)\otimes \gamma_i(c)+m\otimes \partial_{C,i}(c), \quad \forall i\in I.
\]  
It is easy to see that $\theta_{M_C,i}$ commutes with $\theta_{M_C, j}$ for all $1\leqslant i,j\leqslant d$, hence $\nabla_{M_C}$ is integrable. 
The  $\ZZ_p^I$-action on $M_C$ given by $\nabla_{M_C}$ is given by 
\[
\gamma_i(m\otimes c)=\gamma_i(m)\otimes\gamma_i(c), \quad \forall i\in I,
\]
and the natural map $M\to M_C$ induces a morphism of complexes $(M, \nabla_M^{\bullet})\to (M_C, \nabla_{M_C}^{\bullet})$.
The construction $(M,\nabla_{M})\mapsto (M_C,\nabla_{M_C})$ is clearly functorial, thus it defines a scalar extension  functor 
\[
f^*\colon 	\MIC^{\wedge}(B, d)\to \MIC^{\wedge}(C,d_C).
\]
Whenever there is another morphism $g:C\to D$ such that $g^{*}\colon \MIC^{\wedge}(C,d_C)\to \MIC^\wedge(D, d_C)$ is defined, we have a natural equivalence of functors $g^*\circ f^*\simeq (g\circ f)^*$.

\section{Local computation of cohomology of log prismatic crystals}\label{S:cohomology-log-crystals}\label{S:local computation}

In this section, we will prove some local results  for the cohomology of complete prismatic crystals on a semistable formal scheme over a complete algebraically closed non-archimedean field.

\subsection{Notation}\label{S:setup of notation}
In this section, we will take $L=C$ to be an algebraically closed   completed non-archimedean field over $\QQ_p$. Let  $\gothm_C\subseteq \cO_C$ be the maximal ideal of the  ring of  integers $\cO_C$. 

Let $A_{\inf}=\Ainf(\cO_C)$, and let $\phi$ denote the Frobenius endomorphism on $A_{\inf}$.
 Fix a compatible system of primitive $p^n$-th root of unity  $\zeta_{p^n}\in C$ for $n\geqslant1$. We put $\epsilon =(1,\zeta_p, \zeta_{p^2}, \cdots)\in \cO_C^{\flat}$, $q=[\epsilon]$,  $\mu=q-1\in A_{\inf}$. Write $\phi$ for the Frobenius endomorphism on $A_{\inf}$, and put $$\xi=\frac{\mu}{\phi^{-1}(\mu)}, \quad [p]_q:=\phi(\xi)=\frac{q^p-1}{q-1}.$$

As in \S~\ref{S:semistable-setup}, we have the log formal scheme $\Spf(\uO_C)=(\Spf(\cO_C), M_{\cO^{\flat}_C})^a$ attached to $\Spf(\cO_C)$. 
Consider the prelog structure
\[\alpha_{A_{\inf}}\colon M_{\cO_C^{\flat}}:=\cO_{C}^{\flat}\backslash \{0\}\to A_{\inf}, \quad x\mapsto [x^p].
\]
and the associated log prism $\tAinf(\underline {\cO}_{C}):=(A_{\inf}, \ker(\theta_{\cO_C}\circ \phi^{-1})=([p]_q),   M_{\cO_C^{\flat}})^a$.
Then $\tAinf(\underline{\cO}_C)$ is an initial object of the absolute prismatic site $\Spf(\uO_C)^{\opp}_{\prism}$.


\subsection{Local setup}	Let $X=\Spf(R)$ be an affine small semistable $p$-adic formal scheme  over $\Spf(\cO_C)$. We fix  a framing  as in \eqref{E:small-framing}, i.e. a $p$-completely \'etale map of $\cO_C$-algebras
\begin{equation}\label{E:framing-over-C}
	\square\colon 	R^{\square}:=\cO_{C}\langle T_0, \cdots, T_{r}, T_{r+1}^{\pm1}, \cdots, T_d^{\pm 1}\rangle/(T_0\cdots T_r-\pi)  \to R
	\end{equation}
	for some integer $0\leqslant r\leqslant d$ and  $\pi\in \cO_C$ with $0<v_p(\pi)\leqslant 1$. Fix a compatible system $ (\pi^{1/p^n})_{n\geqslant0}$ of $p^n$-th roots of $\pi$  and put $\pi^{\flat}=(\pi, \pi^{1/p}, \cdots)\in \cO_C^{\flat}$. 
	Let $\uX$ be  the canonical log scheme  associated to $X$, and we have the chart 
\begin{equation}\label{E:local-chart-O_C}
\alpha_R\colon M_R=\cO_C^{\flat}\sqcup_{\NN}\NN^{r+1}\to R
\end{equation}
  as defined in \S~\ref{S:local-chart}.
Using the initial object $\tAinf(\uO_C)$ of $\Spf(\uO_C)_{\prism}^{\opp}$, we   identify $\uX_{\prism}$ with the relative prismatic site $(\uX/\tAinf(\uO_C))_{\prism}$. In particular, every object $(A,I, M_{\Spf(A)})$ of $\uX_{\prism}$ is equipped with a canonical map of log prisms $\tAinf(\uO_C)\to (A,I, M_{\Spf(A)})$ so that $I=([p]_q)$. Note also that, by Lemma~\ref{L:local objects}, there exists a $\delta_{\log}$-structure $\alpha_{A}, \delta_{\log}:M_R\to A$ such that $(A,([p]_q), M_R)$ is a prelog prism with  $(A,I, M_{\Spf(A)})=(A,([p]_q), M_R)^a$.

\subsection{The object $\urA(R)$}\label{S:frames}	
We consider the $\delta$-ring 
\[\rA(R^{\square}):=A_{\inf}\langle T_0, \dots, T_{r}, T_{r+1}^{\pm 1}, \dots,  T_{d}^{\pm 1}\rangle /(T_0\cdots T_r-[(\pi^ {\flat})^p])\]
with $\delta(T_i)=0$. Since $R$ is $p$-completely \'etale over $R^{\square}$, by deformation theory, there exists a $(p, [p]_q)$-completely \'etale $\rA(R^{\square})$-algebra $\rA(R)$ such that $\rA(R)/([p]_q)\simeq R$. By \cite[Lemma 2.18]{BS}, the $\delta$-structure on  $\rA(R^{\square})$ lifts uniquely to $\rA(R)$. Moreover, we endow $\rA(R)$ with  the prelog structure
\[
\alpha_{\rA(R)}\colon  M_{R}=M_{\cO^{\flat}_C} \sqcup_{\NN}\NN^{r+1}\to \rA(R^{\square})\to  \rA(R)\]
where 
\begin{itemize}
	\item $M_{\cO^{\flat}_C}\to \rA(R)$ is given by $x\mapsto[x^p]\in A_{\inf}\subset  \rA(R)$, 
	\item $\NN^{r+1}=\bigoplus_{i=0}^r \NN e_i\to \rA(R)$ sends $e_i$ to $T_i$ for $0\leqslant i\leqslant r$. 
\end{itemize}
We define $\delta_{\log}\colon M_{\rA(R)}\to \rA(R)$  to be the zero map, and we get a prelog prism $(\rA(R), ([p]_q), M_{R})$.  Let  $\underline {\rA}(R):=(\rA(R), ([p]_q ), M_{R})^a$ be its associated log prism. 
Together with   the canonical isomorphism $$
(\Spf(\rA(R)/([p]_q)), M_R\xra{\alpha_{\rA(R)}}\rA(R)\to \rA(R)/([p]_q))^a\xra{\sim}\uX,$$
$\underline {\rA}(R)$ becomes an object of $\underline X_{\prism}$.  

Consider the  action of $\Gamma:=\bigoplus_{i=1}^d \ZZ_p\gamma_i$ on $\rA(R^{\square})$ given by 
\begin{equation}\label{E:Gamma-action-Ainf(R)}
\gamma_i(T_j)=\begin{cases}
	q^pT_i & \text{if } j=i,\\
	q^{-p}T_0 & \text{if $j=0$ and $1\leqslant i\leqslant r$},\\
	T_j &\text{otherwise}. 
\end{cases}
\end{equation}
The action is trivial modulo $q^p-1$ and compatible with the $\delta$-structure in the sense that $\delta( \gamma_i(x))=\gamma_i(\delta(x))$ for $x\in \rA(R^{\square})$.
By the $(p,[p]_q)$-complete \'etaleness of $\rA(R^{\square})\to \rA(R)$, there exists a unique extension  of  the  $\Gamma$-action  to  the $\delta$-ring $\rA(R)$  that is trivial modulo $q^p-1$. Note also that the action preserves the log structure on $\rA(R)$ attached to prelog log structure $\alpha_{\rA(R)}$ compatible with the $\delta_{\log}$-structure. Therefore,   $\Gamma$ acts on the object $\underline \rA(R)$  of $\uX_{\prism}$. For $f\in \rA(R)$ and $1\leqslant i\leqslant d$, we put 
\[
\partial_i(f):=\frac{\gamma_i(f)-f}{\mu}.
\]
Then $\partial_i(f)\in [p]_q \rA(R)$, and   $\partial_i$ is an $\Ainf$-linear $\mu$-derivation on $\rA(R)$ in the sense of \S \ref{S:derivation} by Lemma~\ref{L:derivation-gamma-action}.

\begin{lemma}\label{L:cover-final-object}
	The object $\urA(R)$ is a strong covering object of $\uX_{\prism}$ (Definition~\ref{D:versal object}), i.e. for every object $(B,([p]_q),M_{\Spf(B)})$ of $\uX_{\prism}$, the coproduct 
	 $$(D_B, ([p]_q), M_{\Spf(D_B)}):=(B, ([p]_q ), M_{\Spf(B)})\coprod \urA(R)$$  exists  in the category $\underline {X}_{\prism}^{\opp}$, and the canonical map $(B, ([p]_q), M_{\Spf(B)})\to (D_B, ([p]_q), M_{\Spf(D_B)})$ is a flat cover. Moreover,  the formation of $(D_B, ([p]_q), M_{\Spf(D_B)})$ commutes with base change in $(B, ([p]_q), M_{\Spf(B)})$.  
	
	
\end{lemma}
\begin{proof}
	
	The proof is similar to \cite[Lemma 4.2]{Tian}. 
	By Lemma~\ref{L:local objects}, there exists a $\delta_{\log}$-structure $\alpha_{B}, \delta_{\log,B}\colon M_R\to B$ such that $(B, ([p]_q), M_{R})^a=(B, ([p]_q), M_{\Spf(B)})$. 
	 Let $(D, \delta, \alpha_{D}\colon M_D\to D, \delta_{\log,D})$ be the pushout of the maps $(A_{\inf}, M_{\cO^{\flat}_C})\to (B, M_R)$ and $(A_{\inf}, M_{\cO^{\flat}_C})\to (\rA(R), M_R)$ in the category of $\delta_{\log}$-rings. In particular, we have $D=B\widehat\otimes_{A_{\inf}}\rA(R)$, the $(p,[p]_q)$-completed tensor product of $B$  and $\rA(R)$, and $\alpha_{D}=\alpha_B\sqcup \alpha_{\rA(R)}\colon  M_{D}=M_R\sqcup_{M_{\cO_C^{\flat}}}M_R\to D$. 
	 
Put $\overline B=B/([p]_q)$, and $M_{\Spf(\overline B)}$ be the restriction of $M_{\Spf(B)}$ to $\Spf(\overline B)$.  Then  $M_{\Spf(\overline B)}$ is induced by  the composite map $\alpha_{\overline B}\colon M_R\xra{\alpha_B}B\to \overline B$, which  coincides with $M_R\xra{\alpha_R}R\to \overline B$ by Lemma~\ref{L:local objects}. Let $(\widetilde D, M_{\widetilde D})$ be the exactification of  the surjection of prelog rings 
	$
	(D,M_D)\to (\overline B ,  M_{R}\xra{\alpha_{ \overline B}}  \overline B)
	$
	induced by the canonical surjection $(B,M_R)\to (\overline B, M_{\overline B})$ and 
	\[
	(\rA(R),M_R)\twoheadrightarrow  (R,M_R)\rightarrow (\overline B, M_{R}).
	\]
	Note that there exists an induced $\delta_{\log}$-structure on $\widetilde D$. 
	Let $\widetilde J\subset \widetilde D$ be the ideal of $\widetilde D\twoheadrightarrow \overline B$. We define  
	\[
	D_B:=\widetilde D\{\frac{\widetilde J}{[p]_q}\}^{\wedge}
	\]
	to be  the prismatic envelope of $\widetilde D$ with respect to $\widetilde J$, and equip it with the  prelog structure $M_{\widetilde D}\to \widetilde D\to D_B$ and the induced $\delta_{\log}$-structure. Then the associated log prism $(D_B,([p]_q), M_{\Spf(D_B)}):=(D_B,([p]_q), M_{\widetilde D})^a$ together with the composite   map $$(\Spf(D_B/([p]_q)), M_{\Spf(D_B)})\to  (\Spf(\overline B), M_{\Spf(\overline B}))\to \uX$$ gives  the coproduct $(B, ([p]_q), M_{\Spf(B)})\coprod \underline{\rA}(R)$  in $\underline X_{\prism}^{\opp}$. It remains to show that $(B, ([p]_q), M_{\Spf(B)})\to (D_B, ([p]_q), M_{\Spf(D_B)})$ is a flat cover and the formation of $(D_B, ([p]_q), M_{\Spf(D_B)})$ commutes with base change in $(B, ([p]_q), M_{\Spf(B)})$.

	Note that the prelog structure $\alpha_{R}\colon M_R\to \rA(R)$ factors through $\rA(R^{\square})$. Consider the bounded log prism  $ \urA(R^{\square})=(\rA(R^{\square}), ([p]_	q), M_R)^a$. If we apply the construction above with $R$ replaced   by $R^{\square}$, we obtain similarly $\widetilde D^{\square}$,  $\widetilde J^{\square}\subset \widetilde D^{\square}$ and $D_B^{\square}$ in the evident sense.
	Moreover, one has a morphism of $\delta$-pairs $ (\widetilde D^{\square}, \widetilde J^{\square})\to (\widetilde D, \widetilde J)$ such that 
	\begin{itemize}
		\item  the underlying map  $\widetilde D^{\square}\to \widetilde D$ is $(p, [p]_q)$-completely \'etale,
		\item it induces isomorphisms
		$
		\widetilde D^{\square}/\widetilde J^{\square}\simeq \widetilde D/\widetilde J\simeq \overline B.
		$
	\end{itemize}
	By \cite[Prop. 4.11]{Tsuji}, the induced map of prisms $(D_B^{\square}, ([p]_q))\to (D_B, ([p]_q))$ over $(B,([p]_q))$ is an isomoprhism. Therefore, in order to prove the Lemma, we may reduce to the case $R=R^{\square}$. 
	
	We fix a morphism of log $(p, [p]_q)$-adic formal schemes 
	\[
	\tilde f\colon (\Spf(B),M_{\Spf(B)})\to (\Spf(\rA(R)), M_{R})^a
		\]
	over $(\Spf(A_{\inf}),M_{\cO^{\flat}_C}\xra{x\mapsto [x^p]}A_{\inf})^a$
	lifting the structure  map $(\Spf(\overline B), M_{\Spf(\overline B)})\to\uX$. The existence of such an $\tilde f$ is ensured by the formal smoothness of $(\Spf(\rA(R)), M_{R})^a$ over $(\Spf(A_{\inf}),M_{\cO^{\flat}_C}\xra{x\mapsto [x^p]}A_{\inf})^a$.
	
	Let $T_i\in D$ with $0\leqslant i\leqslant d$ still denote the image of $T_i$ via the natural map $\rA(R)\to D$, and write $t_i\in C$ for the image of $T_i$ via $\rA(R)\xra{\tilde f^*} B \to D$. Then it is easy to see that 
	\[
	\widetilde {D}=D\langle U_i^{\pm 1}:1\leqslant i\leqslant d\rangle /(T_i-t_iU_i:1\leqslant i\leqslant d).\]
	Let $u_i\in \widetilde D$ with $1\leqslant i\leqslant d$ denote the image of $U_i$. Then one has $\widetilde J=(u_i-1:1\leqslant i\leqslant d)$ and 
	\begin{equation}\label{E:D_B-formula}
		D_B=\widetilde{D}\{S_i:1\leqslant i\leqslant d\}^{\wedge}/(u_i-1-[p]_qS_i: 1\leqslant i\leq d)_{\delta},
	\end{equation}
	where $\widetilde{D}\{S_i:1\leqslant i\leqslant d\}^{\wedge}$ means the  $(p,[p]_q)$-adically complete free $\delta$-algebra over $\widetilde D$,  and $(u_i-1-[p]_qS_i: 1\leqslant i\leqslant d)_{\delta}$ means the $\delta$-ideal generated by $u_i-1-[p]_qS_i$.
	Note that    $(u_i-1, 1\leqslant i\leqslant d)$ is  $(p, [p]_q)$-completely regular relative to $B$.  By \cite[Prop. 3.13]{BS}, $D_B$ is $(p,[p]_q)$-completely faithfully flat over $B$. This proves that $(B,([p]_q), M_{\Spf(B)})\to (D_B, ([p]_q), M_{\Spf(D_B)})$ is a flat cover in $\uX_{\prism}$. 
	
	Finally, note that the formation of both $D$ and $\widetilde D$ commutes with base change in $(B,([p]_q),M_{\Spf(B)})$. 
	By \cite[Prop. 3.13]{BS}, so does  the formation of $(D_B, ([p]_q), M_{\Spf(D_B)})$.

\end{proof}

We will give a more concrete description for $D_B$.

\begin{lemma}\label{L:basis-for-D_B}
Let $(B, ([p]_q), M_{\Spf(B)})$ be an object of $\uX_{\prism}$, $\overline B=B/([p]_q)$.	Fix a lifting $\tilde f\colon (\Spf(B),M_{\Spf(B)})\to (\Spf(\rA(R)), M_{R})^a$ of  the structural map $(\Spf(\overline B), M_{\Spf(\overline B)})\to \uX$, and write $\tau_i$ with $1\leqslant i\leqslant d$ for the image   in $D_B$ of the element $S_i$ as in \eqref{E:D_B-formula}. 
	For $m=\sum_{k=0}^{l}a_kp^k\in \NN$ with $a_i\in \NN\cap [0,p-1]$,  define 
	\[
	\tau_i^{\{m\}_{\delta}}:=\prod_{k=0}^l \delta^{k}(\tau_i)^{a_k}.
	\]
	Then for any ideal $J\subset \Ainf$ containing some power of $(p,[p]_q)$, $D_B/JD_B$ is a free $B/JB$-module with basis $\{\prod_{i=1}^d \tau_i^{\{m_i\}_{\delta}}: (m_i)\in \NN^d\}$. 
\end{lemma}
\begin{proof}	
	As in the proof of Lemma~\ref{L:cover-final-object}(2), we may assume that $R^{\square}=R$. Then one has 
	\[
	D=\rA(R^{\square})\widehat\otimes_{\Ainf}B\simeq B\langle T_0, \cdots, T_d\rangle/(T_0\cdots T_r-[\pi^{\flat, p}]),
	\]
	and
	$\widetilde D=D\langle U_i^{\pm 1}:1\leqslant i\leqslant d\rangle /(T_i-a_iU_i:1\leqslant i\leqslant d)$ in the notation of the proof of Lemma~\ref{L:cover-final-object}(2).
	It follows easily that, for any ideal $J\subset \Ainf$ as in the statement, the canonical map 
	\[B/JB\simeq \widetilde D/\big(J\widetilde D+\sum_{i=1}^d(u_i-1)\widetilde D\big)
	\]
	is an ismorphism. Then our Lemma follows directly  from   \cite[Prop. 2.1 (4)]{Tsuji}.
\end{proof}

\begin{lemma}\label{L:phi-action-D_B}
	In the situation of Lemma~\ref{L:basis-for-D_B}, assume moreover that $\delta(\tilde f^*(T_i))=0$ for all $1\leqslant i\leqslant d$. 
	Let $\phi\colon D_B\to D_B$ denote the Frobenius action on $D_B$, and $\calK\subset D_B$ denote the closure of the  ideal generated by $\delta^k(\tau_i)$ for $k\geqslant 0$ and $1\leqslant i\leqslant d$. Then we have $\phi(x)\in [p]_qD_B$ for every $x\in \calK$.
\end{lemma}
\begin{proof}
	The proof is similar to \cite[Lemma 5.3]{Tian}. For any integer $k\geqslant0$, let $\calJ_k\subset D_B$ denote the closure of the ideal generated by $([p]_q)^{p^j}\delta^j(\tau_{i})$ for $0\leqslant j\leqslant k$
	and $1\leqslant i\leqslant d$. We prove by induction on $k\geqslant 0$ that there exists $b_{i,k+1}\in D_B^{\times}$ such that 
	\[
	\phi(\delta^k(\tau_i))=\delta^k(\tau_i)^p+p\delta^{k+1}(\tau_i)\equiv b_{i,k+1}([p]_q)^{p^{k+1}}\delta^{k+1}(\tau_i) \mod \calJ_{k},
	\]
	which will certainly imply our Lemma.  We start with the case $k=0$.  The assumption $f(\delta(T_i))=0$ for $1\leqslant i\leqslant d$ implies that $\delta(u_i)=0$ for all $1\leqslant i\leqslant d$.  Applying  $\delta$ to the relation $u_i-1=[p]_q\tau_i$, we obtain 
	\begin{align*}
		\delta([p]_q)\phi(\tau_i)+([p]_q)^p \delta(\tau_i)&=\delta([p]_q\tau_i)\\
		&=\delta(u_i-1)=\frac{1}{p}\sum_{j=1}^{p-1}\binom{p}{j}(u_i-1)^j\in \calJ_0
	\end{align*}
	As $\delta([p]_q)\in A_{\inf}^{\times}$, the  claim for $k=0$ follows immediately. The rest of the induction process works exactly as in the proof of  \cite[Lemma~5.3]{Tian}.
\end{proof}

\begin{proposition}\label{P:dp-envelop-D_B}
	In the situation of Lemma~\ref{L:basis-for-D_B}, assume moreover that $\delta(\tilde f^*(T_i))=0$ for all $1\leqslant i\leqslant d$.   Let  $\bar B\langle X_1,\cdots, X_d\rangle^{PD}$ denote the $p$-adic completion of the free PD-evelope algebra in $d$-variables over $\overline B$.  Then there exists an isomorphism  of topological $\overline  B$-algebras
	\[
	\overline  B\langle S_1,\cdots, S_d\rangle^{PD}\xra{\sim} D_B/[p]_qD_B
	\]		
	defined  by  $S_i\mapsto  \tau_i$ for $1\leqslant i\leqslant d$, where $\tau_i$ still denotes its image in $D_B/[p]_qD_B$ by a slight abuse of notation.
\end{proposition}
\begin{proof}
	For $x\in \calK$, we put $\gamma_p(x):=-\frac{\delta(x)}{(p-1)!}\in \calK$. We claim that 
	\begin{enumerate}
		\item $p!\gamma_p(x)\equiv x^p\mod [p]_qD_B$ for all $x\in \calK$,
		\item $\gamma_p(ax)\equiv a^p \gamma_p(x)\mod [p]_qD_B$
		for all $x\in \calK$,
		
		\item $\gamma_p(x+y)\equiv \gamma_p(x)+\gamma_p(y)+\sum_{i=1}^{p-1}\frac{1}{i!(p-i)!}x^iy^{p-i}$ for all $x,y\in \calK$.
		
	\end{enumerate}
	Indeed, Lemma~\ref{L:phi-action-D_B} implies that $\phi(x)=x^p+p\delta(x)\equiv 0\mod [p]_qD_B$ for $x\in \calK$, from which  assertion (1) follows. Claim (2) follows from $\delta(ax)=a^p\delta(x)+\delta(a)\phi(x)$ and Lemma~\ref{L:phi-action-D_B}. Claim (3) is a direct consequence of the additive property for $\delta$.
	
	Write $\overline\calK\subset D_B/[p]_qD_B$ for the image of $\calK$.  	For $\bar x\in \overline\calK$, we define $\gamma_p(\bar x)\in \overline \calK$ to be the image of $\gamma_p(x)$ for a lift $x\in \calK$ of $\bar x$. 
	As $[p]_q\calK=[p]_qD_B\cap \calK$, it is easy to see that $\gamma_p(\bar x)$ is well defined. 
	By \cite[\href{https://stacks.math.columbia.edu/tag/07GS}{Tag 07GS}]{stacks-project}, there exists a unique PD-structure on $\overline \calK$ such that the $p$-th divided power map  is given by $\gamma_p$. Therefore, there exists a map 
	\[
	\Psi\colon \bar B\langle S_1,\cdots, S_d\rangle^{PD}\to  D_B/[p]_qD_B
	\]
	sending $S_i$ to $\tau_i$ for $1\leqslant i\leqslant d$. For each $m=(m_i)_i\in \NN^d$, write $\tau^{[m]}$ for the image of  $S^{[m]}:=\prod_{i=1}^dS_i^{[m_i]}$.
	Let $\tau^{\{m\}_{\delta}}\in D_B$ be  the element defined in Lemma~\ref{L:basis-for-D_B}, and we still use it to denote  its image in $D_{B}/[p]_qD_B$ by abuse of notation. 
	Then an easy calculation shows that
	\[
	\tau^{\{m\}_{\delta}}=c_m\tau^{[m]} \quad\text{with}\quad  c_m=
	\prod_{i=1}^d(-1)^{l_i(l_i+1)/2}\frac{m_i!}{p^{v_p(m_i!)}}\in \ZZ_{(p)}^{\times}
	\]
	where  $l_i\geqslant 0$ is the minimal integer such that $m_i<p^{l_i}$. 
	It follows immediately from Lemma~\ref{L:basis-for-D_B} that the morphism $\Psi$ above is an isomorphism. 
	
\end{proof}

\subsection{$\Gamma$-action and derivations}  
We keep the assumption and notation of Lemma~\ref{L:basis-for-D_B}. The $\Gamma$-action on the object $\urA(R)$ of $\uX_{\prism}$ induces  by functoriality a natural $B$-linear $\Gamma$-action on $(D_B, ([p]_q), M_{\Spf(D_B)})$. 
Let $u_i$ for $1\leqslant i\leqslant d$ be as in \eqref{E:D_B-formula}. One has $u_i=1+[p]_q\tau_i$ for $1\leqslant i\leqslant d$, $\delta(u_i)=0$ and 
\[
\gamma_j(u_i)=\begin{cases}q^p u_i &\text{if }j=i,\\
	u_i & \text{if }j\neq i.\end{cases}
\]
As $D_B$ has no $[p]_q$-torsion, we deudce that 
\begin{equation}\label{E:gamma-action-tau}
	\gamma_j(\tau_i)=\begin{cases}\tau_i+\mu u_i &\text{if }j=i,\\
		\tau_i & \text{if }j\neq i.\end{cases}
\end{equation}
Note that the $\Gamma$-action is compatible with the $\delta$-strucutre on $D_B$, i.e. we have $\gamma_i(\delta(x))=\delta(\gamma_i(x))$ for all $x\in D_B$ and $1\leqslant i\leqslant d$.

\begin{proposition}\label{P:derivation-D_B}
	Under the above notation, the following holds: 
	\begin{enumerate}
		\item[(1)]  There  exist  unique $B$-linear  $\mu$-derivations $\partial_{D_B,i}\colon D_B\to D_B$  (see \S~\ref{S:derivation} for this notion) with $1\leqslant i\leqslant d$ such that $\partial_{D_B, i}\circ \partial_{D_B, j}=\partial_{D_B, j}\circ \partial_{D_B, i}$,
		\begin{equation}\label{E:partial-tau}
			\partial_{D_B, i}(\tau_j)=\begin{cases}
				u_i & \text{if } j=i,\\
				0, & \text{if }j\neq i,
			\end{cases}
		\end{equation}
		for all $1\leqslant i,j\leqslant d$,
		and 
		\begin{equation}\label{E:partial-delta-x}
			\partial_{D_B,i}(\delta(x))=\beta\partial_{D_B,i}(x)^p+(\mu^{p-1}+p\beta)\delta(\partial_{D_B,i}(x))-\sum_{k=1}^{p-1}\frac{1}{p}\binom{p}{k}\mu^{k-1}\partial_{D_B,i}(x)^kx^{p-k} 
		\end{equation}  for all $x\in D_B$ and $1\leqslant i\leqslant d$, where  $\beta:=\frac{\delta(\mu)}{\mu}=\sum_{k=1}^{p-1}\frac{1}{p}\binom{p}{k}\mu^{k-1}\in A_{\inf}^{\times}$ is a unit.
		
		\item[(2)] For each $i$ with $1\leqslant i\leqslant d$, the  $\mu$-derivation $\partial_{D_B, i}$   on $D_B$ is compatible with  $\partial_{i}$ on  $\rA(R)$ via the natural map $p_A\colon \rA(R)\to D_B$, i.e. $\partial_{D_B, i}\circ p_A=p_A\circ \partial_i$.
		
		\item[(3)]  The 
		natural $\Gamma$-action on $D_B$ is given by $\gamma_i|_{D_B}=\id_{D_B}+\mu\partial_{D_B,i}$. In particular, the natural $\Gamma$-action on $D_B$ is trivial modulo $\mu$. 
		
	\end{enumerate}

\end{proposition}
\begin{proof}
	(1) The uniqueness of $\partial_{D_B, i}$'s is clear since $D_B$ is topologically generated over $B$ by $\delta^{k}(\tau_i)$ for $k\in \NN$ and $1\leqslant i\leqslant d$. It suffices to show their existence. 
	Lemma~\ref{L:basis-for-D_B} shows that the assumption of \cite[Prop. 7.1]{Tsuji} is satisfied. Applying \emph{loc. cit.}, we obtain $\mu u_i$-derivations $\theta_{D_B,i}$ with $1\leqslant i\leqslant d$ on $D_{B}$ such that a similar compatibility with the $\delta$-structure on $D_B$ is satisfied, and we have   
	\begin{equation}
		\theta_{D_B,i}(\tau_j)
		=\begin{cases}
			1 & \text{if }i=j,\\
			0 & \text{if }i\neq j\end{cases},
	\end{equation}
	and  $\theta_{D_B, i}\circ \theta_{D_B,j}=\theta_{D_B, j}\circ \theta_{D_B,i}$ for all $i,j$. Then it is easy to see that  $\partial_{D_B, i}:=u_i \theta_{D_B, i}$ with $1\leqslant i\leqslant d$ satisfy our requirement. 
	
	(2) It suffices to check that $\partial_{D_B, i}(p_A(T_j))=p_A(\partial_i(T_j))$ for all $i,j$. When $i\neq j$, then both sides are equal to $0$. For $i=j$, we have $p_A(T_i)=u_ia_i$ and 
	\[
	\partial_{D_B, i}(p_A(T_i))=\partial_{D_B,i}(u_ia_i)=[p]_qu_ia_i=p_A([p]_qT_i)=p_A(\partial_i(T_i)).
	\]

	(3) We  write $\gamma_i':=\id_{D_B}+\mu \partial_{D_B, i}$ and  check first  that  $\gamma_i'$ is compatible with $\delta$-structure, i.e. $\gamma_i'(\delta(x))=\delta(\gamma'_i(x))$ for all $x\in D_B$. In fact, we have 
	\begin{align*}
		\gamma_i'(\delta(x))&=\delta(x)+\mu\partial_{D_B,i}(\delta(x))\\
		&=\delta(x)+\mu\bigg( \beta\partial_{D_B,i}(x)^p+(\mu^{p-1}+p\beta)\delta(\partial_{D_B,i}(x))-\sum_{k=1}^{p-1}\frac{1}{p}\binom{p}{k}\mu^{k-1}\partial_{D_B,i}(x)^k x^{p-k}\bigg)\\
		&=\delta(x)+\delta(x\partial_{D_B,i}(x))-\sum_{k=1}^{p-1}\frac{1}{p}\binom{p}{k}\mu^k \partial_{D_B,i}(x)^kx^{p-k}\\
		&=\delta(x+\mu\partial_{D_B,i}(x))=\delta(\gamma'_i(x)),
	\end{align*} 
	where the second equality used \eqref{E:partial-delta-x}. Since $D_B$ is topologically generated as a $B$-algebra by $\{\delta^k(\tau_j): k\geqslant 0, 1\leqslant j\leqslant d\}$, it is enough to show that $\gamma_i(\tau_j)=\gamma'_i(\tau_j)$ for all $1\leqslant i,j\leqslant d$. This follows directly from \eqref{E:gamma-action-tau} and 
	\eqref{E:partial-tau}.

\end{proof}
\begin{remark}
	(1)
	If $\mu$ is a nonzero divisor in $D_B$ (or equivalently in $B$), then one can verify directly that the $\Gamma$-action on $D_B$ is trivial modulo $\mu$, and  define  
	\[
	\partial_{D_B,i}(x)=\frac{\gamma_{i}(x)-x}{\mu},
	\]
	for all  $ x\in D_B.$ (cf. Lemma~\ref{L:derivation-gamma-action} and Prop.~\ref{P:module with connection}). Then \eqref{E:gamma-action-tau} is equivalent to \eqref{E:partial-tau} and the formula \eqref{E:partial-delta-x} is equivalent to the compatibility of the $\Gamma$-action with the $\delta$-structure on $D_B$.
	
	(2) In order to have  a description \eqref{E:D_B-formula}, we have chosen a morphism $\tilde f\colon (\Spf(B),M_{\Spf(B)})\to (\Spf(\rA(R)), M_{R})^a$
 lifting  $(\Spf(\overline B), M_{\Spf(\overline B)})\to \uX$. But the $\mu$-derivations $\partial_{D_B, i}$ are independent of the choice of $\tilde f$. Indeed, assume $\tilde f'$ is another such  lifting. We write $a_i'\in D_{B}$ for the image of $T_i$ via $\rA(R)\xra{\tilde f'^*} B\to D_B$, and define $u_i'$ and $\tau'_i$ with $a_i$ replaced by $a_i'$. Then there exists $b_i\in B$ with $a_i=a'_i(1+[p]_qb_i)$,  $u_i'=u_i(1+[p]_qb_i)$ and $\tau'_i=\tau_i+b_iu_i$. Hence, one has 
	\[
	\partial_{D_B, i}(\tau'_i)=\partial_{D_B,i}(\tau_i)+b_i\partial_{D_B, i}(u_i)=u_i+b_i[p]_qu_i=u_i'.
	\]
 the previous remark. 
	
\end{remark}

Let   $\qOmega^1_{D_B/B}$ be the free $D_B$-module of rank $d$ with basis $\{\frac{d u_i}{u_i}: 1\leqslant i\leqslant d\}$. We equip $\qOmega^1_{D_B/B}$ with the $D_B$-bimodule structure given by 
\[
a \cdot (b\frac{du_i}{u_i})\cdot c:=ab\gamma_i(c) \frac{du_i}{u_i}, \quad (a,b,c\in D_B, 1\leqslant i\leqslant d).
\]
We can thus identify $\qOmega^1_{D_B/B}$ with $\Omega_{D_B, \underline \gamma}$ defined in \S~\ref{S:derivation-dga}. Applying the construction of \S~\ref{S:derivation-dga}, we get a differential graded algebra
$(\qOmega_{D_B/B}^{\bullet}, d_{D_B/B}^{\bullet})$ with $d_{D_{B}/B}^k\colon \qOmega^k_{D_B/B}:=\bigwedge^k\qOmega^1_{D_B/B}\to \qOmega^{k+1}_{D_B/B}$ given by 
\[
d_{D_B/B}^k(f \frac{du_{i_1}}{u_{i_1}}\wedge \cdots \wedge\frac{du_{i_k}}{u_{i_k}})= \sum_{i=1}^{d} \partial_{D_B, i}(f)\frac{du_i}{u_{i}}\wedge \frac{du_{i_1}}{u_{i_1}}\wedge \cdots \wedge\frac{du_{i_k}}{u_{i_k}}.
\]
We have the following analogue of Poincar\'e's Lemma in this setup.
\begin{proposition}\label{P:poincare-lemma}
	Let $M$ be a $(p, [p]_q)$-adically complete $B$-module. Then the natural map $M\to M\widehat{\otimes}_BD_B$ induces a quasi-isomorphism of complexes
	$$M\xra{\sim}(M\widehat\otimes_B\qOmega^{\bullet}_{D_{B}/B}, \id_M\widehat\otimes d^{\bullet}_{D_B/B}).$$ 
	Here, $\widehat{\otimes}$ denotes the $(p,[p]_q)$-adic completion of the usual tensor product. 
\end{proposition}
\begin{proof}
	If $M$ is annihilated by some power of $(p, [p]_q)$, then this is exactly  \cite[Theorem 7.8]{Tsuji} (Note that $\theta_{\widehat{D},i}$ in \emph{loc. cit.} corresponds to $\partial_{D_B,i}/u_i$, and all $u_i$ are invertible).
	The general case follows by applying the statement to $M/(p, [p]_q)^k$ for all $k\in  \NN$ and taking the projective limit.
\end{proof}

\subsection{$q$-de Rham complex and $q$-Higgs modules}

Consider the free    $\rA(R)$-module   
\[\qOmega^1_{\rA(R)}:= \bigoplus_{i=1}^d\urA(R)\frac{dT_i}{T_i}\]
equipped  
with a right $\rA(R)$-module structure by 
\[\frac{dT_i}{T_i} f=\gamma_i(f)\frac{dT_i}{T_i}, \quad (f\in \rA(R), 1\leqslant i\leqslant d).\] 
We  identify $\qOmega^1_{\rA(R)}$ with $\Omega_{\rA(R), \underline\gamma}$ considered in \S \ref{S:derivation-dga}. 
Applying the construction there, we get a differential graded $A_{\inf}$-algebra $(\qOmega_{\rA(R)}^{\bullet}, d_q^{\bullet})$ whose 
$k$-th differential map $d^k_q\colon \qOmega^k_{\rA(R)}:=\bigwedge^k\qOmega_{\rA(R)}^1\to \qOmega_{\rA(R)}^{k+1}$ is given by  
\[
d_q^k(f\frac{ dT_{i_1}}{T_{i_1}}\wedge \cdots\wedge\frac{d T_{i_k}}{T_{i_k}})=\sum_{j=1}^{d}\partial_j(f)\frac{dT_j}{T_j}\wedge \frac{ dT_{i_1}}{T_{i_1}}\wedge \cdots\wedge\frac{d T_{i_k}}{T_{i_k}}. 
\]
Following \cite{BMS1}, we call $(\qOmega_{\rA(R)}^{\bullet}, d_q^{\bullet})$ the  \emph{$q$-de Rham complex of $\rA(R)$}.

We recall now the notion of $q$-Higgs modules introduced in \cite{MT} and \cite{Tsuji}.
\begin{definition}[cf. \cite{Tsuji}, Def. 10.5]\label{D:q-Higgs-module}
	A \emph{$q$-Higgs module} over $\rA(R)$ is a $(p,[p]_q)$-adically complete $\rA(R)$-module with a  continuous integrable  connection
	$\Theta\colon E\to E\otimes_{\rA(R)}\qOmega^1_{\rA(R)}$ over $(\rA(R), d_q: \rA(R)\to \qOmega_{\rA(R)}^1)$ in the sense of Definition~\ref{D:connection}.
	In particular, if we write \[\Theta=\sum_{i=1}^{d} \theta_i\otimes \frac{dT_i}{T_i},\]
	then each $\theta_i$ is an $A_{\inf}$-linear   $\mu$-derivation on $E$  with respect to  $\partial_i$, i.e.
	\[
	\theta_i(fx)=\gamma_i(f)\theta_i(x)+\partial_i(f)x
	\]
	for all $f\in \rA(R)$ and $x\in E$, and $\theta_i\theta_j=\theta_j\theta_i$ for all $1\leqslant i,j\leq d$. 
	We will also write the Higgs module $(E,\Theta)$ as $(E,(\theta_i)_{1\leqslant i\leqslant d})$.
	The de Rham complex of a $q$-Higgs module $(E, \Theta)$ over $\rA(R)$ is denoted by $(E\otimes_{\rA(R)}\qOmega^{\bullet}_{\rA(R)}, \Theta^{\bullet})$ (see Definition~\ref{D:connection}(2)). 
\end{definition}
\begin{definition}\label{D:quasi-nilpotence}
	Let $(E, (\theta_{1\leqslant i\leqslant d}))$ be a $q$-Higgs module. For $\underline m=(m_1,\cdots, m_d)\in \NN^d$, we put $\underline m:=\sum_{i=1}^d m_i$ and 
	\[
	\theta^{\underline m}:=\prod_{i=1}^d\theta_i^{m_i}\in \End_{A_{\inf}}(E).
	\]  
	We say that $(E, (\theta_{1\leqslant i\leqslant d}))$  is \emph{topologically quasi-nilpotent}, if  we have
	\[
	\lim_{|\underline m|\to +\infty}\theta^{\underline m}(x)= 0, 
	\] 
	for all $x\in E$.
	
	Let $\qHiggs(\rA(R))=\MIC^{\wedge}(\rA(R), d_q)$ denote the category of $q$-Higgs modules over $\rA(R)$, and $\qHiggs^{\topnil}(\rA(R))$ be full subcategory consisting of the topological quasi-nilpotent objects.

\end{definition}

\begin{remark}\label{R:reduction of Higgs module}
	Note that $\partial_i(f)\in [p]_q \rA(R)$  and $\gamma_i(f)\equiv f\mod [p]_q\rA(R)$ for all $f\in \rA(R)$. It follows that $\theta_i(fx)\equiv f\theta_i(x)\mod [p]_qE$ for all $x\in E$. Therefore,  the reduction  modulo $[p]_q$ of $(E,(\theta_i)_{1\leqslant  i\leqslant d})$, which we denote by $(E/[p]_qE, (\bar \theta_i)_{1\leqslant  i\leqslant d})$,  is a Higgs module over $R$, and $(E,(\theta_i)_{1\leqslant  i\leqslant d})$ is topologically quasi-nilpotent if and only if $(E/[p]_qE, (\bar \theta_i)_{1\leqslant  i\leqslant d})$ is topologically quasi-nilpotent in the evident sense.
\end{remark}

\subsection{Stratified modules and $q$-Higgs modules}\label{S:simplicial-A(R)}
Let $\underline{\rA}^{\bullet}(R)$  be the \v{C}ech nerve of $\underline{\rA}(R)$  over the final object of $\underline {X}_{\prism}$. 
Each term $\underline{\rA}^{n}(R)$ for $n\geqslant 0$ can be described explicitly as follows. Let  $ \underline{\rA}^{{\otimes} (n+1)}(R):=(\rA^{{\otimes} (n+1)}(R), M^{\otimes (n+1)}_{A(R)})$ be the pushout of $(n+1)$-copies of $\underline {\rA}(R)$ over $(A_{\inf}, M_{\cO_C^{\flat}})$ in the category of $(p, [p]_q)$-adically complete prelog rings, i.e.
\[
{\rA}^{{\otimes} (n+1)}(R)= \underbrace{\rA(R)\widehat{\otimes}_{A_{\inf}}\rA(R)\widehat{\otimes}_{A_{\inf}}\cdots\widehat{\otimes}_{A_{\inf}}\rA(R)}_{(n+1)\text{-times}}
\]
and 
\[M^{\otimes (n+1)}_{R}=\underbrace{ M_{R}\sqcup_{M_{\cO_{C}^{\flat}}}\cdots \sqcup_{M_{\cO_{C}^{\flat}}}M_{R}}_{(n+1)\text{-times}}\cong M_{\cO_{C}^{\flat}}\sqcup_{\NN}\NN^{(r+1)(n+1)}  .\]
There is a natural $\delta_{\log}$-structure on $(\rA^{{\otimes} (n+1)}(R), M^{\otimes (n+1)}_{A(R)})$. 
Consider the following factorization  of surjections of prelog rings:
\[
\xymatrix{(\rA^{{\otimes} (n+1)}(R), M^{\otimes (n+1)}_{R})\ar@{->>}[rr] \ar[rd]&&  (\rA(R), M_{R})\\
	& (\widetilde {\rA}^{\otimes(n+1)}(R), \widetilde M_{R}^{\otimes (n+1)})\ar@{-->>}[ru]
}
\]
where the horizontal arrow is the diagonal surjection,  $(\widetilde {\rA}^{\otimes(n+1)}(R), \widetilde M_{R}^{\otimes (n+1)})$ is  its exactification.  Let $\widetilde J(n)\subset \widetilde {\rA}^{\otimes(n+1)}(R)$
be the kernel of the composite  surjection 
\[
\widetilde {\rA}^{\otimes(n+1)}(R)\twoheadrightarrow  \rA(R)\twoheadrightarrow R.
\]
We define $\rA^{n}(R):=\widetilde{\rA}^{\otimes (n+1) }(R)\{\frac{\widetilde J(n)}{[p]_q}\}^{\wedge}$ as the prismatic envelope of $\widetilde{\rA}^{\otimes (n+1) }(R)$ with respect to the ideal $\widetilde J(n)$, and  equip it with the natural induced prelog structure 
\[
\alpha_{\rA^n(R)}\colon M_{\rA^n(R)}:=\widetilde M_{R}^{\otimes (n+1)}\to\widetilde{\rA}^{\otimes (n+1) }(R)\to  \rA^{n}(R)
\]
and the induced $\delta_{\log}$-structure.
Then we have  $\underline \rA^n(R)=(\rA^n(R), ([p]_q),  M_{\rA^n(R)})^a$.

Let $\Strat(\rA^{\bullet}(R))$ be the category of (complete) stratified module  over the cosimplicial ring $\rA^{\bullet}(R)$ (Def.~\ref{D:stratified-modules}). As $\urA(R)$ is a strong covering object of $\uX_{\prism}$ (Lemma~\ref{L:cover-final-object}), Proposition~\ref{P:crystals-strat} implies that the evaluation at $\urA^{\bullet}(R)$ induces an equivalance of categories
\begin{equation}\label{E:equiv-crystal-strat-A(R)}
	\ev_{\rA(\uR)}\colon \CRhat(\uX_{\prism},\cO_{\prism})\xra{\sim} \Strat(\rA^{\bullet}(R))
\end{equation}


As in \S~\ref{S:simplicial notation}, we write $p_i\colon \urA(R)\to \urA^{1}(R)$ for the morphism of log prisms corresponding to $[0]\to [1]$, and $p_{i,j}\colon \urA^1(R)\to \urA^2(R)$ with $0\leqslant i<j\leqslant 2$ for the  morphisms corresponding to the three faces maps $[1]\to [2]$. 

Let $(E,\epsilon)$ be an object of $\Strat(\rA^{\bullet}(R))$. We  construct an   $\rA(R)$-semi-linear action of $\Gamma=\bigoplus_{i=1}^d\ZZ_p\gamma_i$ on $E$ as follows. 
For every $g\in \Gamma$, let $(1,g)\colon \rA^1(R)\to \rA(R)$ be the unique morphism such that $(1,g)\circ p_0=\id_{\rA(R)}$ and $(1,g)\circ p_1=g$. Then we define 
\begin{equation}\label{E:Gamma-action-E}
	\rho_E(g):=(1,g)^*(\epsilon)\colon \quad g^*(E)=(1,g)^* p_1^*(E)\xra{(1,g)^*(\epsilon)} (1,g)^* p_0^*(E)=E.
\end{equation}
The cocycle condition on $\epsilon$ implies  that $\rho_E(1)=\id_E$ and $\rho_E(g_1g_2)=\rho_E(g_1)\circ g_1^*\rho_{E}(g_2)$. If we identify $E$ with $g^*E=E\otimes_{\rA(R),g}\rA(R)$ via $x\mapsto x\otimes 1$, then $\rho_{E}$ is equivalent  to a semi-linear action of $\Gamma$ on $E$.

\begin{lemma}\label{L:crystal-q-Higgs}
	\begin{enumerate}
		\item[(1)] The semi-linear action $\rho_E$ of $\Gamma$ on $E$ is trivial modulo $\mu$, i.e. 
		\[
		\rho_E(g)(x)\equiv x\mod \mu E
		\]
		for all $g\in \Gamma$ and $x\in E$. 
		
		\item[(2)] If $E$ is  $\mu$-torsion free and we put 
		\[
		\theta_{ i}(x):=\frac{\rho_{E}(\gamma_i)(x)-x}{\mu}
		\]
		for $x\in E$ and $1\leqslant i\leqslant d$, then $(E,(\theta_{i})_{1\leqslant i\leqslant d})$ is a  $q$-Higgs module over $\rA(R)$. 
	\end{enumerate}
	
\end{lemma}
\begin{proof}
	
	(1)	For $g\in \Gamma$, let  $1\otimes g: \rA^1(R)\to \rA^1(R)$ denote the $p_0(\rA(R))$-linear action of $g$ on $\rA^1(R)$ induced by functoriality. Then one has $(1,g)=\Delta\circ (1\otimes g): \rA^1(R)\to \rA(R)$. Applying Proposition~\ref{P:derivation-D_B}(3) to $(B,([p]_q), M_B)=\underline\rA(R)$, we see that  $1\otimes g\equiv \id_{\rA^1(R)} \mod \mu \rA^1(R)$ for all $g\in \Gamma$. It follows that 
	$$\rho_E(g)=\Delta^*\circ (1\otimes g)^*(\epsilon)\equiv \Delta^*(\epsilon)=\id_{E}\mod \mu.$$
	
	Assertion (2) follows immediately from  Proposition~\ref{P:module with connection}(2). 
\end{proof}

We will give a more concrete description of the $q$-Higgs module $(E,(\theta_{i})_{1\leqslant i\leqslant d})$. We start with a concrete description of $\rA^1(R)$. Applying  Lemma~\ref{L:basis-for-D_B} to $(B,([p]_q), M_{\Spf(B)})=\urA(R)$ with $\tilde f$ being the identity map, we get 
\[
\rA^1(R)= \bigg\{\sum_{m\in \NN^d}a_m\tau^{\{m\}_{\delta}}\;\bigg|\; a_m\in \rA(R), \lim_{|m|\to +\infty}a_m=0\bigg\},
\] 
where $|m|=\sum_{i=1}^dm_i$ for $m\in \NN^d$ such that 
\begin{itemize}
	\item $p_0\colon \rA(R)\to \rA^1(R)$ is the natural inclusion,
	\item  $p_1\colon \rA(R)\to \rA^1(R)$  sends $T_i$ to $T_i(1+[p]_q\tau_i)$ for $1\leqslant i\leqslant d$, 
	\item and the diagonal map $\Delta\colon \rA^1(R)\to \rA(R)$ is given by $\sum_{m\in \NN^d}a_m\tau^{\{m\}_{\delta}}\mapsto a_0$. 
\end{itemize}  
According to \eqref{E:gamma-action-tau}, the $p_0(\rA(R))$-linear action of  $1\otimes \gamma_i$  for $1\leqslant i\leqslant d$ on  $\rA^1(R)$ is  given by 
\[
(1\otimes \gamma_i)\bigg(\sum_{m\in \NN^d}a_m\tau^{\{m\}_{\delta}}\bigg)=\sum_{m\in \NN^d} a_m \big(\prod_{j\neq i}\tau_j^{\{m_j\}_{\delta}} \big)(\tau_i+\mu u_i)^{\{m_i\}_{\delta}},
\]
where $u_i=1+[p]_q\tau_i$. 

Consider now   the composed map
\[
\epsilon^{\sharp}\colon E\xra{x\mapsto 1\otimes x}p_1^*E\xra{\epsilon}p_0^*E=E\widehat{\otimes}_{p_0, \rA(R)}\rA^1(R).
\]
Then there exist endomorphisms $b_m\in \End_{A_{\inf}}(E)$ for each $m\in \NN^d$ with $b_0=\id_E$ such that 
\begin{equation}\label{E:explicit-strat}
	\epsilon^{\sharp}(x)=\sum_{m\in \NN^d}b_m(x)\otimes \tau^{\{m\}_{\tau}},\quad \forall x\in E
\end{equation}
By the construction of the $\Gamma$-action on $E$, we have 
\begin{equation*}
	\rho_E(\gamma_i)(x)=\sum_{k=0}^{\infty}b_{ke_i}(x) \mu^{\{k\}_{\delta}},\quad \forall x\in E
\end{equation*}
where $e_i\in \NN^d$ is the element whose $i$-th component equals to $1$ and other components equal to $0$. 
Note that $\delta^{k}(\mu)\in \mu A_{\inf}^{\times}$.
It follows that $\mu^{\{k\}_{\delta}}\in \mu A_{\inf}$ for $k\geqslant 1$ and 
\begin{equation}\label{E:formula for theta}
	\theta_i(x)=\sum_{k=1}^{\infty} b_{ke_i}(x)\frac{\mu^{\{k\}_{\delta}}}{\mu}.
\end{equation}

\begin{remark}
	One could also take \eqref{E:formula for theta} as the definition of $\theta_i$, which makes sense even without assuming that $E$ is $\mu$-torsion free. But it seems more difficult to deduce from this that $(E, (\theta_i)_{1\leqslant i\leqslant d})$ is $q$-Higgs module.  
\end{remark}

\begin{proposition}\label{P:top-nilpotence}
	Let $(E,\epsilon)$ be a $\mu$-torsion free object of $\Strat(\rA^{\bullet}(R))$. Then its associated $q$-Higgs module $(E, (\theta_i)_{1\leqslant i\leqslant d})$  by Lemma~\ref{L:crystal-q-Higgs} is topologically quasi-nilpotent (Definition~\ref{D:quasi-nilpotence}).
\end{proposition}
\begin{proof}
	Put $\overline E:=E/[p]_qE$. By Remark~\ref{R:reduction of Higgs module}, it suffices to show that the induced Higgs module $(\overline E, (\bar \theta_i)_{1\leqslant i\leqslant d})$ is topological quasi-nilpotent. 
	First,  Proposition~\ref{P:dp-envelop-D_B} gives   an isomorphism of $R$-algebras
	\[
	R\langle \tau_1, \dots, \tau_d\rangle^{PD}\xra{\sim} \rA^1(R)/([p]_q).
	\]
	Via this isomorphism, the two natural face maps $p_0,p_1\colon R=\rA(R)/([p]_q)\to \rA^1(R)/([p]_q)$ are both identified with the natural inclusion $R\hra R\langle \tau_1, \dots, \tau_d\rangle^{PD} $. 
	Note that  $\urA^2(R)$ is the  pushout  (or fiber product in $\uX_{\prism}$) of $p_0\colon \urA(A)\to \urA^{1}(R)$ over  itself. This implies that 
	\[
	\rA^2(R)=\rA^1(R)\widehat\otimes_{\rA(R)}\rA^1(R).
	\]
	Hence, there exists an isomorphism 
	\[
	R\langle \tau_{i,1}, \tau_{i,2}: 1\leqslant i\leqslant d\rangle\xra{\sim}  \rA^2(R)/([p]_q)
	\]
	such that $p_{0,1}, p_{1,2}, p_{0,2}\colon \rA^1(R)/([p]_q)\to\rA^2(R)/([p]_q)$ are given by   
	\[
	p_{0,1}(\tau_i)=\tau_{i,1}, \quad p_{1,2}(\tau_i)=\tau_{i,2}^2, \quad p_{0,2}(\tau_{i})=\tau_{i,1}+\tau_{i,2}.
	\]
	
	Let   $\bar{\epsilon}^{\sharp}\colon \overline E\to \overline E\widehat{\otimes}_{R}\rA^1(R)/([p]_q)$
	denote the reduction of $\epsilon^{\sharp}$. Then there exist $\bar\theta^{\log}_m\in \End_R(\overline E)$ for all  $m\in \NN^d$ with $\bar\theta^{\log}_0=\id_{\overline E}$ such that 
	\[
	\bar{\epsilon}^{\sharp}(x)=\sum_{m\in \NN^d}^{+\infty}\bar\theta^{\log}_m(x)\otimes \tau^{[m]}, \quad \forall x\in \overline E
	\]
	where  $\tau^{[m]}=\prod_{i=1}^d \tau_i^{[m_i]}$, and $\lim_{|m|\to\infty}\bar\theta^{\log}_{m}(x)= 0$.  By the same argument as  \cite[Prop. 5.11]{Tian}, it follows from the cocycle condition $p_{02}^*(\epsilon)=p_{01}^*(\epsilon)\circ p_{12}^*(\epsilon)$ that $\bar\theta^{\log}_{m_1+m_2}(x)=\bar\theta^{\log}_{m_1}(\bar\theta^{\log}_{m_2}(x))$ for all $x\in \overline E$ and $m_1,m_2\in\NN^d$. It follows therefore that $\bar\theta^{\log}_m=(\bar\theta^{\log})^m:=\prod_{i=1}^{d}(\bar\theta^{\log}_{i})^{m_i}$ for all $m\in \NN^d$, where $\bar\theta^{\log}_i=\bar \theta_{e_i}^{\log}$. 
	
	On the other hand, by the proof of Proposition \ref{P:dp-envelop-D_B},  the image of $\tau^{\{m\}_{\delta}}\in \rA^1(R)$ in  $\rA^1(R)/([p]_q)$  is $c_m\tau^{[m]}$ with \[c_m=\prod_{i=1}^d(-1)^{l_i(l_i+1)/2}\frac{m_i!}{p^{v_p(m_i!)}}\in \ZZ_{(p)}^{\times}\]
	for all $m=(m_i)\in \NN^d$, where $l_i$ is the minimal integer such that $m<p^{l_i}$. 
	By taking the reduction of \eqref{E:explicit-strat} modulo $[p]_q$  and comparing with the formula above for $\bar{\epsilon}(x)$, we get
	$
	\bar b_{m}= c_m \bar{\theta}_m^{\log}
	$
	for all $m\in \NN^d$, 
	where $\bar b_m$ denotes the reduction modulo $[p]_q$ of the endomorphism $b_m\in \End_{A_{\inf}}(E)$. It follows from \eqref{E:formula for theta} that 
	\[
	\bar\theta_i=\sum_{k=1}^{\infty} (\bar\theta_i^{\log})^k \frac{(\zeta_p-1)^{k-1}}{k!}
	\]
	for $1\leqslant i\leqslant d$. Because $\lim_{m\to+\infty}(\bar\theta^{\log})^m(x)=\lim_{m\to+\infty}\bar\theta^{\log}_m(x)=0$ for all $x\in \overline E$, it follows that 
	\[
	\lim_{m\to +\infty}\theta^m(x)=0
	\]
	for all $x\in \overline E$. This proves that the Higgs module $(\overline E,(\bar \theta_i)_{1\leqslant i\leqslant d})$ is topologically quasi-nilpotent. 
	
\end{proof}

For $\scrC\in \{ \CRhat(\uX_{\prism}, \cO_{\prism}), \Strat(\urA^{\bullet}(R)), \qHiggs^{\topnil}(\rA(R))\}$, let $\scrC^{\mu-tf}$ denote the subcategory consisting of objects of $\scrC$ that are $\mu$-torsion free. Then  we have a functor 

\begin{equation}\label{E:strat-Higgs}
	\CRhat(\uX_{\prism}, \cO_{\prism})^{\mu-tf}\simeq 	\Strat(\urA^{\bullet}(R))^{\mu-tf}\to \qHiggs^{\topnil}(\rA(R))^{\mu-tf}.
\end{equation}
where the first equivalence is \eqref{E:equiv-crystal-strat-A(R)}, and the second functor is given by Proposition~\ref{P:top-nilpotence}.

\begin{remark}
	One expects that the  functor $\Strat(\urA^{\bullet}(R))^{\mu-tf}\xra{\sim} \qHiggs^{\topnil}(\rA(R))^{\mu-tf}$ is actually an equivalence of categories. 
	When $X$ is smooth over $\Spf(\cO_C)$, this is  proved by Tsuji in \cite[Thm. 11.20]{Tsuji} (even without  the $\mu$-torsion free condition). 
	Here, we choose to just work with  $\mu$-torsion free objects in $\CRhat(\uX_{\prism}, \cO_{\prism})$,  since this simplifies our discussion and it  suffices already for our purpose.
\end{remark}

\subsection{Frobenius action} Let $\calE$ be an object of $\CRhat(\uX_{\prism}, \cO_{\prism})^{\mu-tf}$. Let $(E, \Theta)$ be the $q$-Higgs module over $\rA(R)$ corresponding to $\calE$. 
Since the $\Gamma$-action on $\rA(R)$ commutes with the $\delta$-structure, there exists a natural induced  tensor $\Gamma$-action $\rho_{\phi^*_{\rA(R)}E}$ on  $\phi_{\rA(R)}^*E=E\widehat\otimes_{\rA(R), \phi_{\rA(R)}}\rA(R)$. We see easily that $$\rho_{\phi_{\rA(R)}^*E}(g)(x\otimes a)-x\otimes a=(\rho_{E}(g)(x)-x)\otimes g(a)+x\otimes (g(a)-a)\in \mu \phi^*_{\rA(R)}E$$ 
for all $x\otimes a\in \phi_{\rA(R)}^*E$ and $g\in \Gamma$. 
Since $\phi_{\rA(R)}$ is $(p,[p]_q)$-completely flat,  $\phi^*_{\rA(R)}E$ is also $\mu$-torsion free, and     
\[
\theta^{\phi}_i\colon =\frac{\rho_{\phi^*_{\rA(R)}E}(\gamma_i)-\id_{\phi^*_{\rA(R)}E}}{\mu}
\]
for $1\leqslant i\leqslant $
is a well defined $\rA(R)$-linear endomoprhism of $\phi^*_{\rA(R)}E$.
Put 
\[
\Theta^{\phi}\colon =\sum_{i=1}^d\theta_i^{\phi}\otimes \frac{d T_i}{T_i}.
\]
We get then a $q$-Higgs module $(\phi_{\rA(R)}^*E, \Theta^{\phi})$, which is nothing but the $q$-Higgs module associated to the  {Frobenius twist} $\phi_{\cO_{\prism}}^*\calE$.

Assume now that $\calE$ is equipped with a Frobenius structure $$\phi_{\calE}\colon (\phi_{\cO_{\prism}}^*\calE)[1/[p]_q]\simeq \calE[1/[p]_q]$$ 
so that $(\calE, \phi_{\calE})$ is an object of $\CRhat(\uX_{\prism}, \cO_{\prism})^{\phi=1}$. Then $\phi_{\calE}$ induces an isomorphism of $\rA(R)[\frac{1}{[p]_q}]$-modules $$\phi_E\colon (\phi_{\rA(R)}^*E)[\frac{1}{[p]_q}]\simeq E[\frac{1}{[p]_q}]$$
that is compatible with the $\Gamma$-actions, and hence with the $q$-Higgs structures on both sides. We get thus an isomoprhism   of the associated  de Rham complexes  after inverting $[p]_q$:
\[
\phi_E^{\bullet}\colon \big((\phi_{\rA(R)}^*E)\otimes \qOmega_{\urA(R)}^{\bullet}, \Theta^{\phi, \bullet}\big)[\frac{1}{[p]_q}]\xra{\sim} \big(E\otimes_{\rA(R)}\qOmega_{\rA(R)}^\bullet, \Theta^{\bullet}\big)[\frac{1}{[p]_q}].
\]
Let $\iota\colon E\to \phi_{\rA(R)}^*E$ be the natural map $x\mapsto x\otimes 1$. For every integer $n\geqslant0$, let $d^n\phi\colon \qOmega^n_{\rA(R)}\to \qOmega^{n}_{\rA(R)}$ be the Frobenius linear map defined by 
\[
d^n\phi \bigg(	f \frac{dT_{i_1}}{T_{i_1}}\wedge \cdots \wedge\frac{dT_{i_n}}{T_{i_n}}\bigg)= [p]_q^n\phi_{\rA(R)}(f)\frac{dT_{i_1}}{T_{i_1}}\wedge \cdots \wedge\frac{dT_{i_n}}{T_{i_n}},
\]
for all $f\in \rA(R)$ and $1\leqslant i_1< \cdots < i_n\leqslant d$. 
Then  one check easily that $\iota\otimes d^{\bullet}\phi$ defines a morphism of attached de Rham complexes:
\begin{equation}\label{E:phi-functorial-deRham}
	\iota\otimes d^{\bullet}\phi\colon\quad  (E\otimes_{\rA(R)}\qOmega^{\bullet}_{\rA(R)}, \Theta^{\bullet})\to \big((\phi_{\rA(R)}^*E)\otimes_{\rA(R)}\qOmega^{\bullet}_{\rA(R)}, \Theta^{\phi,\bullet}\big)
\end{equation}
which is $\phi$-linear as a morphism of complexes of $A_{\inf}$-modules.
By inverting $[p]_q$ and composing with $\phi^\bullet_E$, we get then  a Frobenius action on the   de Rham complexes of $(E,\Theta)$:
\begin{equation}\label{E:Frobenius-Higgs-complex}
	\phi^*\big( E\otimes_{\rA(R)}\qOmega_{\rA(R)}^{\bullet}, \Theta^{\bullet}\big)[\frac{1}{[p]_q}]\xra{\iota\otimes d^{\bullet}\phi}\big((\phi_{\rA(R)}^*E)\otimes_{\rA(R)} \qOmega_{\rA(R)}^{\bullet}, \Theta^{\phi, \bullet}\big)[\frac{1}{[p]_q}]\xra{\phi_E^{\bullet},\,\simeq}  \big(E\otimes_{\rA(R)}\qOmega_{\rA(R)}^{\bullet}, \Theta^{\bullet}\big)[\frac{1}{[p]_q}],
\end{equation}
where $\phi^*$ on the first term means base change along $\phi\colon A_{\inf}\xra{\sim} A_{\inf}$.

\begin{theorem}\label{T:crystal-qHiggs}
	Let $\calE$ be an obejct of $\CRhat(\uX_{\prism}, \cO_{\prism})^{\mu-tf}$, and let $(E, \Theta)$ be the  $q$-Higgs module over $\rA(R)$ attached to $\calE$ via the functor \eqref{E:strat-Higgs}. Then:
	
	\begin{enumerate}
		\item[(1)] There exists a canonical isomorphism 
		\[
		R\Gamma(\uX_{\prism}, \calE)\simeq E\otimes_{\rA(R)}\qOmega_{\rA(R)}^{\bullet}.
		\]
		\item[(2)] If $\calE$ is equipped with a Frobenius structure $\phi_{\calE}\colon \phi_{\cO_{\prism}}^*\calE[\frac{1}{[p]_q}]\simeq \calE[\frac{1}{[p]_q}]$, then the isomorphism in (1) is compatible with the Frobenius action, i.e. the diagram 
		\[
		\xymatrix{\phi^*R\Gamma(\uX_{\prism}, \calE)[\frac{1}{[p]_q}]\ar[r]\ar[d]^{(1)}_{\simeq} &R\Gamma(\uX_{\prism}, \phi_{\cO_{\prism}}^*\calE)[\frac{1}{[p]_q}]\ar[r]^{\phi_{\calE}, \simeq} & R\Gamma(\uX_{\prism}, \calE)[\frac{1}{[p]_q}]\ar[d]^{(1)}_{\simeq}\\
			\phi^*(	E\otimes_{\rA(R)}\qOmega^{\bullet}_{\rA(R)}, \Theta^{\bullet})[\frac{1}{[p]_q}]\ar[rr]^-{\eqref{E:Frobenius-Higgs-complex}} && (E\otimes_{\rA(R)}\qOmega^{\bullet}_{\rA(R)}, \Theta^{\bullet})[\frac{1}{[p]_q}]
		}
		\]
		is commuative, where the left upper arrow is induced by the natural Frobenius linear map $\calE\to \phi_{\cO_{\prism}}^*\calE$ and vertical maps are given by the isomorphism in (1).
	\end{enumerate}
\end{theorem}

\begin{remark}
	When $X$ is smooth over $\Spf(\cO_C)$, Tsuji gives a similar description of $R\Gamma(X_{\prism}, \calE)$ for any object $\calE$ of $\CRhat(X_{\prism}, \cO_{\prism})$ (cf. \cite[Theorem~13.9]{Tsuji}). Our arguments are obviously inspired by his. 
\end{remark}

To prove Theorem~\ref{T:crystal-qHiggs}, we need some preparation.  Let $\Mod(\rA(\rR))$ denote the category of $(p, [p]_q)$-adically complete $\rA(R)$-modules. Following \cite[IV. 3]{Berthelot} and \cite[12.2]{Tsuji},  we  construct introduce a linearization functor 
\[
\calL_{A}\colon \Mod(\rA(R))\to \CRhat(\uX_{\prism}, \cO_{\prism})
\]
as follows. Let $M$ be an object of $\Mod(\rA(R))$. For an object $(B, ([p]_q), M_{\Spf(B)})$  of $\uX_{\prism}$, recall that  $(D_B, ([p]_q), M_{\Spf(D_B)})$  is the  coproduct of $(B, ([p]_q), M_{\Spf(B)})$ with $\urA(R)$ in $\uX_{\prism}^{\opp}$ constructed in Lemma~\ref{L:cover-final-object}.  We put 
\[
\calL_A(M)(B,([p]_q), M_{\Spf(B)}):=M\widehat\otimes_{\rA(R)}D_B.
\]
Now let $f\colon (B,([p]_q), M_{\Spf(B)})\to (C, ([p]_q), M_{\Spf(C)})$ be a morphism in $\uX_{\prism}^{\opp}$. Then Lemma~\ref{L:cover-final-object}(2) implies that $D_B\widehat \otimes_{B}C\cong D_C$. It follows therefore that the canonical transition map 
\[\calL_A(M)(B,([p]_q), M_{\Spf(B)})\widehat{\otimes}_{B}C\xra{\sim} \calL_A(M)(C,([p]_q), M_{\Spf(C)})
\]
is an isomorphism. This shows that $\calL_A(M)$ is indeed an object of $\CRhat(\uX_{\prism}, \cO_{\prism})$. The natural map $M\to M\widehat\otimes_{\rA(R)}D_B$ for any object $(B, ([p]_q), M_{\Spf(B)})$ of $\uX_{\prism}$ induces  a canonical map
\[
M\to \Gamma(\uX_{\prism}, \calL_A(M)).
\]
\begin{proposition}\label{P:cohomology-vanishing}
	For any object $M$ of $\Mod(\rA(R))$, the canonical map above induces an isomorphism 
	\[
	M\xra{\sim}R\Gamma(\uX_{\prism}, \calL_A(M)).
	\]
\end{proposition}
\begin{proof}
	Since the proof is similar to \cite[Prop. 12.5]{Tsuji}, we just give here a sketch of the arguments.  By Lemma~\ref{L:CA-complex-prismatic} and Lemma~\ref{L:cover-final-object}, $R\Gamma(\uX_{\prism}, \calL_A(M))$	 is computed by the \v{C}ech--Alexander complex $\CA(\urA^{\bullet}(R), \calL_{A}(M))$, which is the simple complex attached to cosimplicial abelian group $\calL_A(M)(\urA^{\bullet}(R))=M\widehat\otimes_{\rA(R)}  D_{\urA^{\bullet}(R)} $. Here, $D_{\urA^{\bullet}(R)}$ means the cosimplicial object by applying the functor $(B,([p]_q), M_{\Spf(B)})\mapsto D_B$ to the simplicial object $\urA^{\bullet}(R)$ of $\uX_{\prism}$. 
	
	Note that for every integer $n\geqslant0$, we have an isomorphism 
	\[
	D_{\urA^n(R)}\cong \rA^{n+1}(R)\cong \underbrace{\rA^{1}(R)\widehat\otimes_{\rA(R)}\rA^1(R)\widehat\otimes_{\rA(R)} \cdots \widehat{\otimes}_{\rA(R)}\rA^1(R)}_{\text{$n$ times}},
	\]
	where each $\rA^1(R)$ is viewed as an $\rA(R)$-algebra via $p_0$ (cf. \S \ref{S:simplicial notation} for the notation), and the simplicial object $D_{\urA^\bullet(R)}=\rA^{\bullet+1}(R)$ can be identified with the \v{C}ech nerve of the map $p_0\colon \rA(R)\to \rA^1(R)$.  Note that the diagonal map $\Delta: \rA^1(R)\to \rA(R)$ is  a section  to $p_0$. By \cite[\href{https://stacks.math.columbia.edu/tag/019Z}{Tag 019Z}]{stacks-project}, the cosimplicial abelian group $D_{\urA^\bullet(R)}$ is homotopic to the constant cosimplicial object $\rA(R)$. It follows that $M\widehat\otimes_{\rA(R), p_A}  D_{\urA^{\bullet}}(R)$ is homotopic to $M$, hence $M\xra{\sim} \CA(\urA^{\bullet}(R), \calL_{A}(M))$ is a quasi-isomorphism.
	
\end{proof}
Now let $\calE$ be an object of $\CRhat(\uX_{\prism},\cO_{\prism})^{\mu-tf}$, and let $(E,\Theta)$ be the corresponding object of in $\qHiggs(\rA(R))^{\mu-tf}$ via the functor \eqref{E:strat-Higgs}. 
We will construct a sheaf of complex $\calL_A(E\otimes_{\rA(R)}\qOmega_{\rA(R)}^{\bullet}, \Theta^{\bullet})$ on $\uX_{\prism}$ as follows. 
Let $(B,([p]_q), M_{\Spf(B)})$ be an object of $\uX_{\prism}$. We write   $p_A\colon \urA(R)\to (D_B, ([p]_q), M_{\Spf(D_B)})$ and $p_B\colon  (B,([p]_q), M_{\Spf(B)})\to  (D_B, ([p]_q), M_{\Spf(D_B)})$ for  the two  canonical maps in $\uX_{\prism}^{\opp}$.
Recall that the  $\mu$-derivations $(\partial_{D_B,i})_{1\leqslant  i\leqslant d}$ on $D_B$ is compatible with  $(\partial_i)_{1\leqslant i\leqslant d}$ on $A$ via $p_A$ (Prop. ~\ref{P:derivation-D_B}(2)).	 
By construction in  \S~\ref{S:scalar-extension-MIC}, $p_A:\rA(R)\to D_B$ induces a map of graded differential algebras 
$$dp_A\colon \qOmega_{\rA(R)}^{\bullet}\to \qOmega_{D_B/B}^{\bullet}, \quad \frac{d T_i}{T_i}\mapsto \frac{d u_i}{u_i}$$
and  a scalar extension functor 
$$
p_A^*\colon \qHiggs(\rA(R))=\MIC^{\wedge}(\rA(R), d_q)\to  \MIC^{\wedge}(D_B, d_{D_B/B})
$$
that sends $(E,\Theta)$ to $(p_A^*(E), \nabla_{E_{D_B}})$ with $p_A^*(E)=E\widehat\otimes_{\rA(R)}D_B$ and  $\nabla_{E_{D_B}}:=\Theta\otimes 1+1\otimes d_{D_B/B}$.
We define  
\[
\calL_A(E\otimes_{\rA(R)}\qOmega_{\rA(R)}^{\bullet}, \Theta^{\bullet})(B, ([p]_q), M_B):=(E\widehat\otimes_{\rA(R)} \qOmega_{D_B/B}^{\bullet}, \nabla_{E_{D_B}}^{\bullet}),
\]
the de Rham complex of $(p_A^*(E), \nabla_{E_{D_B}})$. 
If $f:(B, ([p]_q), M_{\Spf(B)})\to (C,([p]_q), M_{\Spf(C)})$ is a morphsim in $\uX^{\opp}_{\prism}$, by the functoriality of the scalar extension functor, $f$ induces a natural morphism of complexes 
\[
(E\widehat\otimes_{\rA(R)} \qOmega_{D_B/B}^{\bullet}, \nabla_{E_{D_B}}^{\bullet})\to (E\widehat\otimes_{\rA(R)} \qOmega_{D_C/C}^{\bullet}, \nabla_{E_{D_C}}^{\bullet}).
\]
Note that	for each fixed integer $n\geqslant 0$, one has an isomorphism 
$$
E\widehat\otimes_{\rA(R)} \qOmega_{D_B/B}^n\cong (E\otimes_{\rA(R)} \qOmega_{\rA(R)}^n)\widehat{\otimes}_{\rA(R)}D_B=\calL_A(E\otimes_{\rA(R)} \qOmega_{\rA(R)}^n)(B, ([p]_B), M_{\Spf(B)}).
$$ 
Hence, $\calL_A(E\otimes_{\rA(R)}\qOmega_{\rA(R)}^{\bullet}, \Theta^{\bullet})$ can be viewed as   a complex of abelian sheaves on $\uX_{\prism}$ whose $n$-th degree  term is $\calL_A(E\otimes_{\rA(R)} \qOmega_{\rA(R)}^n)$. 

For every object $(B, ([p]_q), M_{\Spf(B)})$, the crystal property of $\calE$ implies that  there exist isomorphisms of $D_B$-modules:
\begin{equation}\label{E:isom-E-D_B}
	p_A^*(E)=E\widehat{\otimes}_{\rA(R)}D_B\xra{c_{
			\calE}(p_A),\; \sim}\calE_{D_B}\xleftarrow{\sim,\; c_{\calE}(p_B)}\calE_{B}\widehat\otimes_{B}D_B=p_B^*(\calE_{B}),
\end{equation}
where $\calE_{B}$ (resp.  $\calE_{D_B}$) stands for the evaluation of $\calE$ at $(B,([p]_q), M_{\Spf(B)})$ (resp. at $(D_B, ([p]_q), M_{\Spf(D_B)})$).

\begin{proposition}[cf. \cite{Tsuji}, Prop. 13.24]\label{P:triviality of derivation}
	Via the isomorphism \eqref{E:isom-E-D_B}, the connection $\nabla_{E_{D_B}}$ on $p_A^*(E)$ corresponds to $1\otimes d_{D/B}\colon p_{B}^*(\calE_B)=\calE_B\widehat{\otimes}_{B}D_B\to \calE_B\widehat{\otimes}_{B}\qOmega^1_{D_B/B}$. 
\end{proposition}
\begin{proof}
	
	The problem is clearly  local for the flat topology on $\uX_{\prism}$. Up to replacing $\underline B:=(B,([p]_q), M_{\Spf(B)})$ by a flat cover, we may assume that there exists a map $\urA(R)\to \underline B$ because $\urA(R)$ is a flat cover of the final object of the topos $\uX_{\prism}$ by Lemma~\ref{L:cover-final-object}. 
	Since the formation of $\nabla_{E_{D_B}}$ commutes with base change in $\underline B$,  we are reduced to the case $\underline B=\urA(R)$ and thus $(D_B, ([p]_q), M_{\Spf(D_B)})=\urA^1(R)$. 
	We will make the convention that $p_A=p_1\colon \urA(R)\to \urA^1(R)$ and $p_B=p_0\colon \urA(R)\to \urA^1(R)$. Then the isomorphism \eqref{E:isom-E-D_B} is identified with  the stratification $\varepsilon \colon p_1^*E=p_A^*E\xra{\sim} p_0^*E=p_B^*E$ on $E$ given by the crystal $\calE$ (cf.  Proposition~\ref{P:crystals-strat}). 
	
	As $p_{A}^*E=E\widehat\otimes_{\rA(R), p_1}\rA^1(R)$ is  $\mu$-torsion free, the connection $\nabla_{E_{D_B}}$ is determined by its induced  $\Gamma$-action on  $p_{A}^*E$ by Proposition~\ref{P:module with connection}. Therefore,  it suffices to show that the $\Gamma$-action on $p_A^{*}E$ given by $\nabla_{E_{D_B}}$ corresponds to the $\Gamma$-action on $p_B^*E=E\widehat\otimes_{\rA(R),p_0} \rA^1(R)$ induced by the trivial connection $1\otimes d_{\rA^1(R)/p_0(\rA(R))}$.  
	
	Consider the diagram in $\uX_{\prism}$ 
	\[
	\xymatrix{ && \urA(R)\ar@<1ex>[r]^{p_0}\ar@<-1ex>[r]_{p_1}\ar[d]_{p_A=p_1} & \urA^1(R)\ar[d]^{p_{12}}\\
		\underline{B}=\urA(R) \ar[rr]^-{p_B=p_0} 
		&& \urA^1(R) \ar@<1ex>[r]^{p_{01}} \ar@<-1ex>[r]_{p_{02}} & \urA^{2}(R),	
	}
	\]
	where the square is cocartesian, i.e. $p_{01}$ (resp. $p_{02}$) is the pushout of $p_0$ (resp. $p_1$) via $p_1$.  According to \eqref{E:Gamma-action-E}, the action of $g\in \Gamma$ on $E$ is the base change  via $(1,g)\colon \rA^1(R)\to \rA(R)$ of the stratification $\epsilon \colon p_1^*(E)\xra{\sim}p_0^*(E)$. By the discussion in \S \ref{S:scalar-extension-MIC}, the action of $g$ on $p_A^*E$ determined by $\nabla_{E_{D_B}}$ is then given by the base change via $p_{12}^*(1,g)\colon \rA^2(R)\to \rA^1(R)$ of 
	\[
	p_{12}^*(\epsilon)\colon p_{12}^*(p_1^*E)=p_{02}^*(p_{A}^*E)\xra{\sim } p_{12}^*(p_0^*E)=p_{01}^*(p_A^*E)
	.\] 
	By  the cocylce condition for  $\epsilon$, one has  a commutative diagram  of isomorphisms
	\[
	\xymatrix{p_{02}^*(p_A^*E)\ar[rr]^{p_{12}^*(\epsilon)}\ar[d]_{p_{02}^*(\epsilon)} & & p_{01}^*(p_A^*(E))\ar[d]^{p_{01}^*(\epsilon)}\\
		(p^2_0)^*(E)=	p_{02}^*(p_B^*E)\ar@{=}[rr]^{\id_{(p^{2}_0)^*E}}
		&& p_{01}^*(p_{B}^*E)=(p_0^{2})^*E},
	\]
	where $p^2_0\colon \urA(R)\to \urA^2(R)$ is the morphism corresonding to map of simplices $[0]\to [2]$ sending $0$ to $2$. Via the isomorphism \eqref{E:isom-E-D_B}, the action by $g$ on $p_B^*(E)$ is by definition obtained by the base change via $p_{12}^*(1,g)$ of the bottom horizontal arrow, which is the identify map. It follows  that  the action by $g$ on $p_B^*E=p_0^*E$  corresponds to the trivial derivation $1\otimes d_{\rA^1(R)/p_0(\rA(R))}$.

\end{proof}

\begin{corollary}\label{C:resolution-LAE}
	For every object $(B,([p]_q), M_{\Spf(B)})$ of $\uX_{\prism}$, via the isomorphism \eqref{E:isom-E-D_B}, the canonical map $\calE_B\to p_B^*(\calE_B)$ induces a quasi-isomorphism of complexes 
	\[
	\calE_B\xra{\sim}(E\widehat\otimes_{\rA(R)}\qOmega^{\bullet}_{D_B/B}, \nabla_{E_{D_B}}^{\bullet})
	\]
	which is functorial in $(B,([p]_q), M_{\Spf(B)})$.
\end{corollary}

\begin{proof}
	This follows immediately from Propositions \ref{P:poincare-lemma} and \ref{P:triviality of derivation}.
	
\end{proof}

We have now all the ingredients to finish the proof of Theorem~\ref{T:crystal-qHiggs}.
\begin{proof}[End of the proof of ~\ref{T:crystal-qHiggs}]
	
	(1)	For all integers $m,n\geqslant0$, let
	\[
	\scrC^{m,n}:=\calL_A(E\otimes_{\rA(R)}\qOmega^n_{\rA(R)})(\urA^m(R))\cong E\widehat{\otimes}_{\rA(R)}\qOmega^n_{D_{\rA^m(R)}/\rA^m(R)}, 
	\]
	$d_1^{m,n}\colon \scrC^{m,n}\to \scrC^{m+1,n}$ be the $m$-th differential of  $\CA(\urA^{\bullet}(R),\calL_A(E\otimes_{\rA(R)}\qOmega_{\rA(R)}^n))$,  the \v{C}ech--Alexander for complex for $\calL_A(E\otimes_{\rA(R)}\qOmega^n_{\rA(R)})$, and  $d_2^{m,n}\colon \scrC^{m,n}\to \scrC^{m,n+1}$ be the $n$-th differential of $(E\widehat{\otimes}_{\rA(R)}\qOmega^\bullet_{D_{\rA^m(R)}/\rA^m(R)}, \nabla^{\bullet}_{E_{D_{\rA^m(R)}}})$. By the functoriality of \v{C}ech--Alexander complex, one sees  
	that $d_1^{m,n+1}\circ d_2^{m,n}=d_2^{m+1,n}\circ d_1^{m,n}$ for all $m,n\in \NN$, i.e.   $(\scrC^{\bullet,\bullet}, d_1^{\bullet, \bullet}, d_2^{\bullet, \bullet})$ forms  a double complex in the first quadrant. Write $\underline{s}(\scrC^{\bullet})$ for the simple complex associated to $\scrC^{\bullet, \bullet}$.

	By Corollary~\ref{C:resolution-LAE},  one has a canonical quasi-isomorphism 
	\[
	\calE(\urA^m(R))\xra{\sim} (\scrC^{m,\bullet}, d_2^{m,\bullet})=\calL_{A}(E\otimes_{\rA(R)}\qOmega_{\rA(R)}^{\bullet}, \Theta^{
		\bullet})(\urA^m(R))
	\]
	for a fixed integer $m\geqslant 0$. 
	By functoriality, the induced map  $\calE(\urA^m(R))\to \calE(\urA^{m+1}(R))$ by the morphism of complexes $$d_1^{m,\bullet}\colon (\scrC^{m,\bullet}, d_2^{m,\bullet})\to (\scrC^{m+1,\bullet}, d_2^{m+1, \bullet})$$ is exactly the $m$-th differential in the \v{C}ech--Alexander complex $\CA(\urA^{
		\bullet}(R), \calE)$. We get thus a canonical quasi-isomorphism of complexes 
	\begin{equation}\label{E:double complex 1}
		\CA(\urA^{
			\bullet}(R), \calE)\simeq \underline {s}(\scrC^{\bullet,\bullet}).
	\end{equation}
	On the other hand, for each integer $n\geqslant 0$,  there is a commutative diagram of complexes
	\[
	\xymatrix{E\otimes_{\rA(R)}\qOmega_{\rA(R)}^n\ar[rr]^-{\sim}\ar[d]^{\Theta^{n}} && (\scrC^{\bullet, n}, d_1^{\bullet, n})
		=\CA(\urA^{\bullet}(R),\calL_A(E\otimes_{\rA(R)}\qOmega_{\rA(R)}^n))	\ar[d]^{d_2^{\bullet, n}}\\
		E\otimes_{\rA(R)}\qOmega_{\rA(R)}^{n+1}\ar[rr]^-{\sim} && (\scrC^{\bullet, n+1}, d_1^{\bullet, n+1})	=\CA(\urA^{\bullet}(R),\calL_A(E\otimes_{\rA(R)}\qOmega_{\rA(R)}^{n+1}))
	}
	\]
	where the horizontal quasi-isomorphisms are given by Proposition~\ref{P:cohomology-vanishing} (more precisely by its proof). One deduces thus a quasi-isomorphism 
	\begin{equation}\label{E:double complex 2}
		(E\otimes_{\rA(R)}\qOmega_{\rA(R)}^{\bullet}, \Theta^{\bullet})\xra{\sim} \underline s(\scrC^{\bullet,\bullet}).
	\end{equation}
	Statement (1) of Theorem~\ref{T:crystal-qHiggs} now follows immediately by combining Lemma~\ref{L:CA-complex-prismatic},  \eqref{E:double complex 1}  and \eqref{E:double complex 2}.
	
	(2) We will write $\scrC^{\bullet, \bullet}(\calE)$ for $\scrC^{\bullet, \bullet}$  to emphasize its functorial dependence on $\calE$. Recall that $\qOmega^n_{D_{\rA^m(R)}/\rA^m(R)}$ is a free $D_{\rA^m(\rR)}$-module with basis $\{\frac{d u_{i_1}}{u_{i_1}}\wedge \cdots \wedge \frac{d u_{i_n}}{u_{i_n}}\colon 1\leqslant i_1<\cdots <i_n\leqslant d\}$. We define $d^n\phi_{D_{\rA^m(R)}}\colon \qOmega^n_{D_{\rA^m(R)}/\rA^m(R)}\to \qOmega^n_{D_{\rA^m(R)}/\rA^m(R)}$ by 
	\[
	d^n\phi_{D_{\rA^m(R)}}\bigg(f \frac{d u_{i_1}}{u_{i_1}}\wedge \cdots \wedge \frac{d u_{i_n}}{u_{i_n}}\bigg)=[p]_q^n \phi_{D_{\rA^m(R)}}(f)\frac{d u_{i_1}}{u_{i_1}}\wedge \cdots \wedge \frac{d u_{i_n}}{u_{i_n}},
	\]
	for all $f\in D_{\rA^m(R)}$. If $\iota\colon E\to \phi_{\rA^*(R)}E$ denotes the canonical map $x\mapsto x\otimes1$, then the morphisms
	\[
	\scrC^{m,n}(\calE)=E\widehat\otimes_{\rA(R)}\qOmega^n_{D_{\rA^m(R)}/\rA^m(R)}\xra{\iota\otimes d^n\phi_{D_{\rA^m(R)}}} \scrC^{m,n}(\phi_{\cO_{\prism}}^*\calE)=(\phi_{\rA(R)}^*E)\widehat\otimes_{\rA(R)}\qOmega^n_{D_{\rA^m(R)}/\rA^m(R)}
	\]
	define a Frobenius linear map  of double complexes
	\[
	\scrC^{\bullet,\bullet}(\iota_{\calE})\colon\quad  \scrC^{\bullet, \bullet}(\calE)\to \scrC^{\bullet, \bullet}(\phi_{\calO_{\prism}}^*\calE).
	\]
	Then there exists a commutative diagram of complexes 
	\[
	\xymatrix{
		\CA(\urA^{\bullet}(R),\calE)\ar[r]^-{ \simeq}\ar[d] &\underline {s}(\scrC^{\bullet, \bullet}(\calE))\ar[d]^{\underline{s}\scrC^{\bullet,\bullet}(\iota_{\calE})} & (E\otimes_{\rA(R)}\qOmega^{\bullet}_{\rA(R)}, \Theta^\bullet) \ar[l]_-{\simeq}\ar[d]^{\eqref{E:phi-functorial-deRham}}\\
		\CA(\urA^{\bullet}(R), \phi_{\calO_{\prism}}^*(\calE))\ar[r]^-{\simeq}
		&\underline {s}(\scrC^{\bullet, \bullet}(\phi_{\calO_{\prism}}^*\calE)) & \big((\phi_{\rA(R)}^*E)\otimes_{\rA(R)}\qOmega^{\bullet}_{\rA(R)}, \Theta^{\phi, \bullet}\big)\ar[l]_-{ \simeq}
	}
	\]
	where the left vertical arrow is the  functorial map induced by the Frobenius linear  map of crystals $\iota_{\calE}\colon \calE\to \phi_{\calO_{\prism}}^*(\calE)$, the left two horizontal maps are  \eqref{E:double complex 1}, and the right ones are \eqref{E:double complex 2}. 
	
	On the other hand, the Frobenius structure $\phi_{\calE}$ on $\calE$ induces  a isomorphism of complexes
	\[
	\underline{s}(\scrC^{\bullet, \bullet}(\phi_{\calO_{\prism}}^*\calE))[\frac{1}{[p]_q}]\xra{\sim} \underline{s}(\scrC^{\bullet, \bullet}(\calE))[\frac{1}{[p]_q}]
	\]
	which is compatible with the  functorial isomorphism of \v{C}ech--Alexander complexes via \eqref{E:double complex 1},  and with the functorial isomorphism of de Rham  quasi-isomorphisms via \eqref{E:double complex 2}. Therefore, it follows that, via the isomorphism in Theorem~\ref{T:crystal-qHiggs}(1),   both the upper and bottom horizontal arrows in  \ref{T:crystal-qHiggs}(2) are represented by 
	\[
	\phi^*\bigg(\underline{s}\scrC^{\bullet,\bullet}(\calE)\bigg)[\frac{1}{[p]_q}]\xra{\underline{s}\scrC^{\bullet,\bullet}(\iota_{\calE})} \underline{s}(\scrC^{\bullet, \bullet}(\phi_{\calO_{\prism}}^*\calE))[\frac{1}{[p]_q}]\xra{\sim} \underline{s}(\scrC^{\bullet, \bullet}(\calE))[\frac{1}{[p]_q}].
	\]
	This finishes the proof of Theorem~ \ref{T:crystal-qHiggs}(2).

\end{proof}

\section{Cohomology of log prismatic crystals and Galois cohomology}\label{Section4}

We will keep the same setup and notation as in \S \ref{S:setup of notation}. Namely,  $X=\Spf(R)$ is affine small semistable $p$-adic formal scheme over $\Spf(\cO_C)$ with  framing $\square \colon R^{\square}\to R$ as in \eqref{E:framing-over-C},
 and $\uX=(X,\calM_X)$ is the associated log formal scheme with the canonical log structure. 
 Recall that $\uX$ is attached to the prelog ring $\underline R:=(R, M_R\xra{\alpha_R}R)$ with $\alpha_R$ defined in \eqref{E:local-chart-O_C}. 

\subsection{Pefectoid cover}
\if false
Consider the integral monoid 
\[
M_{R_{\infty}}:=M_{\cO_{C}^{\flat}}\sqcup_{\NN[\frac 1p ]}\NN[\frac1p]^{r+1},
\]
where $\NN[\frac1p]\to M_{\cO_{C}^{\flat}}=\cO_{C}^{\flat}\backslash \{0\}$ sends $p^{-n}$ to $[p^{s/p^n, \flat}]$, and $\NN[\frac1p]\to \NN[\frac1p]^{r+1}$ is the diagonal map. There exists a natural morphism of monoids 
\[
M_{R}=M_{\cO_{C}^{\flat}}\sqcup_{\NN}\NN^{r+1}\to M_{R_{\infty}}
\]
induced by the identity map of $M_{\cO^{\flat}_C}$ and the natural inclusion $\NN^{r+1}\hra \NN[\frac1p]^{r+1}$. 

\begin{lemma}\label{L:flatness-monoid}
	The morphism of monoids $M_R\to M_{R_{\infty}}$  is  flat in the sense of \cite{Bhatt}*{Def. 4.8}.
\end{lemma}

\begin{proof}
	It is clear that $M_R\to M_{R_{
			\infty}}$ is the pushout of $\NN^{r+1}\to \NN[1/p]^{ r+1}$ along  the map $\NN^{r}\to M_{R}$. It suffices to show that $\NN^{r+1}\to \NN[1/p]^{r+1}$ is  flat. But we have 
	\begin{align*}
		\ZZ[\NN^{r+1}]&\cong \ZZ[T_0,\cdots, T_r]\\
		\ZZ[\NN[\frac1p]^{r+1}]&\cong \bigcup_{n\geqslant 1}\ZZ[T_0^{1/p^n},\cdots, T_r^{1/p^n}],
	\end{align*}
	and it is clear that $\ZZ[\NN^{r+1}]\to \ZZ[\NN[\frac1p]^{r+1}]$ is faithfully flat. By \cite[Prop. 4.9]{Bhatt}, the morphism of monoids $\NN^{r+1}\to \NN[\frac1p]^{r+1}$ is flat. 
	
\end{proof}

\fi

For each integer $n\geqslant0$, we consider the $R^{\square}$-algebra
\[
R^{\square}_n:=\cO_{C}\langle T_0^{1/p^n}, \cdots, T_r^{1/p^{n}}, T_{r+1}^{\pm 1/p^n}, \cdots, T_d^{\pm1/p^n}\rangle/ (T_0^{1/p^n}\cdots T_r^{1/p^{n}}-\pi^{1/p^n}), 
\]
and we put 
\[R_{\infty}^{\square}:=\big(\varinjlim_{n} R^{\square}_n\big)^{\wedge},\]
where the completion is for the  $p$-adic topology. 	Explicitly, one has 
\begin{align}\label{E:structure-R-infty-square}
	R_{\infty}^{\square}&=\widehat{\bigoplus}_{\substack{(a_0, \cdots, a_d)\in \ZZ[\frac{1}{p}]_{\geqslant0}^{\oplus (r+1)}\oplus \ZZ[\frac{1}{p}]^{\oplus (d-r)}\\ a_i=0 \text{ for some }0\leqslant i\leqslant r}}\;\cO_C\cdot T_0^{a_0}\cdots T_d^{a_d},\\
	&= \widehat \bigoplus_{\substack{(a_0,\cdots, a_d)\in \big(\ZZ[\frac{1}{p}]\cap[0,1)\big)^{d}\\ a_i=0 \text{ for some $0\leqslant i\leqslant r$}}} R^{\square} \cdot T_0^{a_0}\cdots T_d^{a_d} \nonumber
\end{align}
and via the canonical injection    $R^{\square}_n\to R^{\square}_{\infty}$, $R^{\square}_n$ is identified with the subring  of $R^{\square}_{\infty}$ consisting of the direct summands $\cO_CT_0^{a_0}\cdots T_d^{a_d}$ with $a_i\in \frac{1}{p^n}\ZZ$ for all $0\leqslant i\leqslant n$.
We put 
\[
R_n:=R^{\square}_n\widehat\otimes_{R^{\square}}R, \quad 
R_{\infty}:=R_{\infty}^{\square}\widehat{\otimes}_{R^{\square}}R.
\]
Then we have  a decomposition
\begin{equation}\label{E:decomposition-R-infty}
	R_{\infty}=\widehat\bigoplus_{\substack{(a_0, \cdots, a_d)\in \big(\ZZ[1/p]\cap [0,1)\big)^{d}\\ a_i=0\text{ for some }0\leqslant i\leqslant r}}\;R\cdot T_0^{a_0}\cdots T_d^{a_d}. 
\end{equation}

\subsection{$\Gamma$-action} 	There exists a natural $R^{\square}$-linear action by $\Gamma=\bigoplus_{i=1}^d \ZZ_p\gamma_i$ on $R_{\infty}^{\square}$, and hence on $R_{\infty}$ by scalar extension, given by 

\[
\gamma_i(T_j^{1/p^n})=\begin{cases}\zeta_{p^n} T_j^{1/p^n} &\text{if }j=i,\\
	\zeta_{p^n}^{-1}	T_0^{1/p^n} & \text{if $j=0$ and $1\leqslant i\leqslant r$},\\
	T_j^{1/p^n} & \text{otherwise}.
\end{cases}
\]
for all $1\leqslant i\leqslant d$,  $0\leqslant j\leqslant d$ and $n\geqslant1$. 

From the explicit description of $R_n^{\square}$, we see that $R_n^{\square}=(R_n^{\square}[1/p])^{\circ}$, the subring of power bounded elements in $R^{\square}_{n}[1/p]$. Hence, we have also 
$R_{n}=(R_{n}[1/p])^{\circ}$ for every  $n\geqslant 1$ and $R_{\infty}=(R_{\infty}[1/p])^{\circ}$. It is clear that $\Spa(R_n^{\square}[1/p], R^{\square}_n)\to \Spa(R^{\square}[1/p], R^{\square})$ is an \'etale  Galois cover of adic spaces over $\Spa(C, \cO_C)$ with  group $\Gamma/p^n\Gamma$ for each integer $n\geqslant 0$. 
It follows that  
the affinoid perfectoid space 
\[
X_{\infty,\eta}:=\Spa(R_{\infty}[1/p], R_{\infty})\sim \varprojlim_n \Spa(R_n[1/p], R_n)
\] 
is then a pro\'etale  $\Gamma$-cover of the adic generic fiber $X_{\eta}=\Spa(R[1/p], R)$. With the language of diamonds, we have 
\begin{equation}\label{E:R-infty-Gamma-torsor}
	X_{\eta}^{ \diamond}=X_{\infty, \eta}^{ \diamond}/\Gamma
\end{equation}

\subsection{Decomposition of $\Ainf(R_{\infty})$}

Consider the perfect prism  $(\Ainf(R_{\infty}), ([p]_q))$ corresponding to $R_{\infty}$, which is a natural object of non-logarithmic perfect prismatic site $X^{\perf}_{\prism}$.
%

Let $T_i^{\flat}:=(T_i, T_i^{1/p}, \cdots)\in R_{\infty}^{\square, \flat}$ for $0\leqslant i\leqslant d$, and $U_i:=[T_i^{\flat}]\in \Ainf(R^{\square}_{\infty})$. Then analogous to \eqref{E:structure-R-infty-square}, we have 
\[
\Ainf(R^{\square}_{\infty})=\widehat{\bigoplus}_{\substack{(a_0, \cdots, a_d)\in \ZZ[\frac{1}{p}]_{\geqslant0}^{\oplus (r+1)}\oplus \ZZ[\frac{1}{p}]^{\oplus (d-r)}\\ a_i=0 \text{ for some }0\leqslant i\leqslant r}}\;A_{\inf}\cdot U_0^{a_0}\cdots U_d^{a_d}. 
\]
The $\Gamma$-action on $R_{\infty}^{\square}$ induces a natural $\Gamma$-action on $\Ainf(R^{\square}_{\infty})$: 
\[
\gamma_i(U_j)=\begin{cases}
	q U_i &\text{if }i=j,\\
	q^{-1}U_0 &\text{if $j=0$, $1\leqslant i\leqslant r$}, \\
	U_j & \text{otherwise}.
\end{cases}
\]
for all $0\leqslant i,j\leqslant d$.

Recall  $\rA(R^{\square})$ and the object  $\urA(R)=(\rA(R), ([p]_q), M_{R})^a$ of $\uX_{\prism}$  defined  in \S~\ref{S:frames}, and there exists a $\Gamma$-action on $\urA(R)$. 
We have an injective  map of $\delta$-$A_{\inf}$-algebras 
$\rA(R^{\square})\to \Ainf(R^{\square}_{\infty})$ given  by 
$
T_i\mapsto U_i^p
$ for $0\leqslant i\leqslant d$ compatible with the $\Gamma$-action. We will always identify $\rA(R^{\square})$ with its image in $\Ainf(R^{\square}_{\infty})$. 

By the equivalence of sites $\uX_{\prism}^{\perf}\simeq X_{\prism}^{\perf}$ (Prop.~\ref{P:equiv-perfect-log-site}), there exists a unique $\delta_{\log}$-structure $M_{\Spf(\Ainf(R_{\infty}))}$ on $\Spf(\Ainf(R_{\infty}))$ such that 
$${\Ainf}(\uR_{\infty}):=(\Ainf(R_{\infty}), ([p]_q), M_{\Spf(\Ainf(R_{\infty}))})$$ becomes an object of $\uX^{\perf}_{\prism}$. The map $\rA(R)\to \Ainf(R_{\infty})$ extends uniquely to a morphism $\urA(R)\to {\Ainf}(\uR_{\infty})$ in $\uX^{\opp}_{\prism}$.

Now note that the canonical map
\[
\Ainf(R^{\square}_{\infty})\widehat\otimes_{\rA(R^{\square})}\rA(R)\xra{\sim}\Ainf(R_{\infty})
\]
is an isomorphism. 
Indeed, this is clearly true after   reduction modulo $[p]_q$, and the claim follows from the fact that $\rA(R^{\square})\to \rA(R)$ is $(p,[p]_q)$-completely \'etale. 
Put 
\begin{equation}\label{E:Xi-p}
\Xi_p:=\bigg\{(a_0,\cdots, a_d)\in \big(\ZZ[1/p]\cap [0,p)\big)^d\,\bigg|\, \text{$a_i=0$ for some $0\leqslant i\leqslant r$}\bigg\}.
\end{equation}
Then using the above explicit description of $\Ainf(R^{\square}_{\infty})$, one gets 
a decomposition
\begin{equation}\label{E:decomposition-Ainf-R}
	\Ainf(R_{\infty})=\widehat\bigoplus_{a\in \Xi_p}\;\rA(R)\cdot U_0^{a_0}\cdots U_d^{a_d}
\end{equation}
such that each direct summand is stable  under the $\Gamma$-action.

\begin{lemma}\label{L:final-cover-perfect-prismatic}
(1) The prism $(\Ainf(R_{\infty}), ([p]_q))$ is  isomorphic to the perfection of $(\rA(R), ([p]_q))$.

(2) The perfect prism $(\Ainf(R_{\infty}), ([p]_q))$  is a strong covering object of $X^{\perf}_{\prism}$, i.e. for every object $(B, ([p]_q))$ of $X^{\perf}_{\prism}$, the coproduct $$(D_B^{\perf}, ([p]_q)):=(\Ainf(R_{\infty}), ([p]_q))\coprod (B,([p]_q))$$  exists in the category $X^{\perf, \opp}_{\prism}$, and the canonical map $B\to D_B^{\perf}$ is $(p, [p]_q)$-completely faithfully flat. In particular, $(\Ainf(R_{\infty}), ([p]_q))$ is a flat cover of the final object of $\Shv(X^{\perf}_{\prism})$.

\end{lemma}
\begin{proof}
	(1) We have to show that $\Ainf(R_{\infty})$ is isomorphic to the $(p,[p]_q)$-adic completion of $\colim(\rA(R)\xra{\phi}\rA(R)\xra{\phi}\cdots)$. This follows easily from an explicit computation if $R=R^{\square}$. The case for general $R$ follows from the fact that $\rA(R^{\square})\to \rA(R)$ is $(p, [p]_q)$-completely \'etale and the isomorphism $\Ainf(R^{\square}_{\infty})\widehat\otimes_{\rA(R^{\square})}\rA(R)\xra{\sim}\Ainf(R_{\infty})$.
	
	(2) Let $(B, ([p]_q), M_{\Spf(B)})$ be the object of $\uX_{\prism}^{\perf}$ corresponding to $(B, ([p]_q))$ via the equivalence of categories $\uX^{\perf}_{\prism}\simeq X^{\perf}_{\prism}$ (Prop. \ref{P:equiv-perfect-log-site}).  Let $(D_B, ([p]_q), M_{\Spf(D_B)})$ denote the coproduct of $\urA(R)$ and $(B, ([p]_q), M_{\Spf(B)})$ in $\uX^{\opp}_{\prism}$ given by  Lemma~\ref{L:cover-final-object}, and $(D_B^{\perf}, ([p]_q))$ be the perfection of the prism $(D_B, ([p]_q))$. By the universal property of the perfection of prisms and part (1), the morphism $(\rA(R), ([p]_q))\to (D_B^{\perf}, ([p]_q))$ factors through  $(\rA(R), ([p]_q))^{\perf}=(\Ainf(R_{\inf}), ([p]_q))$  and it is easy to see that $(D_B^{\perf}, ([p]_q))$ is the coproduct of $(\Ainf(R_{\inf}), ([p]_q))$ and $(B, ([p]_q))$ in $X^{\perf}_{\prism}$. 
	Finally, since $B\to D_B$ is $(p, [p]_q)$-completely faithfully flat, so is $B\to D_B^{\perf}$ since $(p, [p]_q)$-complete  faithful flatness is preserved by  colimits.

\end{proof}


\subsection{Cohomology of prismatic crystals and Galois cohomology}\label{S:coh-prismatic-Galois}

Let $\calE$ be an object of  $\widehat{\CR}(\uX_{\prism},\cO_{\prism})^{\mu-tf}$, and $(E,(\theta_{i})_{1\leqslant i\leqslant d})$ be the $q$-Higgs module attached to $\calE$ via the functor \eqref{E:strat-Higgs}. Recall that there exists a semi-linear $\Gamma$-action on $E$ given by $\rho_{E}(\gamma_i)=\id_E+\mu\theta_i$ (cf. \S~\ref{S:simplicial-A(R)}). 
Put $E_{\infty}:=\calE(\Ainf(\uR_{\infty}))$. Then the crystal property of $\calE$ implies that $E_{\infty}=E\widehat\otimes_{\rA(R)} \Ainf(R_{\infty})$, which is equipped with a natural induced semi-linear $\Gamma$-action. 

Let $W(\gothm_{C^{\flat}})$ be the kernel of the surjection $A_{\inf}\to W(k)$ induced by $\cO_{C}^{\flat}\to k$. Then $W(\gothm_{C^{\flat}})$ is generated by $[x]$ with $x\in \cO^{\flat}_{C}$.  An $A_{\inf}$-module $M$ is called \emph{almost zero} if it is killed by  $W(\gothm_{C^{\flat}})$. A morphism of $A_{\inf}$-modules $f: M\to N$ is an \emph{almost isomorphism} if both $\ker(f)$ and $\coker(f)$ are almost zero; in that case, we will write $M\simeq^a N$. 

\begin{theorem}\label{T:coh-crystal-Galois-coh}
Under the above notation, the following holds:
	\begin{enumerate}
		\item[(1)] Recall the morphism of topoi $\varepsilon\colon \Shv(X_{\prism}^{\perf})\simeq \Shv(\uX_{\prism}^{\perf})\to \Shv(\uX_{\prism})$  defined in \eqref{E:perfect-to-nonperfect}, and put  $\calE^{\perf}:=\varepsilon^*\calE$. Then there is a canonical map in the derived category $D(A_{\inf})$
		\[
		\iota_{\calE}\colon R\Gamma(X^{\perf}_{\prism}, \calE^{\perf})\to R\Gamma(\Gamma, E_{\infty})
		\]
		such that the induced map  $H^i(\iota_{\calE})\colon H^i(X^{\perf}_{\prism}, \calE^{\perf})\simeq^a H^i(\Gamma, E_{\infty})$ is an almost isomorphism  of $A_{\inf}$-modules for all $i\in \ZZ$.
		
		\item[(2)] The composite canonical map 
		\[
		R\Gamma(\uX_{\prism}, \calE)\xra{\varepsilon^*} R\Gamma(X^{\perf}_{\prism}, \calE^{\perf})\xra{\iota_{\calE}} R\Gamma(\Gamma, E_{\infty})
		\]
		induces an isomorphism in $D(A_{\inf})$:
		\[
		R\Gamma(\uX_{\prism}, \calE)\xra{\sim}L\eta_{\mu} R\Gamma(\Gamma, E_{\infty}).
		\]
		Here, $L\eta_{\mu}$ is the d\'ecalage functor studied in \cite[\S 6]{BMS1}. In particular, the canonical map in  $D(A_{\inf})$
		\[
		R\Gamma(\uX_{\prism}, \calE[1/\mu])\xra{\sim} R\Gamma(X_{\prism}^{\perf}, \calE^{\perf}[1/\mu])
		\]
		is an isomorphism.
	\end{enumerate}
	
\end{theorem}

\begin{remark}\label{R:Ainf-coh}
	(1)	When $X$ is smooth over $\Spf(\cO_C)$, this Theorem follows from \cite[Theorem 6.1]{MT} and \cite[Theorem 13.31]{Tsuji}. We will  generalize the arguments of \emph{loc. cit.} to our situation.
	
	(2) When $\calE=\cO_{\prism}$, statement (2) of this Theorem actually shows that the $A_{\inf}$-cohomology $R\Gamma_{A_{\inf}}(X)$ defined in \cite{CK} coincides with the log prismatic cohomology $R\Gamma(\uX_{\prism}, \cO_{\prism})$.
	
	(3) We expect that $L\eta_{\mu}R\Gamma(X^{\perf}_{\prism}, \calE^{\perf})\xra{\sim} L\eta_{\mu}R\Gamma(\Gamma, E_{\infty})$ is also a quasi-isomorphism (instead of  being just an almost quasi-isomorphism).
	
	
	
\end{remark}

\begin{proof}
	(1)	Let $(\Ainf(R^{\bullet}_{\infty}),([p]_q))$ be the \v{C}ech nerve of $(\Ainf(R_{\infty}), ([p]_q))$ over the final object of $X_{\prism}^{\perf}$, i.e. for every integer $n\geqslant 0$, $(\Ainf(R_{\infty}^n), ([p]_q))$ is the $(n+1)$-fold self coproduct of $(\Ainf(R_{\infty}), ([p]_q))$ in $X^{\perf, \opp}_{\prism}$. 
	Let $R_{\infty}^{\bullet}=\Ainf(R^{\bullet}_{\infty})/([p]_q)$ denote the corresponding cosimplicial perfectoid ring over $R$. 
	 As $R_{\infty}=R_{\infty}^0$ is $(p,[p]_q)$-completely faithfully flat over $R$, Lemma~\ref{L:final-cover-perfect-prismatic} implies that $R^{n}_{\infty}$ is also $(p,[p]_q)$-completely faithfully flat over $R$ for every integer $n\geqslant 0$.  In particular, $R^n_{\infty}$ is a $p$-torsion free perfectoid ring for each $n\geqslant0$. 
	
	On the  adic generic fiber, we consider similarly the \v{C}ech nerve of the projection of diamonds  $X_{\infty, \eta}^{ \diamond}\to X_{\eta}^{\diamond}$.  
	In view of \eqref{E:R-infty-Gamma-torsor}, we have 
	\begin{align*}
		X^{ n}_{\infty,\eta}&:=\underbrace{X^{\ad}_{\infty, \eta}\times_{X_{\eta}^{\diamond}}\cdots \times_{X^{ \diamond}_{\eta}}X_{\infty,\eta}^{ \diamond}}_{\text{$n+1$-fold product}}\\
		&\cong X^{ \diamond}_{\infty, \eta}\times \underline{\Gamma}^n=\Spa(C^n(\Gamma, R_{\infty}[1/p]),C^n(\Gamma, R_{\infty}))^{\diamond},
	\end{align*}
	where $C^n(\Gamma, R_{\infty})$ denotes the continuous function on $\Gamma^n$ with values in $R_{\infty}$, and the simple complex attached to the cosimplicial object   $C^{\bullet}(\Gamma, R_{\infty})$  is nothing but the standard inhomogeneous cochain complex of $\Gamma$ with values in $R_{\infty}$.

	Using  the equivalence between  perfectoid rings  and  perfect prisms, we see that $R_{\infty}^{n}$ is the initial (integral) perfectoid rings over $R$ with $n$ morphisms (as $R$-algebras) from $R_{\infty}$. Hence, $R_{\infty}^n[1/p]$ is the initial perfectoid  Tate $C$-algebra equipped with $n$ morphisms of $R[1/p]$-algebras from $R_{\infty}[1/p]$. 
	 It follows thus that there exists 
	a morphism of cosimplicial perfectoids rings over $R$ 
	\[
	\iota^{\bullet}\colon R_{\infty}^{\bullet}\to C^{\bullet}(\Gamma, R_{\infty})
	\]
	which induces an isomorphism 
	$R_{\infty}^{\bullet}[1/p]\cong C^{\bullet}(\Gamma, R_{\infty}[1/p])$. Hence, we have   $R_{\infty}^n\subset C^n(\Gamma, R_{\infty}[1/p])^{\circ}=C^n(\Gamma, R_{\infty})$, which is an almost isomorphism of $\cO_{C}$-modulles by \cite[Lemma 3.21]{BMS1}. 
	Applying the $\Ainf$-functor, we get a map of cosimplicial perfect prisms over $(A_{\inf}, ([p]_q))$:
	$$\Ainf(\iota^{\bullet})\colon\quad  (\Ainf(R_{\infty}^{\bullet}), ([p]_q))\xra{\sim}^a(\Ainf(C^\bullet(\Gamma, R_{\infty})), ([p]_q))$$
	such that for each integer $n\geqslant 0$ the underlying map of $A_{\inf}$-algebras $\Ainf(R_{\infty}^{n})\to \Ainf(C^n(\Gamma, R_{\infty}))$ is an almost isomorphism. 
	By evaluating $\calE^{\perf}$,  we get then a morphism of cosimplicial $A_{\inf}$-modules:
	\[
	\iota_{\calE}^{\bullet}\colon	\calE^{\perf}(R_{\infty}^{\bullet})\to \calE^{\perf}(C^\bullet(\Gamma, R_{\infty})),
	\]
	where for a (cosimplicial) integral perfectoid ring $S$ over $R$ we have written $\calE^{\perf}(S):=\calE^{\perf}(\Ainf(S), ([p]_q))$ to simply the notation. 
Then the crystal property of $\calE^{\perf}$ implies that 
	\[
	\calE^{\perf}(C^\bullet(\Gamma, R_{\infty}))\cong \calE^{\perf}(R^{\bullet}_{\infty})\widehat\otimes_{\Ainf(R_{\infty}^{\bullet} )}C^{\bullet}(\Gamma, \Ainf(R_{\infty}))
	\] 
	By the discussion above, we see that  $\iota_{\calE}^n$ is an almost ismorphism of $A_{\inf}$-modules for each integer $n$. In particular, $\iota_{\calE}^\bullet$  induces  a morphism 
	\[
	\iota_{\calE}\colon 
	\underline s(\calE^{\perf}(R^{\bullet}_{\infty}))\to \underline s(\calE(C^\bullet(\Gamma, R_{\infty}))
	\]
	between the associated simple complexes such that $H^i(\iota_{\calE})$ is an almost isomoprhism for all $i\in \ZZ$.

	On the other hand hand,  as $(\Ainf(R_{\infty}), ([p]_q))$ is a strong  covering object  of  $X^{\perf}_{\prism}$ (Lemma~\ref{L:final-cover-perfect-prismatic}), similar arguments as Lemma~\ref{L:CA-complex-prismatic} show that $R\Gamma(X_{\prism}^{\perf}, \calE^{\perf})$ is computed by the simple complex attached to $\calE^{\perf}(R_{\infty}^{\bullet})$. On the other hand, by the crystal property of $\calE$, it is easy to see that the simple complex attached to $\calE(C^\bullet(\Gamma, R_{\infty}))$ is nothing but  $C^{\bullet}(\Gamma, E_{\infty})$, the standard inhomogeneous continuous cochain complex that computes $R\Gamma(\Gamma, E_{\infty})$.  Statement (1) now follows immediately.

	(2) For every $\underline a\in \Xi_p$ (see \eqref{E:Xi-p}), let $\chi_{\underline a}\colon \Gamma\to A_{\inf}^{\times}$ be the character defined by  
	\[
	\chi_{\underline a}(\gamma_i)=\begin{cases}q^{a_i-a_0} &\text{if }1\leqslant i\leqslant r,\\
		q^{a_i}, &\text{if }i>r. 
	\end{cases}
	\] 
	for $1\leqslant i\leqslant d$. Then according to \eqref{E:decomposition-Ainf-R}, we have a decomposition 
	$
	E_{\infty}=E\oplus F_{\infty}
	$ 
	stable under the $\Gamma$-action, where 
	\[
	F_{\infty}=\widehat\bigoplus_{\underline a\in \Xi_p^*} E(\chi_{\underline a}),
	\]
	with $\Xi_p^*:=\Xi_p\backslash \{(0,\cdots, 0)\}$.
	According to Theorem~\ref{T:crystal-qHiggs} and Proposition~\ref{P:module with connection}(2), we have isomorphisms in $D(A_{\inf})$
	\[
	R\Gamma(\uX_{\prism}, \calE)\xra{\sim} E\otimes_{\rA(R)}\qOmega^{\bullet}_{\rA(R)}\xra{\sim}L\eta_{\mu} R\Gamma(\Gamma, E).
	\]
	On the other hand,  by \cite[Lemma 6.20]{BMS1}, one has 
	\[
	L\eta_{\mu}R\Gamma(\Gamma, F_{\infty})\cong \widehat{\bigoplus}_{\underline a\in \Xi^*_p}L\eta_{\mu}R\Gamma(\Gamma, E(\chi_{\underline a})).
	\]
	Therefore, to finish the proof of Theorem~\ref{T:coh-crystal-Galois-coh}, it suffices to show that 
	\[
	L\eta_{\mu}R\Gamma(\Gamma, E(\chi_{\underline a}))=0
	\]
	for all $\underline a\in \Xi_p^*$. According to \cite[Lemma 6.4]{BMS1}, one has  that  
	$$H^i(L\eta_{\mu}R\Gamma(\Gamma, E(\chi_{\underline a})))=H^i(\Gamma, E(\chi_{\underline a}))/H^i(\Gamma, E(\chi_{\underline a}))[\mu], \quad \forall i\in \ZZ.
	$$
	We are reduced  to showing that $H^i(\Gamma, E(\chi_{\underline a}))$ is annihilated by $\mu=q-1$ for all $\underline a\in \Xi^*_p$. It is well known that $R\Gamma(\Gamma, E(\chi_{\underline a}))$ is computed by the Koszul complex (cf. \cite[Lemma~7.3]{BMS1}):
	\[
	\Kos^{\bullet}(E(\chi_{\underline a}); (\gamma_i-1)_{1\leqslant i\leqslant d}):=\bigg(E(\chi_{\underline a})\xra{(\gamma_1-1,\cdots, \gamma_d-1)} \bigoplus_{1\leqslant i\leqslant d} E(\chi_{\underline a})\to \bigoplus_{1\leqslant i_1<i_2\leqslant d}E(\chi_{\underline a})\to \cdots   \bigg).
	\]
	concentrated in degrees $[0,d]$. 
	
	First, we observe that, for any $\underline a\in \Xi^*_p$, there exists $1\leqslant i_0\leqslant d$ such that $\chi_{\underline a}(\gamma_{i_0})=q^{{m}/{p^n}}=\phi^{-n}(q)^m$  for some integers $n\geqslant 0$ and $m$ with $p\nmid m$. Indeed, by the definition of $\Xi^*_p$, there is an integer $i$ with $a_i\in \ZZ[1/p]\cap (0,p)$. If $i>r$, then we take $i_0=i$. If $1\leqslant i\leqslant r$, then   $\chi_{\underline a}(\gamma_i)=q^{a_i-a_0}$ and we can take $i_0=i$ if $a_0\neq a_i$. Otherwise, there exists  $j$ with $1\leqslant j\leqslant r$ with $a_j=0$. Then we can take $i_0=j$. 
	
	The complex $\Kos^{\bullet}(E(\chi_{\underline a}), \gamma_{i_0}-1)$ is given by 
	\[
	E\xra{(1+\mu \theta_{i_0})q^{{m}/{p^n}}-1} E
	\]
	concentrated in degree $0,1$. Now we note that 
	\[
	(1+\mu\theta_{i_0})q^{{m}/{p^n}}-1=\phi^{-n}(\mu)\bigg(\frac{q^{m/p^n}-1}{\phi^{-n}(\mu)}+q^{m/p^n}\frac{\phi(\mu)}{\phi^{-n}(\mu)}\theta_{i_0}\bigg),
	\]
	where  $\frac{q^{m/p^n}-1}{\phi^{-n}(\mu)}$ is invertible in $A_{\inf}$ because  the image 
	$$(\theta_{\cO_C}\circ\phi^{-1})\bigg(\frac{q^{m/p^n}-1}{\phi^{-n}(\mu)}\bigg)=\sum_{i=0}^{m-1}\zeta_{p^{n+1}}^i$$ is invertible in $\cO_{C}$. On the other hand, $\theta_{i_0}$ is topologically quasi-nilpotent by Proposition~\ref{P:top-nilpotence}. 
	Hence,  $\frac{q^{m/p^n}-1}{\phi^{-n}(\mu)}+q^{m/p^n}\frac{\phi(\mu)}{\phi^{-n}(\mu)}\theta_{i_0}$ is an automorphism of $E$, and we have 
	\[
	H^i(\Kos^{\bullet}(E(\chi_{\underline a}), \gamma_{i_0}-1))=\begin{cases}0, &\text{if }i\neq 1,\\
		E/\phi^{-n}(\mu)E, & \text{if } i=1.
	\end{cases}
	\]
	Here, the vanishing of $H^0$ uses the assumption that $E$ is $\mu$-torsion free. 
	
	We relabel the $(\gamma_i)_{1\leqslant i\leqslant d}$ as $(\gamma_{i_j})_{0\leqslant j\leqslant d-1}$, and note that $\Kos^{\bullet}(E(\chi_{\underline a}), (\gamma_{i_j})_{0\leqslant j\leqslant k})$ is the mapping cone of 
	\[
	\Kos^{\bullet}(E(\chi_{\underline a}), (\gamma_{i_j})_{0\leqslant j\leqslant k-1})\xra{\gamma_{i_k}-1}\Kos^{\bullet}(E(\chi_{\underline a}), (\gamma_{i_j})_{0\leqslant j\leqslant k-1})
	\] 
	for all $k\geqslant 1$. 
	It follows  easily by induction on $k$ that $H^i(\Kos^{\bullet}(E(\chi_{\underline a}), (\gamma_{i_j})_{0\leqslant j\leqslant k}))$  is annihilated by $q^{1/p^n}-1=\phi^{-n}(\mu)$ for all $i\in \ZZ$. In particular, $H^i(\Gamma, E(\chi_{\underline a}))$ is annihilated by $\mu$ for all $i\in \ZZ$ and $\underline a\in \Xi^*_p$.

\end{proof}

\section{Prismatic-\'etale comparison theorem on semistable formal schemes over $\Spf(\cO_C)$}\label{Section5}


We will keep the notation in \S~\ref{S:setup of notation}. 
Let $X$ be a semi-stable $p$-adic formal scheme over $\Spf(\cO_C)$ in the sense of Definition~\ref{D:semi-stable}. 
Let $\uX=(X, \calM_X)$ be the associated log formal scheme with canonical log structure. 
Let $X_{\eta}$ be the adic generic fiber, and $X_{\eta, v}^{\diamond}$ be the $v$-site of the diamond associated to $X_{\eta}$. This is the category of perfectoid spaces over $X_{\eta}$ equipped with the $v$-topology  \cite[Def. 8.1]{Sch}.

\begin{lemma}\label{L:perfect-prismatic-v}
 The functor 
\[\kappa_{\sharp}\colon X_{\eta, v}\to X_{\prism}^{\perf}
\]
that sends an affinoid perfectoid space  $\Spa(S, S^+)$ over $X_{\eta}$ to $(\Ainf(S^+), ([p]_q))$ is cocontinuous. 
\end{lemma}

\begin{proof}Let $(\Ainf(S^+), ([p]_q))\to (B, ([p]_q))$ be a flat cover in $X_{\prism}^{\perf}$. We need to prove that there exists a $v$-cover  of perfectoid spaces $\Spa(S',S'^+)\to \Spa(S, S^+)$ such that $\Ainf(S^+)\to \Ainf(S'^+)$ factors through $B$.  Put $\bar B=B/([p]_q)$ which is an integral perfectoid ring $p$-completely faithfully flat over $S^+$. In particular, $\bar B$ is $p$-torsion free. Let $\bar B^+$ be the integral closure of $\bar B$ in $\bar B[\frac{1}{p}]$. Then $\Spa(\bar B[1/p], \bar B^+)\to \Spa(S, S^+)$ is a $v$-cover that satisfies the requirement. 
\end{proof}
Therefore, the functor $\kappa_{\sharp}$ induces a morphism of topoi
\[
\kappa\colon \Shv(X_{\eta, v}^{\diamond})\to \Shv(X_{\prism}^{\perf}).
\]
The direct image functor $\kappa_*$ can be described as follows. Let $\calF$ be an object of $\Shv(X_{\eta}^{\diamond})$, and  $(B, ([p]_q))$ be an object of  $X^{\perf}_{\prism}$. 
Write $\bar B=B/([p]_q)$.  Let $\bar B^{\tf}$ be the maximal $p$-torsion free quotient of $\bar B$ and $\bar B^{\tf, +}$ be the integral closure of $\bar B^{\tf}$ in $\bar B^{\tf}[1/p]$. Then we have 
\begin{equation}\label{E:image-direct-kappa}
\kappa_*\calF(B,([p]_q))=\calF(\bar B^{\tf}[1/p], \bar B^{\tf, +}).
\end{equation}

Let $X_{\eta, \proet}$ be the pro-\'etale site of $X_{\eta}$, and $X_{\et}$ be the  small \'etale site of the formal scheme $X$. Then we have a commutative diagram of morphism of topoi
\begin{equation}\label{E:commutatiive-diag-topoi}
\xymatrix{\Shv(X_{\eta, v})\ar[r]^{\kappa}\ar@/^2pc/[rr]^{\gamma}\ar[d]_{\beta} & \Shv(X_{\prism}^{\perf})\ar[r]^{\varepsilon} \ar[rd]^{u^{\perf}} & \Shv(\uX_{\prism})\ar[d]^-{u}\\
\Shv(X_{\eta, \proet})\ar[rr]^-{\nu} && \Shv(X_{\et}),}
\end{equation}
where $\beta, u, \nu$ are natural projections, $\varepsilon$ is the moprhism induced by the natural inclusion $X^{\perf}_{\prism}\simeq \uX^{\perf}_{\prism}\to \uX_{\prism}$ as defined in  \eqref{E:perfect-to-nonperfect},  $\lambda=\varepsilon\circ \kappa$ and $u^{\perf}=u\circ \epsilon$.

Let $\widehat{\cO}_{X_{\eta}}$ (resp. $\widehat{\cO}^+_{X_{\eta}}$, $\widehat{\cO}^{\flat}_{X_{\eta}}$, $\widehat{\cO}_{X_{\eta}}^{+,\flat}$) denote the presheaf that sends an affinoid perfectoid $\Spa(S,S^+)$ over $X_{{\eta}}$ to  the ring $S$ (resp. $S^+$, $S^{\flat}$, $S^{+,\flat}$). 
Then  according to \cite[Thm. 8.7]{Sch},  $\widehat{\cO}_{X_{\eta}}$, $\widehat{\cO}^+_{X_{\eta}}$, $\widehat{\cO}^{\flat}_{X_{\eta}}$ and  $\widehat{\cO}_{X_{\eta}}^{+,\flat}$ are all $v$-sheaves. We put $\AAA_{\inf, X_{\eta}}:=W(\widehat\cO^{+,\flat}_{X_{\eta}})$.

For a morphism $M\to N$ of $A_{\inf}$-modules, we write $M\simeq^a N$ if it is an almost isomorphism in the sense of \S~\ref{S:coh-prismatic-Galois}. In particular, we write $M\simeq^a0$ if $M$ is an almost $A_{\inf}$-module, i.e. it is annihilated by $W(\gothm_{C^{\flat}})$.

\begin{lemma}\label{L:almost-iso-v-to-qsyn}
	For an object $\calE$ of $\CRhat^{\anfr}(\uX_{\prism}, \cO_{\prism})$ (Definition~\ref{D:complete-prismatic-F-crystal}), let $\calE^{\perf}:=\varepsilon^*(\calE)$ be its restriction to $X^{\perf}_{\prism}$. 
	Then  
	there exists a canonical isomorphism in the derived category of $A_{\inf}$-sheaves on $ \uX_{\prism}$ 
	\[
	R\kappa_*\kappa^{*}(\calE^{\perf}[1/\mu])\simeq \calE^{\perf}[1/\mu],
	\]
	where $\mu=[\epsilon]-1\in A_{\inf}$ is defined in \S~\ref{S:setup of notation}.
\end{lemma}

\begin{proof}
	By Definition~\ref{D:complete-prismatic-F-crystal} , $\calE[1/\mu]$ is a locally free $\cO_{\prism}[1/\mu]$-module of finite rank, hence $\calE^{\perf}[1/\mu]$ is locally free over $\cO_{\prism}^{\perf}[1/\mu]=\varepsilon^*\cO_{\prism}[1/\mu]$. 
	The projection formula implies   a canonical isomorphism 
	\[
	R\kappa_*\kappa^{*}(\calE^{\perf}[1/\mu])\simeq R\kappa_*\kappa^{*}(\cO^{\perf}_{\prism}[1/\mu])\otimes^L_{\cO^{\perf}_{\prism}[1/\mu]} \calE^{\perf}[1/\mu].
	\]
	The problem is thus reduced to the case  $\calE=\cO_{\prism}$.

	Let $(B, ([p]_q))$ be an object of $X^{\perf}_{\prism}$. By \eqref{E:image-direct-kappa}, we have 
	\[
	\kappa_*\kappa^{*}(\cO_{\prism}^{\perf})(B,([p]_q))=\Ainf(\bar B^{\tf, +}).
	\]
	Note that the kernel of the surjection $\bar B\to \bar B^{\tf}$ is annihilated by $\gothm_{C}$ (cf. \cite[\S 2.1.3]{CS}) and $\bar B^{\tf}\xra{\sim}^a \bar B^{\tf, +}$ is an almost isomorphism of $\cO_C$-modules. It follows that $B=\Ainf(\bar B)\xra{\sim}^a \Ainf(\bar B^{\tf, +})$ is  an almost  isomorphism of $A_{\inf}$-modules with respect to the ideal $W(\gothm_{C^{\flat}})$.  Therefore,  $\cO_{\prism}^{\perf}\xra{\sim}^a \kappa_*\kappa^{*}(\cO_{\prism}^{\perf})$ is an almost isomorphism of sheaves of $A_{\inf}$-modules, and in particular it is an isomorphism after inverting $\mu$. 
	
	It remains to show that $R^i\kappa_*\kappa^*(\cO^{\perf}_{\prism}[1/\mu])=0$ for any $i>0$. By \eqref{E:image-direct-kappa}, $R^i\kappa_*\kappa^*(\cO^{\perf}_{\prism}[1/\mu])$ is the sheaf on $X^{\perf}_{\prism}$ associated to  the presheaf 
	$$
	(B, ([p]_q))\mapsto H^i(\Spa(\bar B^{\tf}[1/p],\bar B^{\tf, +})_{v}, \AAA_{\inf, X_{\eta}}[1/\mu]).
	$$
	We claim that  $H^i(\Spa(\bar B^{\tf}[1/p],\bar B^{\tf, +})_{v}, \AAA_{\inf, X_{\eta}})$ is almost zero for all $i>0$, which implies  $H^i(\Spa(\bar B^{\tf}[1/p],\bar B^{\tf, +})_{v}, \AAA_{\inf, X_{\eta}}[1/\mu])=0$ because $\mu\in W(\gothm_{C^{\flat}})$.

It remains to prove  the claim. Put $Y:=\Spa(\bar B^{\tf}[1/p],\bar B^{\tf, +})$ to simply the notation.  First,  	 by induction on $n>0$,   it follows  from  \cite[Thm. 8.8]{Sch} that  $H^i(Y_{v}, W_n(\widehat{\cO}^{+, \flat}_{Y}))$ is almost zero for all  $n,i>0$. 

	Because the topos $\Shv(Y_{v})$ is replete, we have  $$\AAA_{\inf, X_{\eta}}|_Y=\AAA_{\inf, Y}=R\varprojlim_{n}W_n(\widehat\cO^{+, \flat}_{Y})$$ and   hence $R\Gamma(Y_{v}, \AAA_{\inf,Y})=R\lim_n R\Gamma(Y_{v}, W_n(\widehat\cO^{+, \flat}_{Y})$.  One has  a spectral sequence  
	\[
	E_2^{i, j}=R^i\varprojlim_n H^j(Y_{v}, W_n(\widehat\cO^{+, \flat}_{Y}))\Rightarrow H^{i+j}(X_v, \AAA_{\inf, Y})
	\]
which implies  a short exact sequence 
\[
0\to R^1\varprojlim_n H^{i-1}(Y_{v}, W_n(\widehat\cO^{+, \flat}_{Y}))\to H^i(Y_v, \AAA_{\inf, Y})\to \varprojlim_n H^i(Y_{v}, W_n(\widehat\cO^{+, \flat}_{Y}))\to 0
\]
for each integer $i>0$. Since $H^i(Y_{v}, W_n(\widehat\cO^{+, \flat}_{Y}))\simeq^a0$ for all $n,i>0$, we get immediately that $H^i(Y_v, \AAA_{\inf, Y})\simeq^a 0$ for $i\geqslant 2$. When $i=1$, we $H^0(Y_v, W_n(\widehat\cO^{+, \flat}_{Y}))=W_n(\bar B^{\tf, +, \flat})$ and  $R^1\varprojlim_n H^{0}(Y_{v}, W_n(\widehat\cO^{+, \flat}_{Y}))=0$. It follows that  $H^1(Y_v, \AAA_{\inf, Y})\simeq \varprojlim_n H^i(Y_{v}, W_n(\widehat\cO^{+, \flat}_{Y}))$ is almost zero as well.
	 
\end{proof}

\begin{proposition}\label{P:coh-prismatic-perfection}
	Let  $\calE$ be an  object of $\CRhat^{\anfr}(\underline X_{\prism}, \cO_{\prism})$,
	the canonical  map 
	\begin{equation}\label{E:crystal-adjunction}
		Ru_*(\calE)[1/\mu]\xra{\sim} Ru^{\perf}_*(\calE^{\perf}[1/\mu])
	\end{equation}
	is an isomorphism. 
\end{proposition}

\begin{proof}
	Obviously, the problem is local for the \'etale topology on $X$. We may assume that $X=\Spf(R)$ is affine small with a framing \eqref{E:framing-over-C}. As we only need the result after inverting $\mu$, up to replacing $\calE$ by its maximal $\mu$-torsion free quotient, we may assume that $\calE$ itself is $\mu$-torsion free. Then  the result is   exactly Theorem~\ref{T:coh-crystal-Galois-coh}(2).
\end{proof}


\subsection{\'Etale realization of complete prismatic $F$-crystals}\label{S:etale-realization-complete}
Recall that $\Loc(X_{\eta, \et})$ denotes the category of \'etale $\ZZ_p$-local systems over $X_{\eta}$. 
Consider now the composite functor
\begin{equation}\label{E:etale-realization-complete-F-crystal}
T\colon \CRhat^{\anfr}(\uX_{\prism}, \cO_{\prism})^{\phi=1}\to \Vect(\underline X_{\prism}, \cO_{\prism}[1/\calI]^{\wedge}_p)^{\phi=1}\simeq \Loc_{\ZZ_p}(X_{\eta, \et})
\end{equation}
where the first functor is the natural restriction, and the last equivalence is Theorem~\ref{T:etale-realization}.  
For an object  $(\calE, \phi_{\calE})\in \CRhat^{\anfr}(\uX_{\prism}, \cO_{\prism})^{\phi=1}$,  we write 
$T(\calE)=T(\calE, \phi_{\calE})$, which  can be viewed naturally as a sheaf 
on $X_{\eta,\proet}$. We call   it  the \emph{\'etale realization} of $(\calE, \phi_{
	\calE})$. 
The Frobenius structure $\phi_{\calE}\colon (\phi_{\cO_{\prism}}^*\calE)[\frac{1}{[p]_q}]\xra{\sim} \calE[\frac{1}{[p]_q}]$ induces a map 
\begin{equation}\label{E:crystal-Frob-mu}
	\calE[\frac{1}{\mu}]\to \phi^*_{\cO_{\prism}}\bigg(\calE[\frac{1}{\mu}]\bigg)=\phi^*_{\cO_{\prism}}(\calE)[\frac{1}{\phi(\mu)}]\xra{\phi_{\calE},\,\simeq} \calE[\frac{1}{\phi(\mu)}].
\end{equation}

\begin{lemma}\label{L:construction of morphism}
	Let   $(\calE, \phi_{\calE})$ be an object of   $ \CRhat^{\anfr}(\uX_{\prism}, \cO_{\prism})^{\phi=1}$. Then
	
	\begin{enumerate}
		\item[(1)]  One has a canonical isomorphism  of sheaves of $A_{\inf}$-modules on $X_{\eta,v}^{\diamond}$:
		\[
		\beta^{*}\big(T(\calE)\otimes_{\ZZ_p}\AAA_{\inf, X_{\eta}}[1/\mu]\big)\simeq \gamma^{*}(\calE[1/\mu]),
		\]
		which is equivariant under the Frobenius action in the sense that the  diagram 
		\[
		\xymatrix{\beta^{*}\big(T(\calE)\otimes_{\ZZ_p}\AAA_{\inf, X_{\eta}}[\frac{1}{\mu}]\big)\ar[rr]^-{1\otimes \phi_{\AAA_{\inf, X_{\eta}}}}\ar[d]^{\simeq} &&  \beta^{*}\big(T(\calE)\otimes_{\ZZ_p}\AAA_{\inf, X_{\eta}}[\frac{1}{\phi(\mu)}]\big)\ar[d]^{\simeq}\\
			\gamma^{*}(\iota_*\calE[\frac{1}{\mu}])\ar[rr]^-{\eqref{E:crystal-Frob-mu}} && \gamma^{*}(\calE[\frac{1}{\phi(\mu)}])}
		\]
		is commutative.
		\item[(2)] The canonical adjunction map 
		\[
		T(\calE)\otimes_{\ZZ_p}\AAA_{\inf, X_{\eta}}[1/\mu]\xra{\sim} R\beta_*\beta^{*}\big(T(\calE)\otimes_{\ZZ_p}\AAA_{\inf, X_{\eta}}[1/\mu]\big)
		\]
		is an isomorphism. 
		
	\end{enumerate}
\end{lemma}

\begin{proof}
	(1) By definition of the functors $\beta$ and  $\gamma$, it suffices to show that, for any  affinoid perfectoid space $\Spa(S,S^+)$ over $X_{\eta}$, one has an isomorphism  of $\Ainf(S^+)$-algebras
	\[
	T(\calE)\otimes_{\ZZ_p}\Ainf(S^+)[1/\mu]\simeq \calE^{\perf}(\Ainf(S^+), ([p]_q))[1/\mu]
	\]
	functorial in $\Spa(S,S^+)$. 
	The arguments are exactly the same as \cite[Lemma 3.13]{GR}. 
	
	(2) Since  any \'etale $\ZZ_p$-local system can be trivialized locally for the pro\'etale topology on $X_{\eta}$, we may reduce to the case  that $T(\calE)=\ZZ_p$. 
	By \cite[Theorem 8.7 and 8.8]{Sch},  one has an isomorphism of sheaves on $X_{\eta, \proet}$:  \[\widehat{\cO}_{X_{\eta}}^+\simeq \beta_*\beta^{*}(\widehat{\cO}_{X_{\eta}}^+)\]
	and $R^i\beta_{*}\beta^*(\widehat{\cO}_{X_{\eta}}^+)$ is almost zero for all $i>0$. Applying the $\AAA_{\inf}$-functor, we get  an isomorphism 
	\[
	\AAA_{\inf, X_{\eta}}\simeq 	\beta_*\beta^{*}(\AAA_{\inf, X_{\eta}}),
	\]
	and the $A_{\inf}$-module $R^i\beta_*\beta^{*}(\AAA_{\inf,X_{\eta}})$ is almost zero for all $i>0$ with respect to the ideal $W(\gothm_{\cO^{\flat}_{C}})\subset A_{\inf}$. As $\mu\in W(\gothm_{\cO^{\flat}_C})$, it follows that 
	\[\AAA_{\inf, X_{\eta}}[1/\mu]\xra{\sim } R\beta_*\beta^{*} (\AAA_{\inf, X_{\eta}}[1/\mu])
	\]
	is an isomorphism.
\end{proof}




\if false
\begin{lemma}\label{L:exact-sequence-Ainf}
	One has an exact sequence of abelian sheaves on $X_{\eta, \qproet}^{\diamond}$:
	\[
	0\to \ZZ_p\to \AAA_{\inf, X_{\eta}}[1/\mu]\xra{\phi-1}\AAA_{\inf, X_{\eta}}[1/\mu]\to0.
	\]
\end{lemma}
\begin{proof}
	It is well known that one has an Artin--Schreier exact sequence 
	\[
	0\to \ZZ_p\to \AAA_{\inf, X_{\eta}}\xra{\phi-1}\AAA_{\inf, X_{\eta}}\to 0.
	\]
	Now let $\frac{x}{\mu^n}\in (\AAA_{\inf, X_{\eta}}[1/\mu])^{\phi=1}$ for some zero local section $x$ of $\AAA_{\inf ,X_{\eta}}$ and integer $n\geqslant0$. Then we have
	$\phi(x)=(\frac{\phi(\mu)}{\mu})^nx$. As $[p]_q=\frac{\phi(\mu)}{\mu}$ is topologically nilpotent in $A_{\inf}$, we must have $n=0$ and $x\in \AAA_{\inf, X_{\eta}}^{\phi=1}=\ZZ_p$.

\end{proof}

\fi

We have the  following comparison theorem for semi-stable formal schemes over $\cO_C$.

\begin{theorem}\label{T:etale-comparison-O_C}
	Let  $X$ be a semistable $p$-adic formal   scheme over $\Spf(\cO_C)$ as in Definition~\ref{D:semi-stable} with associated log formal scheme $\uX$, and $(\calE, \phi_{\calE})$ be an object of $\CRhat^{\anfr}(\underline X_{\prism}, \cO_{\prism})^{\phi=1}$. Then:
	\begin{enumerate}
		\item[(1)] There exists a canonical isomorphism   
		\[
		Ru_*(\calE[1/\mu])\xra{\sim} R\nu_*(T(\calE)\otimes_{\ZZ_p}\AAA_{\inf,X_{\eta}}[1/\mu]).
		\]
		Moreover, this isomorphism is	equivariant under the Frobenius action in the sense that the following diagram is commutative
		\[
		\xymatrix{
			Ru_*(\calE[\frac{1}{\mu}])\ar[rr]^{\eqref{E:crystal-Frob-mu}}\ar[d]^{\simeq} && 	Ru_*(\calE[\frac{1}{\phi(\mu)}]) \ar[d]^{\simeq}\\
			R\nu_*(T(\calE)\otimes_{\ZZ_p}\AAA_{\inf,X_{\eta}}[\frac{1}{\mu}])\ar[rr]^-{1\otimes \phi_{\AAA_{\inf, X_{\eta}}}}	&&R\nu_*(T(\calE)\otimes_{\ZZ_p}\AAA_{\inf,X_{\eta}}[\frac{1}{\phi(\mu)}])
		}	\]
		where  the lower horizontal one is induced by the endomorphism $1\otimes \phi_{\AAA_{\inf, X_{\eta}}}$ of the sheaf $T(\calE)\otimes_{\ZZ_p}\AAA_{\inf,X_{\eta}}$.
		
		\item[(2)] Assume moreover that $X$ is  proper over $\Spf(\cO_C)$. Then one has a canonical isomorphism 
		\[
		R\Gamma(\uX_{\prism}, \calE[1/\mu])\simeq R\Gamma(X_{\eta,\et}, T(\calE))\otimes_{\ZZ_p}A_{\inf}[1/\mu],
		\] 
		equivariant under the Frobenius action in a similar sense as in  (1). In particular, one has a canonical isomorphism 
		\[
		H^i(\uX_{\prism}, \calE[1/\mu])\simeq H^i(X_{\eta,\et}, T(\calE))\otimes_{\ZZ_p}A_{\inf}[1/\mu],
		\]
	for all integers $i\in \ZZ$.	
	\end{enumerate}
\end{theorem}

\begin{proof}
	
	(1) Putting together all the previous results, one has a sequence of canonical isomoprhisms 
	\begin{align}
		Ru_*(\calE[1/\mu])&\xra{\sim}Ru^{\perf}_*(\calE^{\perf}[1/\mu])\nonumber \\
		&\cong  Ru_*R\gamma_*\gamma^{*}(\calE[1/\mu]) \nonumber \\ 
		&\cong R\nu_*R\beta_*	\beta^{*}\big(T(\calE)\otimes_{\ZZ_p}\AAA_{\inf,X_{\eta}} [1/\mu]\big)\nonumber \\
		&\cong R\nu_*(T(\calE)\otimes_{\ZZ_p}\AAA_{\inf,X_{\eta}}[1/\mu])\nonumber
	\end{align}
	where the first isomorphism is Proposition~\ref{P:coh-prismatic-perfection}, the second follows from  Lemma~\ref{L:almost-iso-v-to-qsyn}, and the last two isomorphism are established in  Lemma~\ref{L:construction of morphism}.
	Since all the isomorphisms above are functorial in $\calE$, 
	the equivariance under the Frobenius action follows by functoriality from Lemma~\ref{L:construction of morphism}(1).

	(2)	Taking the cohomology over $X_{\et}$, one deduces from (1) that 
	\begin{align*}
		R\Gamma(\uX_{\prism},\calE[1/\mu] )&\simeq 	R\Gamma(X_{\et}, Ru_{*}(\calE[1/\mu]))\\
		&\simeq  R\Gamma(X_{\et}, R\nu_*(T(\calE)\otimes_{\ZZ_p}\AAA_{\inf,X_{\eta}}[1/\mu]))\\
		&\simeq R\Gamma(X_{\eta, \proet}, T(\calE)\otimes_{\ZZ_p}\AAA_{\inf,X_{\eta}}[1/\mu]).
	\end{align*}
	By  a limit version of  Scholze's primitive comparison theorem \cite[Theorem~5.1]{Sch1},  for each integer $i\geqslant0$, one gets an almost isomorphism of  $A_{\inf}$-modules:
	\[
	H^i(X_{\eta, \et}, T(\calE))\otimes A_{\inf}\simeq^a H^i(X_{\eta,\proet}, T(\calE)\otimes_{\ZZ_p}\AAA_{\inf, X_{\eta}}).
	\]
	After inverting $\mu$, one finds that the  canonical map  
	\[
	R\Gamma(X_{\eta, \et}, T(\calE))\otimes A_{\inf}[1/\mu]\xra{\sim} R\Gamma(X_{\eta, \proet}, T(\calE)\otimes_{\ZZ_p}\AAA_{\inf, X_{\eta}}[1/\mu]).
	\]
	is an isomorphism in $D(A_{\inf})$. The statement on the cohomology groups $H^i$ follows from the flatness of  $A_{\inf}$ over $\ZZ_p$. 
\end{proof}

\section{Breuil--Kisin cohomology of analytic log prismatic crystals}\label{S:analytic-cyrstals}

In this section, we apply our theory over $\cO_C$ in previous sections to semi-stable formal scheme over a discrete valuation ring. 

Let $K$ denote a complete discrete valuation field over $\QQ_p$ with perfect residue field $k$,  $\cO_K$ be its ring of integers. 
Fix a uniformizer $\pi\in K$, and let $E(u)\in W(k)[u]$ be the Eisenstein polynomial of $\pi$ over $W(k)[1/p]$. We put $\gothS=W(k)[[u]]$, and equip it with the lift of Frobenius endomorphism $\phi$ such that $\phi(u)=u^p$.  





\subsection{Breuil--Kisin log prism}\label{S:BK-log-prism} 
Let  $X=\Spf(R)$ be   affine  small semi-stable formal scheme with a framing  as in \eqref{E:small-framing}:
\[
\square \colon R^{\square}:=\cO_K\langle T_0,\cdots, T_r, T_r^{\pm 1},\cdots, T_d^{\pm1}\rangle/(T_0\cdots T_r-\pi)\to R.
\]
Let $\uX$ be the associated log formal scheme with caononical log structure. 
 Recall that  $\uX$ is induced by the prelog structure  on $R$
\[
\alpha_{R}\colon M_R:=\NN^{r+1} \cdot e_i\to R^{\square}\to R
\]
given by  $e_i\mapsto T_i$  (cf. \S~\ref{S:local-chart}).
Consider the $(p,u)$-adically complete topological ring 
\[
R^{\square}_{\gothS}:=\gothS\langle T_0,\cdots, T_r, T_{r+1}^{\pm 1}, \cdots, T_d^{\pm 1} \rangle/(T_0\cdots T_r-u),
\]
and the prelog structure 
\[
\alpha_{R^{\square}_{\gothS}}\colon  \NN^{r+1}\to R^{\square}_{\gothS}
\]
given by $e_i\mapsto T_i$ for $0\leqslant i\leqslant r$.
We set $\delta(u)=\delta(T_i)=0$ for $0\leqslant i\leqslant d$ and $\delta_{\log}=0$ so that $(R^{\square}_{\gothS},\delta, \NN^{r+1},  \delta_{\log})$ becomes a $\delta_{\log}$-ring. 
It is clear that $R^{\square}_{\gothS}/(E(u))\cong R^{\square}$. By deformation theory, there exists a unique $(p,u)$-completely \'etale map 
$
R^{\square}_{\gothS}\to R_{\gothS}
$ 
whose reduction modulo $E(u)$ gives  the framing  $R^{\square}\to R$. 

By \cite[Lemma 2.13]{Kos}, the $\delta_{\log}$-structure on $(R^{\square}_{\gothS}, \NN^{r+1})$ extends uniquely to the prelog ring $(R_{\gothS}, \NN^{r+1}\to R^{\square}_{\gothS}\to R_{\gothS} )$ so that $(R_{\gothS}, (E(u)), \NN^{r+1})=(R_{\gothS}, \delta, (E(u)), \NN^{r+1}, \delta_{\log})$ is a bounded prelog prism. Together with the isomorphism of log schemes
\[
\uX\simeq (\Spf(R), \NN^{r+1}\to R_{\gothS}\to R=R_{\gothS}/(E(u)))^a,
\]
we get a bounded log prism $\uR_{\gothS}\colon=(R_{\gothS}, (E(u)), \NN^{r+1})^a$  on $\uX_{\prism}$, which we call \emph{Breuil--Kisin log prism}. We recall some basic properties of Breuil--Kisin log prism in \cite{DLMS}.

\begin{lemma}[Du--Liu--Moon--Shimizu]\label{L:BK-log-prism}
	The object $\uR_{\gothS}$ is a strong covering object in $\uX_{\prism}$ (Definition~\ref{D:versal object}), i.e. for every object  $(B, I_B, M_{\Spf(B)})$ in $\uX_{\prism}$, the coproduct   $$(B_{\gothS}, I_{B_{\gothS}}, M_{\Spf(B_{\gothS})})\colon =(B, I_B, M_B)\coprod \uR_{\gothS}$$  
	in $\uX_{\prism}^{\opp}$ exists, and   the canonical map $$(B, I_B, M_{\Spf(B)})\to (B_{\gothS}, I_{B_{\gothS}}, M_{\Spf(B_{\gothS})})$$
	is a  flat cover. 

\end{lemma}

\begin{proof}
	This is exactly \cite[Lemma~2.8, 2.9]{DLMS}.
\end{proof}

\subsection{Stratified modules}	Consider now the \v{C}ech nerve $\uR_{\gothS}^{\bullet}\colon =(R_{\gothS}^\bullet, (E(u)), M_{R_{\gothS}}^{\bullet})$ of $\uR_{\gothS}$ over the final object of $\Shv(\uX_{\prism})$, i.e.  each $\uR_{\gothS}^n:=(R_{\gothS}^n, (E(u)), M_{R_{\gothS}}^{n})$ with $n\geqslant0$ the $(n+1)$-fold self product of $(R_{\gothS}, (E(u)), \NN^{r+1})^a$ in $\uX_{\prism}$. 
Let $\Strat(R_{\gothS}^{\bullet})$ be the category of (complete) statified modules over $\uR_{\gothS}^{\bullet}$ (Definition~\ref{D:stratified-modules}), and let $\Strat^{\anfr}(R^{\bullet}_{\gothS})$ be the full subcategory of $\Strat(R^\bullet_{\gothS})$ consisting of obejcts $(E,\varepsilon)$ such that the restriction to $\Spec(R_{\gothS})\backslash V(p,u)$ of the quasi-coherent sheaf on $\Spec(R_{\gothS})$ attached to $E$ is locally free of finite rank.

If $(M, \epsilon)$ is an object of $\Strat(R_{\gothS}^{\bullet})$, then the Frobenius twist $\phi_{R_{\gothS}}^*M$ is equipped with a natural stratification 
\[
\phi_{R^1_{\gothS}}^*(\epsilon)\colon p_1^*(\phi_{R_{\gothS}}^*E)=\phi_{R^1_{\gothS}}^*p_1^*(M)\xra{\sim} \phi_{R^1_{\gothS}}^*p_0^*(M)=p_0^*(\phi_{R_{\gothS}}^*M). 
\]
Here, the maps $p_0,p_1\colon R_{\gothS}\to R^{1}_{\gothS}$ are the maps corresponding to the two face maps of simplices $[0]\to [1]$ (See \S~\ref{S:simplicial notation}). 
We call $(\phi_{R_{\gothS}}^*M, \phi_{R^1_{\gothS}}^*(\varepsilon))$ the Frobenius twist of $(M, \epsilon)$. 

\begin{definition}\label{D:Frobenius-strat}  A Frobenius structure on  an object  $(M,\epsilon)$  of $\Strat^{\anfr}(R^{\bullet}_{\gothS})$ is  an isormorphism 
	\[
	\phi_{M}\colon \phi^*_{R_{\gothS}}(M)[\frac{1}{E(u)}]\xra{\sim} M[\frac{1}{E(u)}].
	\]
	compatible with the stratification, i.e.
	the diagram 
	\[
	\xymatrix{
		p_1^*\phi_{R^1_{\gothS}}^*(M)[\frac{1}{E(u)}]\ar[rr]^{p_1^*(\phi_M)}_{\simeq}\ar[d]_{\phi_{R^1_{\gothS}}^*(\epsilon)}^{\simeq} && p_1^*M[\frac{1}{E(u)}]\ar[d]^{\epsilon}_{\simeq}\\
		p_0^*\phi_{R^1_{\gothS}}^*(M)[\frac{1}{E(u)}]\ar[rr]^{p_0^*(\phi_M)}_{\simeq} && p_0^*(M)[\frac{1}{E(u)}]
	}
	\]
	is commutative.
	A morphism of stratified modules with Frobenius structure  is a morphism  in $\Strat^{\anfr}(R_{\gothS}^{\bullet})$ compatible with the Frobenius structures on the source and the target.
	
	Let $\Strat^{\anfr}(R^{\bullet}_{\gothS})^{\phi=1}$ denote the category of stratified modules over $R_{\gothS}$ with  Frobenius structure.  
\end{definition}

As  $\uR_{\gothS}$ is a strong covering object in $\uX_{\prism}$ (Lemma~\ref{L:BK-log-prism}), Proposition~\ref{P:crystals-strat} implies that the evaluation at $\uR^{\bullet}_{\gothS}$ induces an equivalence of categories 
\begin{equation*}
	\CRhat(\uX_{\prism}, \cO_{\prism})\simeq \Strat(R_{\gothS}^{\bullet}):\quad \calE\mapsto (\calE(\uR_{\gothS}), \epsilon_{\calE}).
\end{equation*}
It is clear that if an object $\calE$ of $\CRhat(\uX_{\prism}, \cO_{\prism})$ lies in the subcategory $\CRhat^{\anfr}(\uX_{\prism}, \cO_{\prism})$, then $(\calE(\uR_{\gothS}), \epsilon_{\calE})$ lies in $\Strat^{\anfr}(R_{\gothS}^{\bullet})$.
 If $(\calE,\phi_{\calE})$ is an object of $\CRhat(\uX_{\prism}, \cO_{\prism})^{\phi=1}$, then $\phi_{\calE}$  induces an Frobenius structure $\phi_{\calE(\uR_{\gothS})}$ on $(\calE(\uR_{\gothS}),\epsilon_{\calE})$ so that $(\calE(\uR_{\gothS}),\epsilon_{\calE}, \phi_{\calE(\uR_{\gothS})})$ becomes an object of $\Strat^{\anfr}(R_{\gothS}^{\bullet})^{\phi=1}$. Hence, one has  equivalences of categories 
\begin{equation}\label{E:equiv-crystal-strat-BK}
	\CRhat^{\anfr}(\uX_{\prism}, \cO_{\prism})^{(\phi=1)}\simeq \Strat^{\anfr}(R_{\gothS}^{\bullet})^{(\phi=1)}
\end{equation}

\begin{lemma}\label{L:completeness-base-change}
	Let $(A,I)$ be a prism with $A$ noeotherian, $M$ be a finitely generated $A$-module. Then for any  $(p,I)$-completely flat map of bounded prisms  $(A,I)\to (B,IB)$,  $M\otimes_{A}B$ is $(p,I)$-adically complete.
\end{lemma}
\begin{proof}
	
	It is clear that $M\otimes_AB$ is a  finite presented $B$-module. Since   derived $(p,I)$-complete modules are stable under quotients,  it follows  that $M\otimes_AB$ is derived $(p,I)$-complete.  By \cite[Lemma 3.6]{BS}, $I^p$ is principal. Write $I^p=(\xi)$ for some nonzero divisor $\xi\in A$.  
	By the explicit computation of the derived completion  \cite[\href{https://stacks.math.columbia.edu/tag/091Z}{Tag 091Z}]{stacks-project},  one has an isomorphism 
	\[
	M\otimes_AB \simeq R\lim_n \Kos^{\bullet}(p^n, \xi^n; M\otimes_AB),
	\]
	where $\Kos^{\bullet}(p^n,\xi^n; M\otimes_AB)$ is the Koszul complex for $M\otimes_AB$ and the sequence $(p^n, \xi^n)$. Note also that $B$ is flat over $A$ by  \cite[\href{https://stacks.math.columbia.edu/tag/0912}{Tag 0912}]{stacks-project}. Hence we have $\Kos^{\bullet}(p^n,\xi^n; M\otimes_AB)= \Kos^{\bullet}(p^n,\xi^n; M)\otimes_AB$. But the same argument as \cite[\href{https://stacks.math.columbia.edu/tag/0921}{Tag 0921}]{stacks-project} shows that the pro-object $\{\Kos^{\bullet}(p^n,\xi^n; M):n\in\NN\}$ is equivalent to $\{M/(p^n,\xi^n):n\in \NN\}$ in the derived category $D(A)$. It follows that 
	\begin{align*}
		M\otimes_AB &\simeq R\lim_n \big(\Kos^{\bullet}(p^n, \xi^n; M)\otimes_AB\big)\\
		&\simeq R\lim_n \big(M/(p^n,\xi^n)\otimes_AB\big)\\
		&\simeq \lim_n (M\otimes_AB)/(p^n,\xi^n)(M\otimes_AB).
	\end{align*}
\end{proof}

\begin{remark}
	Let $(E, \epsilon)$ be an object of $\Strat(R_{\gothS}^{\bullet})$. By Definition~\ref{D:stratified-modules}, the stratification $\epsilon$ is an isomorphism of $R^1_{\gothS}$-modules   $p_1^*(E)\xra{\sim}p_0^*(E)$, where $p_i^*(E)$  for $i=1,2$ is the $(p, E(u))$-adically complete  tensor product $E\widehat\otimes_{R_{\gothS}, p_i}R^1_{\gothS}$. However, if $M$ is finitely generated over $R_{\gothS}$, then Lemma~\ref{L:completeness-base-change} shows that $p_i^*(M)$ coincides with the usual tensor product $M\otimes_{ R_{\gothS}, p_i}R^1_{\gothS}$. 
\end{remark}

\subsection{Analytic prismatic $F$-crystals and semi-stable $\ZZ_p$-local systems}
Let $X$ be a semi-stable $p$-adic formal scheme over $\Spf(\cO_K)$ (Def.~\ref{D:semi-stable}).
We recall the notion of analytic prismatic $F$-crystal on $\uX_{\prism}$ introduced in \cite[Def. 3.1]{GR} and  \cite[Def. 3.3]{DLMS}.
For every object $(A, I_A, M_{\Spf(A)})$ of $\uX_{\prism}$, let $\Vect^{\an}(A, I_A)$ be the  category of vector bundles over $\Spec(A)\backslash V(p,I_A)$, and $\Vect^{\an}(A, I_A)^{\phi=1}$ be the category of pairs $(\calE_A, \phi_{\calE_A})$ where $\calE_A$ is an object of $\Vect^{\an}(A, I_A)$ and 
\begin{equation*}\label{E:Frob-struct-analytic}
	\phi_{\calE_A}\colon \phi_A^*(\calE_A)[\frac{1}{I_A}]\xra{\sim} \calE[\frac{1}{I_A}].
\end{equation*}
Then the category of analytic prismatic  ($F$-)crystals on $\uX_{\prism}$ is defined as 
\[
\Vect^{\an}(\uX_{\prism}, \cO_{\prism})^{(\phi=1)}:=\lim_{(A,I_A, M_{\Spf(A)})\in \uX_{\prism}}\Vect^{\an}(A,I_A)^{(\phi=1)}.
\]
For an object $\calE$  (resp. $(\calE, \phi_{\calE})$)  of $\Vect^{\an}(\uX_{\prism}, \cO_{\prism})$ (resp. of $\Vect^{\an}(\uX_{\prism}, \cO_{\prism})^{\phi=1}$), we  write $\calE_A$ (resp. $(\calE_A, \phi_{\calE_A})$) its image in $\Vect^{\an}(A,I_A)$  (resp. in $\Vect^{\an}(A,I_A)^{\phi=1}$) for any object $(A,I_A, M_{\Spf(A)})$ of $\uX_{\prism}$.

The restriction to the analytic locus defines a funtor 
\begin{equation}\label{E:restriction-analytic-crystal}
	j^*\colon \CRhat^{\anfr}(\uX_{\prism}, \cO_{\prism})^{(\phi=1)}\to \Vect^{\an }(\uX_{\prism}, \cO_{\prism})^{(\phi=1)}, 
\end{equation}

\begin{proposition}\label{P:can-extension-analytic}
	Let $X$ be a separated semi-stable $p$-adic formal scheme over $\Spf(\cO_K)$. Then the   restriction functor $j^* $
	admits a right adjoint $j_*$ such that the adjunction transformation of functors  $\id\xra{\sim} j^*j_*$ is an equivalence. Moreover, such a funtor $j_*$ is unique up to canonical equivalence.
	
\end{proposition}
\begin{proof}
Our  proof is motivated by \cite[Thm. 5.10]{GR} and \cite[Lemma 3.8]{DLMS}.
The uniqueness of $j_*$ is an easy formal consequence of the adjunction property and the requirement that $\id\xra{\sim} j^*j_*$. 
	
	We  only need to prove the existence  for $j_*\colon  \Vect^{\an }(\uX_{\prism}, \cO_{\prism})\to \CRhat^{\anfr}(X_{\prism}, \cO_{\prism})$, because  a Frobenius structure on an object of $\CRhat^{\anfr}(\uX_{\prism}, \cO_{\prism}, )$ is given over the analytic locus. 	
	Consider first the case when $X=\Spf(R)$ is  small affine semi-stable. 
	We fix a framing map $\square \colon R^{\square}\to R$ as in \eqref{E:small-framing}.
	Consider the Breuil--Kisin log prism  $\uR_{\gothS}$ and the \v{C}ech nerve $\uR^{\bullet}_{\gothS}$ as  in \S~\ref{S:BK-log-prism}.
	Similarly to Definition~\ref{D:stratified-modules}, let $\Strat(R_{\gothS}^{\bullet}\backslash V(p, E(u)))$ denote the category of pairs $(\calF^{\an}_{R_{\gothS}}, \varepsilon_{\calF^{\an}})$ where 
	\begin{itemize}
		\item $\calF^{\an}_{R_\gothS}$ is an object of  $\Vect^{\an}(R_{\gothS}, (E(u)))$, 
		\item $\epsilon_{\calF^{\an}}\colon p_1^*(\calF^{\an}_{R_{\gothS}})\xra{\sim} p_0^*(\calF^{\an}_{R_{\gothS}})$ is an isomorphism of in $\Vect^{\an}(R^1_{\gothS}, (E(u)))$ such that $p_{0,2}^*(\psi)\simeq p_{0,1}^*(\psi)\circ p_{1,2}^*(\psi)$.
	\end{itemize}
	Then similarly to Proposition~\ref{P:crystals-strat},  \cite[Thm. 7.8]{Mat}  implies that  the evaluation at $\uR_{\gothS}$ induces a natural equivalence of categories  $$\Vect^{\an}(\uX_{\prism}, \cO_{\prism})\simeq \Strat(R_{\gothS}^{\bullet}\backslash V(p, E(u))).$$
	Then $j^*$ corresponds to  the natural restriction functor 
	\[
	j^*_{R_{\gothS}}\colon \Strat^{\anfr}(R^{\bullet}_{\gothS})\to \Strat(R^{\bullet}_{\gothS}\backslash V(p, E(u))).
	\]
	It suffices to construct a functor $j_{R_{\gothS},*}\colon \Strat(R^{\bullet}_{\gothS}\backslash V(p, E(u)))\to \Strat^{\anfr}(R^{\bullet}_{\gothS})$ right adjoint to $j^*_{R_{\gothS}}$ such that the adjunction $\id\xra{\sim} j^*_{R_{\gothS}}j_{R_{\gothS},*}$ is an equivalence.
	
	Let $(\calF^{\an}_{R_{\gothS}}, \epsilon_{\calF^{\an}})$ be an object of $\Strat(R_{\gothS}^{\bullet}\backslash V(p, E(u)))$. 
	For any integer $n\geqslant 0$, denote  by $$j^n\colon \Spec(R^n_{\gothS})\backslash V(p,E(u))\hra \Spec(R^n_{\gothS})$$  the natural inclusion. 
	Note that $R_{\gothS}$ is regular noetherian and $V(p,E(u))$ has codimension $2$ in $\Spec(R_{\gothS})$. It  follows from \cite[\href{https://stacks.math.columbia.edu/tag/0BK3}{Tag 0BK3}]{stacks-project} that  $j^0_*\calF^{\an}_{R_{\gothS}}$  is a coherent sheaf on $\Spec(R_{\gothS})$. Hence,  
	\[
	M:=\Gamma(\Spec(R_{\gothS}), j^0_*\calF^{\an}_{R_{\gothS}})=\Gamma(\Spec(R_{\gothS})\backslash V(p, E(u)), \calF^{\an}_{R_{\gothS}})
	\]
	is a finitely generated $R_{\gothS}$-module, thus automatically $(p,E(u))$-adically complete. As $p_i\colon R_{\gothS}\to R^1_{\gothS}$ for $i=0,1$ is flat by Lemma~\ref{L:BK-log-prism} and  \cite[\href{https://stacks.math.columbia.edu/tag/0912}{Tag 0912}]{stacks-project}, the flat base change formula for quasi-coherent sheaves \cite[\href{https://stacks.math.columbia.edu/tag/0GN5}{Tag 0GN5}]{stacks-project} gives  a canonical isomorphism 
	\[
	c_{M}(p_i)\colon M\otimes_{R_{\gothS}, p_i} R^1_{\gothS}\xra{\sim} \Gamma(\Spec(R^1_{\gothS}), j^1_*(p_1^*(\calF^{\an}_{R_{\gothS}}))).
	\]
	Note that  $M\otimes_{R_{\gothS}, p_i} R^1_{\gothS}$ is already $(p, E(u))$-adically complete by Lemma~\ref{L:completeness-base-change}, hence it coincides with the $(p,E(u))$-adically complete base change $p_i^*(M)$ in the sense of \S~\ref{S:simplicial notation}.
	Put  $\epsilon_M:= c_M(p_0)^{-1}\circ \Gamma(\Spec(R^1_{\gothS}), j^1_*\epsilon_{\calF^{\an}})\circ c_M(p_1)$.
	If $p_i^2\colon \uR_{\gothS}\to \uR^2_{\gothS}$ for $i\in [2]$ denotes the natural degeneracy map corresponding to $[0]\to [2]$ given by  $0\mapsto i$, then each $p_i^2\colon R_{\gothS}\to R^2_{\gothS}$ is  $(p,E(u))$-adically faithfully flat, and the flat base change  implies that the  canonical   map
	\[
	p^{2,*}_{i}(M)\xra{\sim} \Gamma(\Spec(R^2_{\gothS}), j^2_*(p^{2,*}_i(\calF^{\an}_{\uR_{\gothS}})) )
	\]
	is an isomorphism. From this, we deduce easily that 
	$p_{0,2}^*(\epsilon_M)=p_{0,1}^*(\epsilon_M)\circ p_{1,2}^*(\epsilon_M)$ so that  $j_{R_{\gothS},*}(\calF^{\an}_{R_{\gothS}}, \epsilon_{R_{\gothS}})\colon =(M, \epsilon_M)$ is an object of $\Strat^{\anfr}(R_{\gothS}^{\bullet})$.
	We get thus a functor 
	\[
	j_{R_{\gothS},*}\colon \Strat(R^{\bullet}_{\gothS}\backslash V(p, E(u)))\to \Strat^{\anfr}(R^{\bullet}_{\gothS}).
	\]
	By the usual adjoint property for the pairs $(j_{n}^*, j_{n,*})$ with $n\geqslant 0$ for coherent sheaves, we see easily that $j_{R_{\gothS},*}$ is  right adjoint to  $j^*_{R_{\gothS}}$,  and the adjunction transformation $\id \xra{\sim}j^*_{R_{\gothS}}j_{R_{\gothS}, *}$ is an equivalence. This finishes the proof of the Proposition when $X$ is small affine semi-stable.

	We consider now the general case. Let $(X_{\alpha})_{\alpha\in A}$ be a $p$-completely \'etale cover of  $X$ by small affine semistable formal schemes $X_{\alpha}=\Spf(R_{\alpha})$. We fix a framing $\square\colon R_{\alpha}^{\square}\to R_{\alpha}$ for each $\alpha
	\in A$. The discussion above gives a  functor $j_{\alpha,*}\colon \Vect^{\an}(\uX_{\alpha, \prism}, \cO_{\prism})\to \CRhat^{\anfr}(\uX_{\alpha, \prism}, \cO_{\prism})$ right adjoint to the restriction functor. To finish the proof of the Proposition, it suffices to show that $j_{\alpha,*}$ glue together to  a functor $j_*\colon \Vect^{\an}(\uX_{ \prism}, \cO_{\prism})\to \CRhat^{\anfr}(\uX_{ \prism}, \cO_{\prism})$. 
	
	Let $\calF^{\an}$ be an object of $\Vect^{\an}(\uX_{\prism}, \cO_{\prism})$, and $\calF^{\an}_{\alpha}$ be its restriction to $\uX_{\alpha, \prism}$. We need to show that, for any $\alpha,\beta\in A$, if $X_{\alpha, \beta}=X_{\alpha}\times_{X}X_{\beta}$, then  $j_{\alpha,*}(\calF^{\an}_{\alpha})|_{\uX_{\alpha, \beta, \prism}}=j_{\beta,*}(\calF^{\an}_{\beta})|_{\uX_{\alpha, \beta, \prism}}$
	As $X$ is separated, $X_{\alpha, \beta}=\Spf(R)$ is affine small semi-stable. It is equipped with two framings $\square_{?}\colon R_{?}^{\square}\to R_{?}\to R$ for $?\in \{\alpha, \beta\}$.
	Let $\uR_{\gothS}^{\alpha}$ and $\uR_{\gothS}^{\beta}$ be the  Breuil--Kisin log prisms in $\uX_{\alpha, \beta}$ corresponding to $\square_{\alpha},\square_{\beta}$  respectively, and $\uR^{\alpha,\bullet}_{\gothS}$ and $\uR^{\beta,\bullet}_{\gothS}$ be the corresponding \v{C}ech nerves.  Let
	$$\uR_{\gothS}^{\alpha, \beta}=\uR_{\gothS}^ {\alpha}\coprod \uR_{\gothS}^{ \beta}$$ be  the  coproduct in $\uX_{\alpha, \beta, \prism}^{\opp}$. Then $\uR^{\alpha, \beta}_{\gothS}$ is also a cover of the final object of $\uX_{\opp}$. 
	Let  $\uR^{\alpha, \beta, \bullet}_{\gothS}$ be the \v{C}ech nerve of $\uR_{\gothS}^{\alpha, \beta}$ over the initial object of $\uX_{\alpha, \beta, \prism}^{\opp}$.
	We can define similarly the category $\Strat^{\anfr}(R^{\alpha, \beta}_{\gothS})$ in an evident sense, and the same argument as Proposition~\ref{P:crystals-strat} shows that  the evaluation at $\uR_{\gothS}^{\alpha, \beta, \bullet}$ have an equivalence of cateogries 
	\[
	\CRhat^{\anfr}(\uX_{\prism}, \cO_{\prism})\simeq \Strat^{\anfr}(R^{\alpha, \beta,\bullet}_{\gothS}).
	\] 
	We have a morphism of 
	cosimplicial objects in $\uX_{\alpha, \beta, \prism}^{\opp}$: $
	q_{?}^{\bullet}\colon \uR_{\gothS}^{?, \bullet}\to \uR_{\gothS}^{\alpha, \beta, \bullet}
	$
	which induces a natural pull-back functor 
	\[
	q_?^*\colon \Strat^{\anfr}(R^{?,\bullet}_{\gothS})\to \Strat^{\anfr}(R^{\alpha,\beta,\bullet}_{\gothS}).
	\]
	Let $(M_{?}, \varepsilon_{M_{?}})=j_{R^?_{\gothS}, *}(\calF^{\an}_{R^{?}_{\gothS}}, \varepsilon_{\calF_{?}^{\an}})$ for $?\in \{\alpha, \beta\}$ be the object of $\Strat^{\anfr}(R^{?,\bullet}_{\gothS})$ corresponding to $j_{\alpha,*}(\calF^{\an}_{\alpha})$. 
	To finish the proof of the Proposition, it is enough to show that $q_{\alpha}^*(M_{\alpha}, \varepsilon_{M_{\alpha}})=q_{\beta}^*(M_{\beta}, \varepsilon_{M_{\beta}})$.
	
	By the construction  in the first step, we have 
	$M_{?}=\Gamma(\Spec(R^{?}_{\gothS})\backslash V(p,E(u)), \calF^{\an}_{R^{\alpha}_{\gothS}})$ for $?\in \{\alpha,\beta\}$.  Note that for each integer $n\geqslant 0$, the underlying map $q_{?}^n\colon R_{\gothS}^{?,n}\to R_{\gothS}^{\alpha, \beta, ?}$ is  $(p,E(u))$-completely flat by Lemma~\ref{L:BK-log-prism}.  
	Put 
	\[
	M_{\alpha, \beta}:=\Gamma(\Spec(R^{\alpha,\beta}_{\gothS})\backslash V(p,E(u)), \calF^{\an}_{R^{\alpha,\beta}_{\gothS}}).
	\]
	The flat base change for quasi-coherent sheaves implies that
	\[
	M_{\alpha}\otimes_{R_{\gothS}^{\alpha}, q_{\alpha}^0}R^{\alpha, \beta}_{\gothS}\simeq M_{\alpha, \beta}\simeq M_{\beta}\otimes_{R_{\gothS}^{\alpha}, q_{\beta}^0}R^{\alpha, \beta}_{\gothS},
	\] 
	and  both $q^{1,*}_{\alpha}(\epsilon_{M_{\alpha}})$  and  $q^{1,*}_{\beta}(\epsilon_{M_\beta})$ are identified with the natural stratification on $M_{\alpha, \beta}$ given by the composed isomorphism 
	\[
	p^{\alpha, \beta, *}_1(M_{\alpha, \beta})\simeq \Gamma(\Spec(R^{\alpha, \beta, 1}_{\gothS})\backslash V(p, E(u)), \calF^{\an}_{R^{\alpha,\beta,1}_{\gothS}})\simeq p_0^{\alpha, \beta, *}(M_{\alpha, \beta})
	\]
	This shows that $q_{\alpha}^*(M_{\alpha}, \epsilon_{M_{\alpha}})=q_{\beta}^*(M_{\beta}, \epsilon_{M_{\beta}})$, and hence finishes the proof of the Proposition.
\end{proof}

\begin{definition}\label{D:can-ext-analytic}
	For an object $\calE$ of $\Vect^{\an}(\uX_{\prism}, \cO_{\prism})$, we call $j_*\calE$ given by Proposition~\ref{P:can-extension-analytic}  the \emph{canonical extension} of  $\calE$. 
\end{definition}

\begin{remark}\label{R:torsion-freeness}	
	If $X=\Spf(R)$ is affine small semi-stable over $\Spf(\cO_K)$ with framing \eqref{E:small-framing} and $\calE\in \Vect^{\an}(\uX_{\prism}, \cO_{\prism})$, then it follows from the construction of $j_*$ that
	\[
	(j_*\calE)(\uR_{\gothS})=\Gamma(\Spec(\uR_{\gothS})\backslash V(p,E(u)), \calE_{R_{\gothS}}).
	\]
	It follows that  $j_*\calE(\uR_{\gothS})$ is a finitely generated $\uR_{\gothS}$-module such that the following holds:
	\begin{itemize}
		\item Both $(j_*\calE)(\uR_{\gothS})[1/p]$ and $(j_*\calE)(\uR_{\gothS})[1/E(u)]$ are finite projective module over $\uR_{\gothS}[1/p]$ and $\uR_{\gothS}[1/E(u)]$ respectively.
		
		\item   $(j_*\calE)(\uR_{\gothS})=(j_*\calE)(\uR_{\gothS})[1/p]\cap (j_*\calE)(\uR_{\gothS})[1/E(u)]$.
	\end{itemize}
	In particular, $(j_*\calE)(\uR_{\gothS})$ is a torsion free $R_{\gothS}$-module. 
\end{remark}

\if false

\begin{lemma}
	Let $\CRhat^{\anfr}(\uX_{\prism}, \overline{\cO}_{\prism})$ be the full subcategory of $\CRhat^{\anfr}(\uX_{\prism}, {\cO}_{\prism})$ consisting of objects $\calE$ that are annihilated by $\calI_{\prism}$. Then for every object $\calE$ of $\Vect^{\an}(\uX_{\prism}, \cO_{\prism})$,  $j_*(\calE)/\calI_{\prism}j_*(\calE)$ is an object of $\CRhat^{\anfr}(\uX_{\prism}, \overline{\cO}_{\prism})$.
\end{lemma}

\begin{proof}
	By fpqc-descent, we see easily that
	\[
	(j_*(\calE)/\calI_{\prism}j_*(\calE))(B,J,M_B)=j_*(\calE)(B,J,M_B)/Jj_*(\calE)(B,J,M_B)
	\]
	for every object $(B,J,M_B)$ of $\uX_{\prism}$. To finish the proof of the Lemma, we need to show that $j_*(\calE)(B,J,M_B)/Jj_*(\calE)(B,J,M_B)$ is $p$-adically complete. Since derived $(p, J)$-complete modules are stable under quotients, we see that $j_*(\calE)(B,J,M_B)/Jj_*(\calE)(B,J,M_B)$ is derived $p$-complete. By \cite[\href{https://stacks.math.columbia.edu/tag/0BKG}{Tag 0BKG}]{stacks-project}, it is $p$-adically complete if and only if it has bounded $p^{\infty}$-torsion. Since this property can be checked after a $(p, J)$-completely faithfully flat base change in $B$, we may assume that $X=\Spf(R)$ is affine small with a framing \eqref{E:framing-semi-stable}, and there exists a morphism of log prisms $\uR_{\gothS}\to (B,J,M_B)$. By the crystal property for $j_*(\calE)$, we are reduced 
\end{proof}

\fi


\subsection{Relative log prismatic site}
 Applying the construction of \S~\ref{S:BK-log-prism} to $\Spf(\cO_K)$, we get the absolute Breuil--Kisin log prism $\underline{\gothS}=(\gothS, (E(u)),  \NN\xra{1\mapsto u} \gothS)^a$.
We have the relative log prismatic site $(\uX/\ufS)_{\prism}$ (cf. \S~\ref{S:relative-prismatic-site}). We have the natural  categories of prismatic crystals  $\CRhat((\uX/\ufS)_{\prism}, \cO_{\prism})$, $\CRhat^{\anfr}((\uX/\ufS)_{\prism}, \cO_{\prism})^{(\phi=1)}$ defined  similarly as Definition~\ref{D:crystals} and \ref{D:complete-prismatic-F-crystal}. 
There exist natural restriction functors 
\begin{align*}
	\CRhat(\uX_{\prism}, \cO_{\prism})&\to \CRhat((\uX/\ufS)_{\prism}, \cO_{\prism})\\
	\CRhat^{\anfr}(\uX_{\prism}, \cO_{\prism})^{(\phi=1)}&\to \CRhat^{\anfr}((\uX/\ufS)_{\prism}, \cO_{\prism})^{(\phi=1)}
\end{align*}
By abuse of notation, we will use the same notation to denote the image in $\CRhat((\uX/\ufS)_{\prism}, \cO_{\prism})$ of an object of $\CRhat(\uX_{\prism}, \cO_{\prism})$.

\subsection{Local description}
Assume now that  $X=\Spf(R)$ is affine small semi-stable of relative dimension $d$ with a framing $\square\colon R^{\square}\to R$ \eqref{E:small-framing}, and let $\uR_{\gothS}=(R_{\gothS}, (E(u)), \NN^{r+1})^a$ be the Breuil--Kisin log prism for the framing $\square$ (see \ref{S:BK-log-prism}).  
There exists a  morphism of log prisms $\ufS\to \uR_{\gothS}$ induced by the natural inclusion $\gothS\hra R_{\gothS}$ and the diagonal map of monoids $\NN\hra \NN^{r+1}$.  Hence,  the Breuil--Kisin log prism $\uR_{\gothS}$ becomes an object of $(\uX/\ufS)_{\prism}$.  
By \cite[Lemma~3.16]{MW}, $\uR_{\gothS}$ is a strong covering object in $(\uX/\ufS)_{\prism}$ in an evident similar sense as Definition~\ref{D:versal object}.
Let $\uR_{\gothS}^{\rel,\bullet}$ denote the \v{C}ech nerve of $\uR_{\gothS}$ over the final object of $\Shv((\uX/\ufS)_{\prism})$.
For an abelian sheaf $\calF$ on $(\uX/\ufS)_{\prism}$, we write $\CA(\uR^{\rel, \bullet}_{\gothS}, \calF)$ the simple complex attached to cosimplicial abelian group $\calF(\uR_{\gothS}^{\rel, \bullet})$. 

\begin{lemma}\label{L:CA-complex-relative}
	Let $\calF$ be an object of $\CRhat((\uX/\ufS)_{\prism}, \cO_{\prism})$, then  $R\Gamma((\uX/\ufS)_{\prism}, \calF)$ is computed by the complex $\CA(\uR^{\rel, \bullet}_{\gothS}, \calF)$. 
\end{lemma}

\begin{proof}
	This arguments are the same as Lemma~\ref{L:CA-complex-prismatic}, because $\uR_{\gothS}$ is a cover of the final object of $\Shv((\uX/\ufS)_{\prism})$ \cite[Lemma 3.16]{MW}. 
\end{proof}

\if false
We have the following description for the cosimplicial $\cO_K$-algebra $R^{\rel,\bullet}_{\gothS}/(E(u))$. 

\begin{lemma}[Min--Wang]\label{L:BK-relative}
	For every integer $n\geqslant 0$, there exists an  isomorphism
	\[
	R^{\rel, n}_{\gothS}/(E(u))\simeq R\langle Y_{j,i}:1\leqslant j\leqslant d, 1\leqslant i\leqslant n\rangle^{PD},
	\]
	where the right hand side is the $p$-adic completion of the free divided power polynomial ring over $R$ in varibles $Y_{j,i}$'s, and it is $R$ if $n=0$ by convention. 
	Moreover, if $\delta^n_{k}\colon R_{\gothS}^{\rel, n}\to R^{\rel,n+1}_{\gothS}$ with $0\leqslant  k\leqslant n+1$ denote the structure map corresponding to the order-preserving injection of simplices $[n]\to [n+1]$ that omits $k$, then 
	\[
	\delta^n_k(Y_{j,i})=\begin{cases}Y_{j+1, i}-Y_{1,i} &\text{if }k=0,\\
		Y_{j,i} &\text{if }1\leqslant j\leqslant k,\\
		Y_{j+1,i} &\text{if }0<k\leqslant j.
	\end{cases}
	\]
	
\end{lemma}

\begin{proof}
	This  is \cite[Lemma 3.19]{MW}.
\end{proof}

\fi

Let $\widehat\Omega^1_{R,\log}$ denote completed logarithmic differential module of $R$ relative to $\cO_K$. This is a free $R$-module with basis $\{\frac{dT_i}{T_i}:1\leqslant i\leqslant d\}$. 

\begin{definition}
	A complete log Higgs module over $R$ is a $p$-adically complete module $M$ together with an $R$-linear morphism 
	\[\Theta\colon M\to M{\otimes}_R \widehat{\Omega}_{R,\log}^1\]
	such that $\Theta\wedge\Theta=0$. If we write $\Theta=\sum_{i=1}^d \theta_i \otimes\frac{dT_i}{T_i}$, then $\Theta\wedge\Theta=0$ is equivalent to $\theta_i\theta_j=\theta_j\theta_i$ for all $i,j$. 
	
	For any $\underline{m}=(m_i)_i\in \NN^d$, put 
	\[
	\theta^{\underline m}=\prod_{i=1}^d\theta_i^{m_i}\in \End(M)
	\]
	We say that  $(M, \Theta)$ is topologically quasi-nilpotent if for every $x\in M$, $\theta^{\underline m}(x)$ tends to $0$ as $|\underline m|:=\sum_{i=1}^d m_i$ tends to infinity.
	
	Let $\Higgs^{\topnil}(R)$ denote the category of topologically quasi-nilpotent complete log Higgs module over $R$. 
\end{definition}
Note that a complete log Higgs module is a module with integrable connection over $(R,  R\xra{0} \widehat{\Omega}^1_{R,\log})$ defined in Definition~\ref{D:connection}. 
For a complete log Higgs module  $(M,\Theta)$, we denote by $(M\otimes_{R}\widehat\Omega^1_{R,\log}, \Theta^{\bullet})$ the associated de Rham complex. 

Consider the sheaf  $\overline{\cO}_{\prism}:=\overline{\cO}_{\prism}/E(u)$  on $(\uX/\ufS)_{\prism}$, and let $\CRhat(\uX_{\prism},\overline {\cO}_{\prism})$ denote the full subcategory of $\CRhat(\uX_{\prism}, {\cO}_{\prism})$ consisting of objects that are annihilated by $E(u)$. 

\begin{proposition}\label{P:Higgs-equiv-BK}
	Assume that $X=\Spf(R)$ is small affine semi-stable over $\Spf(\cO_K)$ with framing \eqref{E:small-framing}. Then 
	\begin{enumerate}
		\item For every object	$\calE$  of $\CRhat((\uX/\ufS)_{\prism},\overline{\cO}_{\prism})$, there exists a topologically quasi-nilpotent Higgs connection on $\calE(\uR_{\gothS})$
		\[
		\Theta_{\calE}\colon \calE(\uR_{\gothS})\to \calE(\uR_{\gothS})\otimes_R \widehat{\Omega}^1_{R,\log}
		\] 
		such that  the construction $\calE\mapsto (\calE(\uR_{\gothS}), \Theta_{\calE})$ induces an equivalence of categories
		\[
		\CRhat((\uX/\ufS)_{\prism}, \overline{\cO}_{\prism})\simeq \Higgs^{\topnil}(R).
		\]
		
		\item Via the equivalence in (1), we have an isomorphism in the derived category $D(\cO_K)$
		\[
		R\Gamma((\uX/\ufS)_{\prism}, \calE)\simeq (\calE(\uR_{\gothS})\otimes_{R}\widehat{\Omega}^{\bullet}_{R,\log}, \Theta_{\calE}^{\bullet})
		\]
		for every object $\calE$ of $\CRhat((\uX/\ufS)_{\prism}, \overline{\cO}_{\prism})$.
		
	\end{enumerate}
\end{proposition}

\begin{proof}
	The arguments  are  the same as  \cite[Thm. 3.21, 3.22]{MW}. In fact, the statements of \emph{loc. cit.} are only for locally free $\overline{\cO}_{\prism}$-crystals, but the essential ingredient is \cite[Lemma~3.19]{MW}, which is independent of this assumption.  
\end{proof}

Now let us return to the general separated semi-stable case, i.e.  we assume no longer that $X$ is affine small.

\begin{definition}\label{D:BK-cohomology}
	
	Let $\calE\in \Vect^{\an}(\uX_{\prism}, \cO_{\prism})$ be an analytic prismatic crystal on $\uX_{\prism}$, and $j_*\calE$ be its canonical extension. We put 
	\[
	R\Gamma_{\gothS}(X, \calE):=R\Gamma((\uX/\ufS)_{\prism}, j_*\calE),
	\]  
	and called it \emph{Breuil--Kisin cohomology of $\calE$}. For all $i\in \ZZ$, we write   $H^i_{\gothS}(X,\calE)$ for  the $i$-th cohomology of $R\Gamma_{\gothS}(X, \calE)$.
	
\end{definition}
For an object $(\calE, \phi_{\calE})$ of $\Vect^{\an}(\uX_{\prism}, \cO_{\prism})^{\phi=1}$, 
$\phi_{\calE}$ indcues   a natural  Frobenius endomorphism $R\Gamma_{\gothS}(X, \calE)$
\begin{equation}\label{E:Frob-end-BK-cohomology}
	\bigg(\phi_{\gothS}^*R\Gamma_{\gothS}(X, \calE)\bigg)[\frac{1}{E(u)}]\to R\Gamma_{\gothS}(X, \phi_{\cO_{\prism}}^*\calE)[\frac{1}{E(u)}]\xra{\phi_{\calE}} R\Gamma_{\gothS}(X, \calE)[\frac{1}{E(u)}].
\end{equation}

In general, for a noetherian ring $A$ and integers $a,b\in \ZZ$ with $a\leqslant b$, let $D^{[a,b]}_{\coh}(A)$ denote the full subcategory of the derived category $D(A)$ consisting of the objects $L$ such that each $H^i(L)$ is a finitely generated $A$-module and  $H^i(L)=0$ if $i\notin [a,b]$; equivalently,  an object $L\in D(A)$ lies in $D^{[a,b]}_{\coh}(A)$ if and only if it is quasi-isomorphic to a complex $E^{\bullet}$ concentrated in degrees $[a,b]$ with each term $E^i$ a finitely generated $A$-module.

\begin{lemma}\label{L:nakayama}
	Let $A$ be a noetherian ring that is $f$-adically complete for a nonzero divisor $f\in A$. Let $L\in D(A)$ be a derived $f$-complete object such that $L\otimes_A^LA/(f)$ lies in $D^{[a,b]}_{\coh}(A/(f))$ for some integers $a\leqslant b$.  Then $L$ lies in  $D^{[a,b]}_{\coh}(A)$.
\end{lemma}
\begin{proof}
	The distinguished triangle $L\xra{\times f} L\to L\otimes_A^LA/(f)\to $ induces a short 	exact sequence 
	\[
	0\to H^i(L)/f H^i(L)\to H^i(L\otimes_A^L A/(f))\to H^{i+1}(L)[f]\to 0
	\]
	for every integer $i\in \ZZ$. It follows that  $H^i(L)/fH^i(L)$ for all $i\in \ZZ$ is a finitely generately $A/(f)$-module, and vanishes if $i\notin [a,b]$.  Note that $H^i(L)$ for $i\in \ZZ$ is derived $f$-complete since so is $L$. It follows from the derived Nakayama \cite[\href{https://stacks.math.columbia.edu/tag/09BA}{Tag 09BA},\href{https://stacks.math.columbia.edu/tag/09B9}{Tag 09B9}]{stacks-project} that $H^i(L)$ is a finitely generated $A$-module, and $H^i(L)=0$ if $i\notin [a,b]$. 
\end{proof}

We have the following finiteness theorem for the Breuil--Kisin cohomology of analytic prismatic crystals. 
\begin{theorem}\label{T:BK-analytic-cohomology}
	Let $X$ be a proper semistable $p$-adic formal scheme over $\Spf(\cO_K)$ of relative dimension $d$. 
	Let $\calE$ be an object of $\Vect^{\an}(\uX_{\prism}, \cO_{\prism})$. Then   $R\Gamma_{\gothS}(X, \calE)$ lies in $D^{[0, 2d]}_{\coh}(\gothS)$.
\end{theorem}

\begin{proof}
	Note that $R\Gamma_{\gothS}(X, \calE)$ is derived $(p,E(u))$-complete. 
	By Lemma~\ref{L:nakayama}, it is sufficient to show that  $R\Gamma_{\gothS}(X,\calE)\otimes^L_{\gothS} \cO_K$ lies in $D^{[0,2d]}_{\coh}(\cO_K)$.
	The short exact sequence of sheaves on $(\uX/\ufS)_{\prism}$
	$$0\to j_*(\calE)\xra{\times E(u)}j_*(\calE)\to j_*(\calE)/E(u)j_*(\calE)\to 0$$  induces a canonical isomorphism 
	\[
	R\Gamma_{\gothS}(X,\calE)\otimes^L_{\gothS} \cO_K\cong R\Gamma((\uX/\ufS)_\prism, j_*(\calE)/E(u)j_*(\calE)).
	\] 
	Let $\nu\colon (\uX/\ufS)_{\prism}\to X_{\et}$ denote the canonical projection of sites. Then we have 
	\[R\Gamma((\uX/\ufS)_\prism, j_*(\calE)/E(u)j_*(\calE))\simeq R\Gamma(X_{\et}, R\nu_*(j_*(\calE)/E(u)j_*(\calE))).
	\]
	By the finiteness theorem of coherent sheaves, it is enough to show that $R^i\nu_*(j_*(\calE)/E(u)j_*(\calE))$ for all $i\in \ZZ$ is a coherent $\cO_X$-sheaf and vanishes if $i\notin [0,d]$. The problem being clearly local for the \'etale topology on $X$,  we may assume that $X=\Spf(R)$ is affine small with framing \eqref{E:small-framing}. Then the desired statement is equivalent to saying that $R\Gamma((\uX/\ufS)_{\prism}, j_*(\calE)/E(u)j_*(\calE))$ lies in $D^{[0,d]}_{\coh}(R)$.

	Put $\overline{j_*\calE}:=\varprojlim_{n} j_*(\calE)/(p^n,E(u)j_*(\calE))$, which is an object of $\CRhat((\uX/\ufS)_{\prism}, \overline \cO_{\prism})$. 
	We have a canonical map $j_*(\calE)/E(u)j_*(\calE)\to \overline{j_*\calE}$. As $j_*(\calE)(\uR_{\gothS})$ is finitely generated over  $R_{\gothS}$ (Remark~\ref{R:torsion-freeness}),  we deduce  from Lemma~\ref{L:completeness-base-change} canonical  isomorphisms of cosimplicial $\cO_K$-modules\footnote{I don't know whether $j_*(\calE)/E(u)j_*(\calE)\simeq \overline{j_*\calE}$. The problem is that for a general object $(B,(E(u)), M_{\Spf(B)})$ of $(\uX/\ufS)_{\prism}$, the $B$-module $j_*(\calE)(B, E(u),M_{\Spf(B)})/E(u)$ may not be $p$-adically complete, even though it is always derived $p$-complete. This pathology does not happen for $\uR_{\gothS}^{\rel, \bullet}$.}:
	\[
	(j_*(\calE)/E(u)j_*(\calE) ) (\uR_{\gothS}^{\rel, \bullet})\simeq j_*(\calE)(\uR_{\gothS}^{\rel, \bullet})\otimes_{\gothS}\cO_K\simeq  \overline{j_*\calE}(\uR_{\gothS}^{\rel, \bullet}).
	\]
	Applying  Lemma~\ref{L:CA-complex-relative} to  both $j_*(\calE)$ and $\overline{j_*\calE}$, we obtain isomorphisms in $D(\cO_K)$:
	\[
	R\Gamma((\uX/\ufS)_{\prism}, j_*(\calE)/E(u)j_*(\calE))\simeq R\Gamma((\uX/\ufS)_{\prism}, \overline{j_*(\calE)})\simeq\CA(\uR^{\rel, \bullet}_{\gothS}, j_*\calE)\otimes_{\gothS}\cO_K.
	\]
	On the other hand, Proposition~\ref{P:Higgs-equiv-BK} shows that 
	\[
	R\Gamma((\uX/\ufS)_{\prism},\overline{j_*(\calE)})\simeq (\overline{j_*\calE}(\uR_{\gothS})\otimes_{R}\widehat{\Omega}^{\bullet}_{R, \log}, \Theta^{\bullet}_{\overline{j_*\calE}}),
	\]
	where every term in the complex on  the right hand side is a finitely generated   $R$-module. In particular, we have $R\Gamma((\uX/\ufS)_{\prism}, j_*(\calE)/E(u)j_*(\calE))\simeq R\Gamma((\uX/\ufS)_{\prism}, \overline{j_*(\calE)})\in D^{[0,d]}_{\coh}(R)$.  
\end{proof}
\begin{remark}
	As $\gothS$ is a noetherian regular local ring of dimension $2$, in the situation of  Theorem~\ref{T:BK-analytic-cohomology},  $R\Gamma_{\gothS}(X, \calE)$ is actually a perfect complex of $\gothS$-modules with tor-amplitude  $[-2, 2d]$. 
\end{remark}

\section{Prismatic-\'etale comparison  theorem for analytic prismatic crystals}\label{Section:comparison}

We fix an algebraic closure $\overline K$ of $K$, and put $C=\widehat{\overline{K}}$. 
Fix a compatible system of $p^n$-th roots $\pi^{1/p^n}\in \cO_C$ of $\pi$. Put $\pi^{\flat}=(\pi, \pi^{1/p}, \dots)\in \cO_C^{\flat}$. 

Recall the log prism $\tAinf(\underline\cO_C)=(A_{\inf}, ([p]_q), M_{\cO_{C}^{\flat}})^a$ defined in \S~\ref{S:setup of notation}. 
Consider the morphism of prelog rings 
\[
\iota\colon (\gothS, \NN\xra{1\mapsto u}\gothS)\to (A_{\inf},  M_{\cO_C^{\flat}})
\]
with $\gothS\to A_{\inf}$
given by $u\mapsto [\pi^{\flat, p}]$, and $ \NN\to M_{\cO^{\flat}_C}=\cO_{C}^{\flat}\backslash\{0\} $ by $1\mapsto \pi^{\flat, p}$. 
Note that $\iota$ is compatible with the natural  $\delta_{\log}$-structures on both sides, and $\iota(E(u))A_{\inf}=([p]_q)=\ker(\theta_{\cO_C}\circ\phi^{-1})$.  Hence, it induces a morphism of log prisms 
\begin{equation}\label{E:gothS-Ainf}
	\ufS=(\gothS, (E(u)), \NN)^a\to \tAinf(\underline{\cO}_C)=(A_{\inf}, ([p]_q), M_{\cO_C^{\flat}})^a.
\end{equation}
Note that the underlying map $\iota\colon \gothS\to A_{\inf}$ is faithfully flat by \cite[Lemma 4.30]{BMS1}.

\subsection{Base change of relative prismatic sites} Let $X$ be a separated semi-stable $p$-adic formal scheme over $\Spf(\cO_K)$.
Put $X_{\cO_C}:=X\times_{\Spf(\cO_K)}\Spf(\cO_C)$. We have the associated log formal scheme $\underline X_{\cO_C}=(X_{\cO_C}, \calM_{X_{\cO_C}})$ with the canonical log structure. 
There exists a commutative diagram of log formal schemes %
\[
\xymatrix{ \uX_{\cO_C}\ar[r]^{\lambda}\ar[d] & \uX\ar[d]\\
	(\Spf(A_{\inf}), M_{\cO_C^{\flat}})^a\ar[r] & (\Spf(\gothS), \NN)^a.
}
\]
Recall that we have an equivalence of sites $\uX_{\cO_C, \prism}\simeq (\uX_{\cO_C}/\tAinf(\underline{\cO}_C))_{\prism}$,  because $\tAinf(\underline{\cO}_C)$ is an initial object of $\Spf(\uO_C)_{\prism}^{\opp}$.
We define a   functor  of  prismatic sites
\[
\lambda_{\sharp}\colon \uX_{\cO_C, \prism}\simeq (\uX_{\cO_C}/\tAinf(\uO_C))_{\prism}\to (\uX/\ufS)_{\prism}
\]
as follows. Let $(\Spf(B), ([p]_q), M_{\Spf(B)}; f)$ be an object of $(\uX_{\cO_C}/\tAinf(\underline{\cO}_C))_{\prism}$, i.e. $(B, ([p]_q), M_{\Spf(B)})$ is a log prism over $\tAinf(\uO_C)$, and  $f\colon (\Spf(\bar B), M_{\Spf(\bar B)})\to \uX_{\cO_C}$  is a strict morphism of log formal schemes, where  $\bar B=B/([p]_q)$ and $M_{\Spf(\bar B)}$ is the striction of $M_{\Spf(B)}$ to $\Spf(\bar B)$. In particular, one has $M_{\Spf(\bar B)}\simeq f^*\calM_{X_{\cO_C}}$.
 Then the canonical surjection  $M_{\Spf(B)}\to i_{B,*}M_{\Spf(\bar B)}$ is a torsor under  $1+[p]_q\cO_{\Spf(B)}$. 
Note that $\lambda^{*}\calM_{X}\subset \calM_{X_{\cO_C}}$ is a submonoid.  Let  $M'_{\Spf(\bar B)}\subset M_{\Spf(\bar B)}$  be the subsheaf of monoids corresponding to $f^*\lambda^*\calM_X\subset f^*\calM_{X_{\cO_C}}$ under the isomorphism, and let $M'_{\Spf(B)}\subset M_{\Spf(B)}$ be the preimage of $M'_{\Spf(\bar B)}$. Then there exists a natural  induced $\delta_{\log}$-structure on $M'_{\Spf(B)}$ so that   $(B, ([p]_q), M'_{\Spf(B)})$ becomes a log prism. By construction of $M'_{\Spf(B)}$, one sees easily that  the composite  morphism $(\Spf(\bar B), M'_{\Spf(\bar B)})\to \uX_{\cO_C}\to \uX$ is strict, and the map of log prisms $\ufS\to \tAinf(\uO_C)\to (B, ([p]_q), M_{\Spf(B)})$ factors though $ (B, ([p]_q), M'_{\Spf(B)})$. Therefore, $(\Spf(B), ([p]_q), M'_{\Spf(B)}; \lambda\circ f)$ is a well defined object of $(\uX/\ufS)_{\prism}$, and we define 
\[
\lambda_{\sharp}(\Spf(B), ([p]_q), M_{\Spf(B)};f)=(\Spf(B), ([p]_q), M'_{\Spf(B)}; \lambda\circ f).
\]
It is obvious that $\lambda_{\sharp}$ is both continuous and cocontinuous, and it induces thus a morphism of topoi 
\[
\lambda\colon \Shv(\uX_{\cO_C, \prism})\simeq \Shv(\uX_{\cO_C}/\tAinf(\uO_C))_{\prism}\to \Shv(\uX/\ufS)_{\prism}
\]
such that we have $\lambda^*\calF=\calF\circ \lambda_{\sharp}$  for any sheaf $\calF\in \Shv(\uX/\ufS)_{\prism}$.
In particular,  $\lambda^*$ induces a pullback functor of complete prismatic crystals:
\[
\lambda^*\colon \CRhat((\uX/\ufS)_{\prism}, \cO_{\prism})\to \CRhat(\uX_{\cO_C, \prism}, \cO_{\prism}).
\]

\subsection{Comparison of local framings}
Assume that $X=\Spf(R)$ is affine small with framing $\square\colon R^{\square}\to R$ as in  \eqref{E:small-framing}. Then the base change to $\cO_C$ induces a framing for $R_{\cO_C}$:
\[
\square_{\cO_C}\colon R^{\square}_{\cO_C}:=\cO_C\langle T_1,\cdots, T_r, T_{r+1}^{\pm 1},\cdots, T_d^{\pm 1}\rangle/(T_0\cdots T_r-\pi)\to R_{\cO_C}.
\] 
Let $\uR_{\gothS}=(R_{\gothS}, (E(u)), \NN^{r+1})^a$ be the Breuil--Kisin log prism for the framing  $\square$ defined in  \S~\ref{S:BK-log-prism}, and $\urA(R)=(\rA(R), ([p]_q), M_{R_{\cO_C}})^a$ be the log prism in the site $\uX_{\cO_C, \prism}$  constructed in \S~\ref{S:frames}. 
There exists a commutative diagram  of log prisms 
\begin{equation}\label{E:diag-local-framings}
	\xymatrix{
		\ufS\ar[r]\ar[d] & \tAinf(\underline{\cO_C})\ar[d]\\
		\uR_{\gothS}\ar[r] & \urA(R),
	}
\end{equation}
where the upper horizontal map is \eqref{E:gothS-Ainf}, the two vertical arrows are natural structural maps, and the lower horizontal one  is defined as follows:
\begin{itemize}
	\item The underlying map of $\delta$-algebras $R_{\gothS}\to \rA(R)$ is induced by 
	\begin{align*}
		R_{\gothS}^{\square}&=\gothS\langle T_0,\dots, T_r, T_{r+1}^{
			\pm1},\dots, T_d^{\pm1}\rangle/(T_0\cdots T_r-u)\\
		&\longrightarrow \rA(R^{\square}_{\cO_C})=A_{\inf}\langle T_0,\dots, T_r, T_{r+1}^{
			\pm1},\dots, T_d^{\pm1}\rangle/(T_0\cdots T_r-[\pi^{\flat}]^p)\end{align*} is given by $u\mapsto[\pi^{\flat}]^p$ and $(T_0,T_1, \dots, T_d)\mapsto (T_0, T_1 \dots, T_d)$.
	
	\item The underlying map of log structures between $\uR_{\ufS}$ and $\urA(R)$ is induced by the natural map of monoids
	\[
	\NN^{r+1}\to M_{R_{\cO_C}}=M_{\cO_{C}^{\flat}}\sqcup_{\NN} \NN^{r+1}.
	\]
		
\end{itemize}
It is easy to see that $\uR_{\gothS}\to \urA(R)$ induces a morphism $\uR_{\gothS}\to \lambda_{\sharp}(\urA(R))$ in $(\uX/\ufS)_{\prism}^{\opp}$. 


Let $\uR_{\gothS}^{\rel, \bullet}=(R^{\rel, \bullet}_{\gothS}, (E(u)), M^{\rel, \bullet}_{R_{\gothS}})$ (resp. $\urA^{\bullet}(R_{\cO_C})=(\rA^{\bullet}(R_{\cO_C}), ([p]_q), M^{\bullet}_{\rA(R)})$) be the \v{C}ech nerve of $\uR_{\gothS}$ (resp. $\urA(R)$) over the final object of $\Shv((\uX/\ufS)_{\prism})$ (resp. of $\Shv(\uX_{\cO_C, \prism})$). 
As $\lambda_{\sharp}$ commutes with products and fiber products,  $\uR_{\gothS}\to \urA(R)$ induces  a morphism of cosimplicial log prisms $\uR_{\gothS}^{\rel, \bullet}\to \lambda_{\sharp}\urA^{\bullet}(R)$ in $(\uX/\ufS)_{\prism}^{\opp}$.

\begin{lemma}\label{L:comparison-framing}
	The natural base change map 
	\[
	R_{\gothS}^{\rel, \bullet}\widehat{\otimes}_{\gothS}A_{\inf}\xra{\sim} \rA^{\bullet}(R_{\cO_C}) 
	\]
	is an isomorphism of  cosimplicial $\delta$-$A_{\inf}$-algebras. Here, $\widehat{\otimes}$ means the classical $(p, [p]_q)$-adic completion of the usual tensor product. 
\end{lemma}

\begin{proof}
	
	The statement follows from explicit constructions of $\uR^{\rel, \bullet}_{\gothS}$ in  and  $\urA^{\bullet}(R)$ (cf. \cite[\S 3.2.1]{MW} and \S~\ref{S:simplicial-A(R)}).
\end{proof}


\begin{proposition}\label{P:base-change-coh-O_C}
	Let $X$ be a separated semi-stable $p$-adic formal scheme over $\Spf(\cO_K)$, and $\calE$ be an object of $\Vect^{\an}(\uX_{\prism}, \cO_{\prism})$. Then the natural base change map 
	\[
	(R\Gamma_{\gothS}(X, \calE)\otimes_{\gothS}^L A_{\inf})^{\wedge}\xra{\sim} R\Gamma(\uX_{\cO_C, \prism}, \lambda^*(j_*\calE))
	\]
	is an isomorphism, where $(-)^{\wedge}$ is the derived completeion with respect to  the ideal $(p, [p]_q)$. Moreover, if $X$ is proper over $\Spf(\cO_K)$, the isomorphism holds without the derived completion $(-)^{\wedge}$. 
	
\end{proposition}

\begin{proof}
	By \'etale descent, we reduce to the case when $X=\Spf(R)$ is affine small semi-stable with framing \eqref{E:small-framing}. Then by Lemma~\ref{L:CA-complex-relative},  $R\Gamma_{\gothS}(X, \calE)$ is represented by the \v{C}ech--Alexander complex  $\CA(\uR^{\bullet}_{\gothS}, j_*\calE)$, and similarly $R\Gamma(\uX_{\cO_C, \prism}, \lambda^*(j_*\calE))$  is represented by $\CA(\urA^{\bullet}(R), \lambda^*(j_*\calE))$ by Lemmas~\ref{L:CA-complex-prismatic}, \ref{L:cover-final-object}.

	Recall that $(j_*\calE)(\uR_{\gothS})$ is finitely generatd over $R_{\gothS}$ (Remark~\ref{R:torsion-freeness}). Let $p^n_0\colon \uR_{\gothS}= \uR^{\rel, 0}_{\gothS}\to \uR^{\rel, n}_{\gothS}$ be the canonical map corresponding to the inclusion of simplices $[0]\hra [n]$. The underlying map of $\gothS$-algebras $p_n^0: R_{\gothS}\to R^{\rel, n}_{\gothS}$ is $(p, E(u))$-faithfully flat (\cite[Lemma~3.16]{MW}) hence flat by \cite[\href{https://stacks.math.columbia.edu/tag/0912}{Tag 0912}]{stacks-project}. 
	Then  the crystal property of $j_*(\calE)$ and  \cite[\href{https://stacks.math.columbia.edu/tag/0EEU}{Tag 0EEU}]{stacks-project} imply that the $(p, [p]_q)$-derived completion of $(j_*\calE)(\uR_{\gothS}^{\rel, n})\otimes_{\gothS}A_{\inf}\simeq (j_*\calE)(\uR_{\gothS})\otimes_{R_{\gothS}, p_0^n}R^{\rel, n}_{\gothS}\otimes_{\gothS}A_{\inf}$ coincides with its classical $(p, [p]_q)$-adic completion $(j_*\calE)(\uR_{\gothS}^{\rel, n})\widehat\otimes_{\gothS}A_{\inf}$. On the other hand,  by Lemma~\ref{L:comparison-framing} and the crystal property of $j_*\calE$, we have a canonical isomorphism of cosimplicial $A_{\inf}$-modules
	\[
	\lambda^*(j_*\calE)(\urA^{\bullet}(R))=(j_*\calE)(\lambda_{\sharp}\urA^{\bullet}(R))\simeq (j_*\calE)(\uR^{\rel,\bullet}_{R})\widehat\otimes_{\gothS}A_{\inf}.
	\]
	It follows therefore  that  
	\begin{align*}
		(R\Gamma_{\gothS}(X, \calE)\otimes_{\gothS}^L A_{\inf})^{\wedge}&\cong \CA(\uR_{\gothS}^{\rel, \bullet}, j_*\calE)\widehat{\otimes}_{\gothS}A_{\inf}\\
		&\cong\CA(\urA^{\bullet}(R), \lambda^*(j_*\calE))\\
		&\cong R\Gamma(\uX_{\cO_C,\prism}, \lambda^*(j_*\calE)).
	\end{align*}
	This proves the first part of the Proposition. The second part follows from \cite[\href{https://stacks.math.columbia.edu/tag/0EEV}{Tag 0EEV}]{stacks-project} and the fact that  $R\Gamma_{\gothS}(X,\calE)$  is pseudo-coherent (Theorem~\ref{T:BK-analytic-cohomology}).
\end{proof}

\begin{remark}\label{R:comp-Frob-action}
	In the situation of Proposition~\ref{P:base-change-coh-O_C}, assume moreover that $\calE$ is equipped with an Frobenius structure $\phi_{\calE}\colon (\phi_{\cO_{\prism}}^*\calE)[\frac{1}{\calI_{\prism}}]\simeq \calE[\frac{1}{\calI_{\prism}}]$. Then there is  a natural induced Frobenius structure on $(j_*\calE)|_{\uX_{\cO_C}}$, and the isomorphism given by  \ref{P:base-change-coh-O_C} is equivariant under the natural Frobenius actions on both sides: 
	\[
	\xymatrix{\big(\phi_{\gothS}^*R\Gamma_{\gothS}(X,\calE)\otimes^L_{\gothS}A_{\inf}\big)^{\wedge}[\frac{1}{[p]_q}]\ar[r]\ar[d]^{\simeq} &\big(R\Gamma_{\gothS}(X, \calE)\otimes^L_{\gothS}A_{\inf}\big)^\wedge[\frac{1}{[p]_q}]\ar[d]^{\simeq}\\
		\big(\phi_{A_{\inf}}^*	R\Gamma(\uX_{\cO_C, \prism}, \lambda^*(j_*\calE))\big)[\frac{1}{[p]_q}]\ar[r] & R\Gamma(\uX_{\cO_C, \prism}, \lambda^*(j_*\calE))[\frac{1}{[p]_q}]
	}
	\]
	where the horizontal maps are induced by the Frobenius structure $\phi_{\calE}$ and vertical arrows are isomorphisms given by Proposition~\ref{P:base-change-coh-O_C}.
\end{remark}

\subsection{\'Etale realization of analytic prismatic crystals} 
Recall the category of  prismatic Laurent $F$-crystals $\Vect(\uX_{\prism}, \cO_{\prism}[\frac{1}{\calI_{\prism}}])^{\phi=1}$ (Definition~\ref{D:prismatic-Laurent-crystal}) and $\Loc_{\ZZ_p}(X_{\eta,\et})$, the category  of \'etale 
$\ZZ_p$-local systems on the generic adic fiber $X_{\eta}$. 
Consider  the composite functor
\[
T\colon \Vect^{\an}(\uX_{\prism}, \cO_{\prism})^{\phi=1}\to \Vect(\uX_{\prism},\laurent)^{\phi=1}\simeq \Loc_{\ZZ_p}(X_{\eta,\et}),
\]
called \emph{ \'etale realization functor}, where the first functor is the natural restriction, and  the second equivalence is Theorem~\ref{T:etale-realization}. By \cite[Prop. 3.21]{DLMS},  $T$ is fully faithful.
Let $\Loc^{\st}_{\ZZ_p}(X_{\eta, \et})$ denote the essential image of the functor $T$. By \cite[Cor. 5.2]{DLMS}, $\Loc^{\st}_{\ZZ_p}(X_{\eta, \et})$ consists exactly of semi-stable  \'etale  $\ZZ_p$-local systems on $X_{\eta, \et}$  in the  sense that they are attached to certain crystalline $F$-crystals on  the absolute log-crystalline site  of $\uX_1=\uX\times_{\Spf(\cO_K)}\Spf(\cO_K/p)$ (See \cite[Def. 3.40]{DLMS}).  When $X=\Spf(\cO_K)$, $\Loc^{\st}_{\ZZ_p}(X_{\eta, \et})$ is exactly the category of  Galois stable $\ZZ_p$-lattices in semi-stable representations of $\Gal(\overline K/K)$. 

Now we can put everything together to prove Theorem~\ref{T:comp-prismatic-etale-analytic}.

\begin{proof}[End of Proof of \ref{T:comp-prismatic-etale-analytic}]

It is clear that, for an object $(\calE, \phi_{\calE})$ of $\Vect^{\an}(\uX_{\prism}, \cO_{\prism})^{\phi=1}$, the restriction to $X_{C, \et}$ of $T(\calE)$ coincides with the \'etale realization of the object $(\lambda^*j_*(\calE), \phi_{\lambda^*j_*(\calE)})$ via the functor \eqref{E:etale-realization-complete-F-crystal}. Then Theorem~\ref{T:comp-prismatic-etale-analytic} follows immediately from Proposition~\ref{P:base-change-coh-O_C} and  Theorem~\ref{T:etale-comparison-O_C}. The last part of \ref{T:comp-prismatic-etale-analytic}(2) follows from the isomorphism 
\[
R\Gamma_{\gothS}(X, \calE)\otimes_{\gothS}A_{\inf}[\frac{1}{\mu}]\simeq R\Gamma(X_{C,\et}, T(\calE))\otimes A_{\inf}[\frac{1}{\mu}]
\]
and the fact that $A_{\inf}$ is flat over $\ZZ_p$ and $\gothS$. 

\end{proof}

\begin{remark}
 To each object $(\calE, \phi_{\calE})$ of $\Vect^{\an}(\uX_{\prism}, \cO_{\prism})^{\phi=1}$, Du--Liu--Moon--Shimizu  associate in \cite[Construction 3.39]{DLMS} a crystalline $F$-isocrystal $\calE^{\cris}$ on the absolute log crystalline site $\uX_{1, \cris}$ of $\uX_1= \uX\otimes_{\cO_K} \cO_K/p$.  It is also expected that there exists a prismatic-crystalline comparison theorem to relate  $R\Gamma_{\gothS}(X, \calE)$  to the cohomology of $\calE^{\cris}$.  Combining this conjectural comparison theorem with Theorem~\ref{T:comp-prismatic-etale-analytic}, one would get a comparison crystalline-\'etale comparison theorem for semi-stable local systems that generalizes Tsuji's proof of $C_{st}$-conjecture \cite{Tsuji1}.
  \end{remark}



\end{document}